\documentclass[9pt]{amsart}
\usepackage{amsmath}
\usepackage{amssymb}
\usepackage{amsfonts}
\usepackage{mathrsfs}
\usepackage{color}
\usepackage{graphics}
\usepackage[dvips]{graphicx}
\usepackage{wrapfig}
\usepackage{upgreek}
\usepackage[colorinlistoftodos, dvistyle, textsize=tiny]{todonotes} %for fancy comments
\newtheorem{theorem}{Theorem}[section]
\newtheorem{lemma}{Lemma}[section]
\newtheorem{proposition}{Proposition}[section]
\newtheorem{corollary}{Corollary}[section]
\newtheorem{definition}{Definition}[section]
\newtheorem{remark}[theorem]{Remark}

\newcommand{\p}{\partial}
\newcommand\R{\mathbb R}

\newcommand\Wu{\mathbb W_H}
\newcommand{\eqdef}{\overset{def}{=}}

\newcommand\Ga{\Gamma}

\newcommand\dist{\text{dist}}
\renewcommand\det{\text{det}}

\newcommand\ee{\eps}
\newcommand\eps{\varepsilon}
\newcommand\ol{\overline}

\newcommand{\e}{\ee}
\renewcommand\phi{\varphi}

\newcommand\cal{\mathcal}

\newcommand\W{\mathbb W}
\newcommand{\Wup}{\overline \W_{\rm H}}
\newcommand{\Wlo}{\underline \W_{\rm E}}

\newcommand\pr{{\bf Proof.\ }}
\newcommand{\U}{\mathcal U}
\newcommand{\V}{\mathcal V}
\newcommand{\s}{\mathscr{S}}
\renewcommand{\r}{\mathscr{ R}}

\renewcommand\phi\varphi % use curly phi
\newcommand\pz{\mathpzc}                                                                                                             %
\renewcommand\a{\pz A}
\renewcommand\b{\pz B}
\numberwithin{equation}{section}
\newcommand\kap{\kappa}% the deviation constant 1-1/eccen^2
\newcommand{\Id}{{\rm{Id}}}

\numberwithin{equation}{section}

\DeclareMathAlphabet{\mathpzc}{OT1}{pzc}{m}{it}
\newtheorem*{AAA}{Theorem A}                                                                   %
\newtheorem*{BBB}{Theorem B}
\newtheorem*{CCC}{Theorem C}
\newtheorem*{DDD}{Theorem D}

\def\H{\mathcal H^n}
\usepackage{picins,graphicx}
\usepackage{sidecap}
\usepackage{subfig}

\begin{document}

\title[Surfaces refracting parallel lighting]{An inverse problem for the refractive surfaces with parallel lighting
%and Existence and Regularity for the surfaces refracting co-axial light rays
}
\author{Aram L. Karakhanyan}
\address{Maxwell Institute for Mathematical Sciences and School of Mathematics,
University of Edinburgh, James Clerk Maxwell Building, King's Buildings, Mayfield Road, EH9 3JZ,  Edinburgh, Scotland }
\email{aram.karakhanyan@ed.ac.uk}
\thanks{$2010$ {\it Mathematics Subject Classification.\/} Primary
 35J96, Secondary  78A05, 78A46, 78A50.}
\maketitle

%\numbering
%\footnote{AMS Classifications: }

\begin{abstract}
In this article we examine the regularity of two types of weak
solutions to a Monge-Amp\`ere type equation which emerges in
a problem of finding surfaces that refract parallel light rays emitted from source domain and striking a  given
target after refraction. Historically, ellipsoids and hyperboloids of revolution were the first surfaces to be considered in this context.
The mathematical formulation commences with deriving the energy conservation equation
for sufficiently smooth surfaces, regarded as graphs of functions to be sought, and then studying the existence and regularity
of two classes of suitable weak solutions constructed from  envelopes of hyperboloids or ellipsoids of revolution. Our main result in this article states that under suitable conditions on source  and target domains and respective intensities these weak solutions are locally smooth.
\end{abstract}

\setcounter{tocdepth}{1}
\tableofcontents

\section{Introduction}

%%%%%%%%%%%%%%%%%%%%%%%%%%%%%%%%%%%%%%%%%%%
%%%%%%%%%%%%%%%%%%%%%%%%%%%%%%%%%%%%%%%%%%%

Let $\U\subset\R^n$ be a bounded domain with smooth boundary
and $u:\U\to \R$  a smooth function. By $\Gamma_u$ we denote the
graph of $u$. Let $\gamma$  denote the  unit normal of $\Gamma_u$.
We think of $\Ga_u$ as  a surface that dissevers two distinct media.

From each $x\in \U$ we issue a ray $\ell_x$
parallel to $e_{n+1}-$the unit direction of the
$x_{n+1}$ axis in $\R^{n+1}$. Then $\ell_x$ strikes  $\Gamma_u$,
the surface separating the two media I and II, refracts into the second media II and strikes the
receiver surface $\Sigma$, see Figure 1. Let $Y$ be the unit direction of the refracted ray.
%%%%%%%%%%%%%%%%%%%%%%%%%
\begin{figure}
 \includegraphics[scale=0.5]{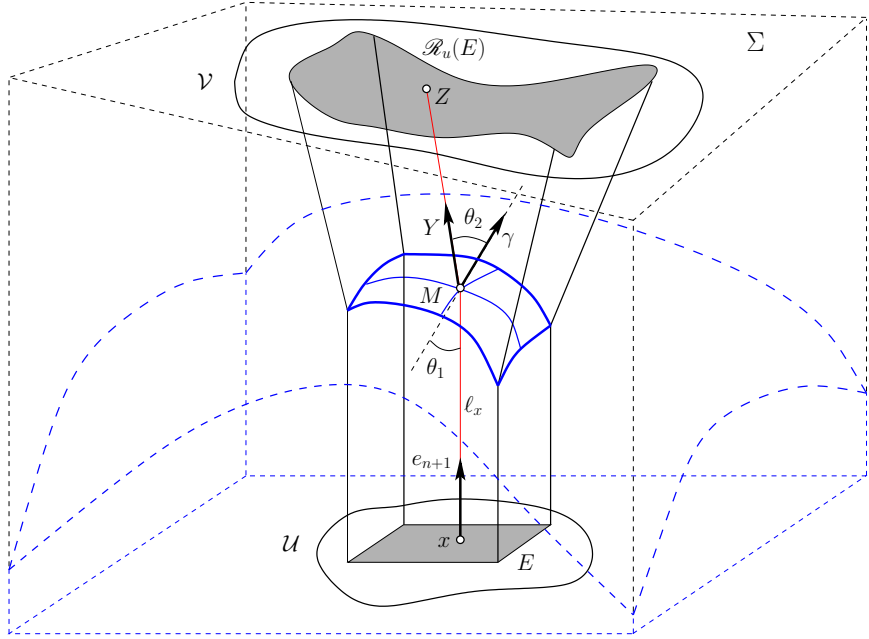}
 \caption{The blue doted lines confine the boundary of media $\rm I$.}
\end{figure}

%%%%%%%%%%%%%%%%%%%%%%%%%
If $\gamma$ is the unit normal at $M=(x,u(x))\in \R^{n+1}$  where  $\ell_x$ strikes $\Gamma_u$
then from the refraction law we have
\begin{eqnarray}\label{angles}
  \frac{\sin\theta_1}{\sin\theta_2}=\frac{n_2}{n_1},
\end{eqnarray}
where $n_1, n_2$ are the refractive indices of  the  media I and II respectively,
dissevered by the interface $\Gamma_u$, $\theta_1$ and $\theta_2$ are  the angles
between $\ell_x$ and $\gamma$, and between $Y$ and $\gamma$, respectively, see Figure 1.

\smallskip
Suppose that the intensity of light on $\U$ is $f\geq 0$ and let $\V$ be the set of points
where the refracted rays strike the receiver $\Sigma$. Denote by $g\geq 0$
the gain intensity on $\V$. For each $\U'\subset \U$ let $\V'$ be the set of points where
the rays, issued from $\U'$ and refracted off $\Gamma_u$,
 strike $\Sigma$. Thus $u$ generates the refractor mapping
 $$Z_u: \U\longrightarrow \V$$
and the illuminated domain on $\Sigma$,
corresponding to $\U'\subset \U$,  is $\V'=Z_u(\U')$. If $\Gamma_u$ is a perfect refractor, then one would have the  energy balance equation (in local form)
\begin{equation}\label{eq-energ-bal}
 \int\limits_{\U'}f=\int\limits_{\V'=Z_u(\U')}g.
\end{equation}

\smallskip

The main problem that we are concerned with  is formulated below:

\subsection*{Problem}{\it Assume that we are given a smooth surface $\Sigma$ in $\R^{n+1}$, a pair of bounded smooth
 domains $\U\subset \Pi=\{X\in\R^{n+1} : X^{n+1}=0\}$ and $\V\subset \Sigma$ and
 a pair of nonnegative, integrable functions $f:\U\to\R$ and $g:\V\to \R$ such that the energy balance condition holds
 \begin{equation}\label{glob-balance}
 \int_\U f=\int_\V gd\H.
\end{equation}

 Find a
 function $u:\U\to\R$ such that the following two conditions are fulfilled }

\begin{align}\label{problem}
\tag{\textit{\bf RP}}
\left\{
\begin{array}{lll}
 \int\limits_{\U'}f=\int\limits_{Z_u(\U')}g,\ \mbox{\it for\ any\ measurable}\ \U'\subset \U \\ %\qquad \mbox{and}\qquad
  %\displaystyle
  Z_u(\U)=\V.
 \end{array}
 \right.
\end{align}

\medskip

Problems of this kind appear in geometric optics \cite{Appl} page 315.
In the 17th century
Descartes posed a similar problem with target set $\V$ being a single point, say $\V=\{Z_0\}$. 
It was observed that the ellipsoids and hyperboloids of revolution with focal axis  parallel to $e_{n+1}$ will solve this problem if $Z_0$ is one of the foci. The case of general target $\V$
can be treated via approximation argument, namely by constructing a solution from ellipsoids or hyperboloids for finite set $\V=\{Z_1, \dots, Z_m\}$ and then letting $m\rightarrow \infty$.
Moreover, the eccentricity of these surfaces is fixed and determined by the refractive indices $n_1$ and $n_2$.
To see this we take advantage of some well-known facts from geometric optics and
record them here for further reference, see \cite{Mountford}.
Let $H(x)=Z^{n+1}-a\eps-\frac ab\sqrt{b^2+|x-x_0|^2}$ be the lower sheet of a
hyperboloid of revolution with focal axis passing through the point $x_0\in \U$ and parallel to
$e_{n+1}$, see Section \ref{sec-ell-hyp}.
Similarly, we define the lower half of an ellipsoid of revolution $E(x)=Z^{n+1}-a\eps-\frac ab\sqrt{b^2-|x-x_0|^2}$.
If $n_1$ and $n_2$ are the refractive indices of media I and II respectively then
\begin{eqnarray}\label{eq-ref-law}
 \eps=\frac{n_1}{n_2}=\frac{\sin\theta_2}{\sin\theta_1}=
\left\{\begin{array}{lll}
        \frac{\sqrt{a^2-b^2}}{a}<1\ \  \textrm{for ellipsoids},\\
        \frac{\sqrt{a^2+b^2}}{a}>1\ \  \textrm{for hyperboloids}.
       \end{array}
\right.
\end{eqnarray}
Here $\varepsilon$ is the eccentricity, see \cite{Mountford}. Since $\ee$ is fixed we can
drop the dependence of $E$ and $H$ from $b=a\sqrt{|\eps^2-1|}$ and take
\begin{eqnarray}
E(x, a, Z)=Z^{n+1}-a\eps -a\sqrt{1-\frac{(x-z)^2}{a^2(1-\eps^2)}}, &\quad \rm{if}\ \eps<1,\\\label{hypo-intro}
 H(x, a, Z)=Z^{n+1}-a\eps -a \sqrt{1+\frac{(x-z)^2}{a^2(\eps^2-1)}}, &\quad \rm{if}\ \eps>1.
\end{eqnarray}
We also define
the constant
\begin{equation}\label{def-kap}
 \kap=\frac{\eps^2-1}{\eps^2}
\end{equation}
which will prove to be useful, in a number of computation to follow.

%%%%%%%%%%%%%%%%%%%%%%%%%%%%%%%%%%%%%%%%%%%%%%%%%%%%%%%%%%%%
%                                                          %
%                SECTION                                   %
%                                                          %
%%%%%%%%%%%%%%%%%%%%%%%%%%%%%%%%%%%%%%%%%%%%%%%%%%%%%%%%%%%%
\smallskip
\section{Main theorems}
Let $\Sigma$ be the receiver surface defined implicitly
\begin{equation}\label{sgma-intro}
\Sigma=\{Z\in \R^{n+1}: \psi (Z)=0\}
\end{equation}
where $\psi: \R^{n+1}\to \R$ is a smooth function.
If $u\in C^2(\U)$ then the first condition in (\ref{problem}), after using change of variables, results a Monge-Amp\`ere type equation
for $u$, whereas the second one plays the role of boundary condition for $u$. More precisely we have the following

\begin{AAA}
 Let $u\in C^2(\U)$ be a solution to (\ref{problem}). Then
 \begin{itemize}
\item[$\bf 1^\circ$] $Y=\eps\left(\frac{\kap Du}{1+q}, 1-\frac{\kap}{1+q}\right)$ is the unit direction of refracted ray,

\item[$\bf 2^\circ$] $u$ solves the equation
\begin{equation}\label{eq-main-th}
\left| \det\left[\frac{q+1}{t\eps\kap}\left\{\Id-\kap\eps^2Du\otimes Du\right\}+D^2u\right]\right|=\left|-\eps q \left[\frac{q+1}{t\eps\kap}\right]^n
\frac{\nabla \psi\cdot Y }{|\nabla \psi|}\frac fg\right|,
\end{equation}
 \end{itemize}
 where
\begin{equation}\label{def-q}
 q(x)=\sqrt{1-\kap(1+|Du|^2)}, \quad\kap=\frac{\eps^2-1}{\eps^2}
\end{equation}
and $t$ is the stretch function
defined in (\ref{def-Z}) via an implicit relation $\psi(x+e_{n+1}u(x)+Yt)=0$.
\end{AAA}
If the receiver $\Sigma$ is a plane then taking $\psi(Z)=Z\cdot \xi+\xi_1$ we find that 
$t=-[Y\cdot \xi_0]^{-1}(x+u(x)e_{n+1}+\xi_1)$. In particular for 
the horizontal plane $X^{n+1}=m$, with some constant $m>0$,  one has
$$t=\frac{m-u}{Y^{n+1}}=(m-u)\frac{q+1}{\eps(1-\kap+q)}.$$
Quadric $\Sigma$ is another example of receiver for which $t$ can be computed explicitly.
In general $t$ is a function of $x, u(x)$ and $Du(x)$ which may not have simple explicit form. 
However, in terms of applications  the case of planar receiver is of particular interest, since the 
flat screens are  easy to construct. 
The method of the stretch function was introduced in
\cite{KW-1, KW-2} to treat the near-field reflection problem.
The equation for a near-field \textit{refraction} 
problem with point source is derived in  \cite{GutH}, \cite{K-AJM}.

\medskip

Next, we need to introduce the notion of weak solution of
(\ref{eq-main-th}). It will allow us to develop the existence theory along the lines of the classical Monge-Amp\`ere equation.  To this end, we say that $u:\U\to \R$ is upper (resp. lower) admissible with respect to $\V$ if for any
$x\in \U$ there is a hyperboloid $H(\cdot, a, Z)$ (resp. ellipsoid $E(\cdot, a, Z)$) with focus
$Z\in \V$ such that $H(\cdot, a, Z)$ (resp. $E(\cdot, a, Z)$) touches $u$ from above (resp. below) at $x$.
Such $H(\cdot, a, Z)$ (resp. $E(\cdot, a, Z)$) is called supporting hyperboloid (resp. ellipsoid) of
$u$ at $x$. To fix the ideas we consider the class of upper admissible function and denote it by $\Wup(\U, \V)$. The class of lower admissible functions is denoted by $\Wlo(\U, \V)$.
For each $u\in \Wup (\U,\V)$ we define the mapping
$\s_u: \V\to \U$  by
\begin{equation*}
 \s_u(Z)=\{x\in \U : \exists a>0\ \mbox{such \ that}\ H(\cdot,a, Z )\ \mbox{is\ a \ supporting \ hyperboloid\ of} \ u\  \mbox{at} \ x\},
 \end{equation*}
 and take 
 \begin{equation*}
  \beta_{u, f}(E)=\int_{\s_u(E)}f(x)dx, \quad E\subset\V. 
\end{equation*}
%%%
Furthermore,  we also consider the mapping 
$\r_u:\U\to \V$ defined by 
\begin{equation*}
 \r_u(x)=\{Z\in \V : \mbox{there \ is\ a \ supporting \ hyperboloid}\  H(\cdot,a, Z ) \ \mbox{ of} \ u\  \mbox{at} \ x\}
\end{equation*}
and associate the following set function 
\begin{equation*}
\alpha_{u, g}(E)=\int_{\r_u(E)}gd\H, \quad E\subset \V.
\end{equation*}
Notice that for smooth $u$, the mapping  $\s_u$ is the inverse of $\r_u$.

With the aid of these set functions $\alpha_{u,g}$ and $\beta_{u,f}$ we can introduce two notions of
weak solution to (\ref{problem}), called $A$ and $B$ type weak solutions, respectively. 
It is not hard to see that $\beta_{u,f}$ is in fact $\sigma$-additive
measure, while for $\alpha_{u,g}$ it is less obvious.   Towards proving this the major obstruction   
is to show that $\r_u$ is one-to-one modulo a set of vanishing $\H$ measure on $\Sigma$. This is 
circumvented by introducing the Legendre-like transformation $v(z)$ of an admissible function $u(x)$ in Section \ref{sec-A-type} defined as an upper envelope of some function of $\dist(Z, X)$ for $Z\in \V$ and $X\in \U$. 
In order to infer that $v(z)$ is  semi-concave  (which in turn will lead to $\sigma$-additivity of $\alpha_{u,g}$) 
we  assume that (\ref{dist-cond})  is fulfilled. That done, one can show that an $A$-type weak solution exists in the 
sense of Definition \ref{bndr-A-def}.  Note that once we found the Legendre-like transformation then 
our problem can be treated as a prescribed Jacobian type equation discussed in \cite{Trud-jaco}. 
However one still has to check all conditions formulated there in order to trigger the theory.
Furthermore, the construction of locally smooth solution for (\ref{problem}) is very complicated and require a 
careful analysis of Dirichlet's problem. This issues are addressed in Lemma \ref{lem-apprx} and Section \ref{sec-Dirichlet}. 

\smallskip 

If, for a moment, we take the existence of $A$-type weak solution for granted,
the question about its regularity is even more complex. 
To set stage for the weak solutions  we assume that
$\Sigma=\{Z\in \R^{n+1}: \psi (Z)=0\}$ and $\psi: \R^{n+1}\to \R$ being  a smooth function.
Clearly, some conditions must be imposed on $\psi$ to guarantee, among other things, that the
right hand side of the equation (\ref{eq-main-th}) is well defined, at least  for smooth solutions.

To this end we enlist the following conditions to be used in the construction of weak solutions and proving their smoothness. 
\begin{eqnarray}\label{vis-cond}
  &&\nabla \psi(Z)\cdot(X-Z)>0 \quad  \forall X\in(\U\times[0,m_0]),\forall Z\in\Sigma\ {\rm and\ for\ some\ large\ constant}\ m_0>0,\\\label{dist-cond}
 &&\dist(\U, \V)>0,\\\label{Rconv-cond}
&& \V\ {\rm is} \ R-{\rm convex\ with \ respect \ to } \ \U, \ {\rm see \ Definition} \ \ref{def-R-conv}, \\\label{dens-cond}
 && f, g>0,\\\label{A3-cond00}
&& \frac1t\left[\frac{t\eps\kap}{q+1}\right]^2 {\rm II}+\frac{\kap}q\frac{\psi_{n+1}
}{|\nabla \psi|} \left(\Id+\kap \frac{p\otimes p}{q^2}\right)<0,\quad {\rm if}\  \kap>0,
\end{eqnarray}
where $\rm II$ is the second fundamental form of $\Sigma$ and $p=Du(x)$. The subdomain of $\U \times[0, \infty)$ where 
(\ref{vis-cond})-(\ref{A3-cond00}) are simultaneously satisfied  is called the {\it regularity domain} $\mathcal D$. 

It is worthwhile to explain the meaning of these conditions: the first one 
(\ref{vis-cond}) means that the reflected rays do not 
strike $\Sigma$ tangentially, otherwise $\Sigma$ would not detect the 
gain intensity at the tangential points, i.e. at the points where   $\nabla \psi(Z)\cdot(X-Z)=0$.
On the technical level, however, it 
allows to apply the inverse function theorem to recover the stretch function $t=t(x, u, Du)$. 
It is worth pointing out that (\ref{vis-cond}) holds for a large class of surfaces $\Sigma$.
To see this it is enough to notice that there is a positive constant $c(\eps)$, depending only on 
$\eps$ such that $Y^{n+1}\in[c(\eps), 1]$. In other words the unit directions $Y$ of 
refracted rays remain within the cone $c(\eps)\leq Y^{n+1}\leq 1$.
%\begin{remark}
Indeed,if $u$ is differentiable at $x$ then  $Y^{n+1}(x)=\eps [1+(\sqrt{\cos^2\theta_1-\kap}-\cos\theta_1)\gamma^{n+1}]$ from 
refraction law, see (\ref{eq-Y}) and Figure 1. Here $\kap =1-\eps^{-2}$. If $u$ is not differentiable at $x$, we interpret  
$\gamma$ as one of 
the normals of supporting  planes of admissible $u$ at $x\in \U$ since $u$ is concave (resp. convex) 
if $u$ is upper (resp. lower) admissible, see Section \ref{sec-ell-hyp}. Thus if $p\in \p u(x)$, where $\p u$ is the subdifferential of 
$u$ at $x$, then $\gamma=\left(\frac{p}{\sqrt{1+|p|^2}}, \frac1{\sqrt{1+|p|^2}}\right)$ and $\cos\theta_1=\frac1{\sqrt{1+|p|^2}}$. Consequently if 
$u$ is lower admissible then $Y^{n+1}\geq \e$ if $\kap<0$, i.e. $\eps<1$ and hence $c(\eps)=\eps<1$.
On the other hand if $\kap>0$ then for any $u\in \mathbb H(\U, \V)$ we have 
\begin{equation}\label{awkw-cond}
\sup_\U|Du|<\frac{1}{\sqrt{\eps^2-1}}\quad {\rm if}\ \kap=\frac{\eps^2-1}{\eps^2}>0.
\end{equation}
This simply follows from the fact that supporting hyperboloids control  the 
magnitude of the gradient of $u$, see Lemma \ref{lem-lips}. But in its turn $|DH|$, for any hyperboloid $H$ given by  (\ref{hypo-intro}),
satisfies the estimate (\ref{awkw-cond}). Because $\gamma^{n+1}=\cos\theta_1=\frac1{\sqrt{1+|p|^2}}$ (see Figure 1 and the derivation of (\ref{unit-Y})) and $|p|< \frac1{\sqrt{\eps^2-1}}$
we infer that  
$$(\sqrt{\cos^2\theta_1-\kap}-\cos\theta_1)\gamma^{n+1}=\frac{-\kap}{\sqrt{\cos^2\theta_1-\kap}+\cos\theta_1}\gamma^{n+1}>-\kap$$
and consequently $Y^{n+1}>\eps(1-\kap)=\eps^{-1}<1$. Thus for $\eps>1$ we can take $c(\eps)=\eps^{-1}$.
From here we see that (\ref{vis-cond}) holds for any horizontal receiver 
$Z^{n+1}=m$, for large $m>0$. More generally if $\Sigma$  is concave in $Z^{n+1}$ direction 
and the normal mapping of $\Sigma$ is strictly inside of the cone $c(\eps)<Y^{n+1}$ on the unit sphere then
(\ref{vis-cond}) holds true. This leads to the following cone condition for the unit directions of refracted rays
\begin{equation}\label{sml-cone}
0<c(\eps)\leq Y^{n+1}\leq 1.
\end{equation}
%\end{remark} 

The second condition (\ref{dist-cond}) assures that the 
Legendre-like transformation $v(z)$ for an admissible function $u$ is well defined as an envelope of 
$C^1$ smooth functions, in particular $v(z)$ is semi-concave and hence differentiable
almost everywhere, see Section \ref{sec-A-type}. This yields that $\alpha_{u, g}$ is a Radon measure.

The next two conditions  (\ref{Rconv-cond}) and (\ref{dens-cond}) assure that 
$B$-type solution is also of $A$-type and therefore one gets the existence of $A$-type weak solutions 
in some indirect way using the methods of \cite{C92}, \cite{Urbas}. That done, we can approximate 
$\V$ by R-convex domains and show the existence of $A$-type weak solutions without assuming (\ref{Rconv-cond}), see Theorem C4.

Last  condition (\ref{A3-cond00}), which is crucial for regularity of weak solutions, 
deserves special attention because it is the most sophisticated one.
In fact the next theorem is entirely devoted to the verification of (\ref{A3-cond00}).

\begin{BBB}
 Let $u$ be a $C^2(\U)$ 
 solution of (\ref{eq-main-th}) and $\mbox{II}$ be the second fundamental form of $\Sigma=\{Z\in \R^{n+1} : Z^{n+1}=\phi(z)\}$. 
 Let  $\kap>0, {\rm II}< 0$  then there is a constant $\lambda >0$ such that if $\dist(\U, \V)>\lambda$ then (\ref{A3-cond00}). 
 Furthermore, for  $\kap<0, {\rm II}>0$  there is a constant $\hat \lambda >0$ such that if $\dist(\U, \V)>\hat \lambda$ then the opposite inequality in 
 (\ref{A3-cond00}) holds for all $x\in \U$.
  \end{BBB}

%%%%%%%%%%%%%%%%%%%%%%%%%%%%%%%%%%%%%%
%%%%%%%%%%%%%%%%%%%%%%%%%%%%%%%%%%%%%%
%%%%%%%%%%%%%%%%%%%%%%%%%%%%%%%%%%%%%%

\smallskip
If $\Sigma$ is a graph, say $Z^{n+1}=\phi(z)$ then (\ref{A3-cond00}) can be rewritten as  
$$
\frac1t\left[\frac{t\eps\kap}{q+1}\right]^2 D^2\phi+\frac{
\kap}q \left(\Id+\kap \frac{p\otimes p}{q^2}\right)<0,\quad {\rm if}\  \kap>0.
$$
In lieu of (\ref{sml-cone}) this assumption on $\Sigma$ is not restrictive. In addition, Theorem B suggests that it is convenient to think of $\Sigma$ as an unbounded  convex (resp. concave) surface without boundary if $\kap <0$ (resp. $\kap>0$) by extending 
$\phi $ to $\R^n$ as a convex function $\widetilde \phi$ such that $\phi(z)\to \pm\infty$ as $|z|\to\infty$. We will take advantage of such extension of $\phi$ (and hence $\Sigma$) in Section \ref{ssec-confocal} and Lemma \ref{lem-apprx}, see also Remark \ref{rem-extns}.

\smallskip
%%%%% Theorem C
Now we are ready to formulate our main existence result.
\begin{CCC}
\begin{itemize}
\item[1] If $f, g\ge 0$ and   (\ref{glob-balance}) and \eqref{vis-cond} hold then there is a $B$-type weak solution provided that the condition below
\begin{equation}\label{vis-cond-B}
Z^{n+1}\geq \left[\frac2{\eps-1}+\frac1{\sqrt{\eps^2-1}}\right]\rho(z)
\end{equation}
is satisfied. Here $\rho(z)=\inf\{R>0 : \U\subset B_{R}(z)\}$ is the maximal visibility radius from $z\eqdef\widehat Z\in\widehat\V$,
\item[2] if (\ref{vis-cond}) and (\ref{dist-cond}) hold then $\alpha_{u,g}$ is countably additive,
\item[3] if (\ref{vis-cond})-(\ref{Rconv-cond}), (\ref{vis-cond-B}) hold and $f\ge 0$ while $g>0$ then $B$-type weak solution is also of $A$ type,
\item[4] if we remove the $R-$convexity assumption but require the positivity of densities  (\ref{dens-cond}) and
  (\ref{vis-cond})-(\ref{dist-cond}), (\ref{vis-cond-B}) then again any $B$-type weak solution is also of $A$-type.
\end{itemize}
\end{CCC}
The proof of Theorem C1 is by polyhedral approximation and utilising the confocal expansion of hyperboloids as described in 
Section \ref{ssec-confocal}.   In this regard  the condition (\ref{vis-cond-B}) in Theorem C1 says that 
one can construct a $B$-type weak solution if there is sufficient span between $\Pi$ and $\Sigma$. 
%%%
The existence of $B$-type weak solutions, constructed from an envelope of ellipsoids of revolution can be found in \cite{GutT}. 
%%%
%%%

Our last result concerns with the smoothness of $A$-type weak solutions. We use the well-known method of  comparing the mollified weak solution with that of Dirichlet's problem to the slightly modified equation in a small ball $B$. 
To this end one first has to obtain $C^{2, \alpha}, \alpha\in(0, 1]$ estimates in $\ol B$ for the solutions of mollified equations and after that  making sure  that uniform $C^{2}$ estimates hold in, say, $\frac12 B$. Then passing to limit and using the comparison principle the result will follow. 
The construction of weak solutions to Dirichlet's problem
is based on  Perron's method and follows the approach developed by Xu-Jia Wang in \cite{W-96} where 
a far field reflector design problem is studied.  Our research is inspired by \cite{W-96} and subsequent developments in \cite{KW-1}, \cite{KW-2} \cite{K-ref}. For more recent results on this problem see \cite{Liu}. 
The global $C^2$ estimates for the solution of Dirichlet's problem for the regularised equation follow from \cite{Turd-chin-arxiv} whereas
the local  uniform estimates in $\frac12 B$ are established in \cite{MTW}, see also \cite{Liu-Trud}
for the global regularity of  near field reflector problem with point source of a light.
Thus we have the following theorem 
 
\begin{DDD}
 Let $f, g$ be $C^2$ smooth functions such that $ \lambda\leq f,g\leq \Lambda$ for some constants $\Lambda>\lambda>0$ and the conditions (\ref{vis-cond})-
 (\ref{A3-cond00})
  are satisfied.
  Then $A$-type weak solutions of (\ref{problem}) are locally $C^2$ regular in $\U$.
\end{DDD}
The conditions (\ref{vis-cond})-
 (\ref{A3-cond00}) cannot be relaxed as one may easily construct counterexamples to regularity in the spirit of those in \cite{KW-1}, \cite{KW-2}. For instance let us examine (\ref{Rconv-cond}) (see also Remark \ref{rem-zzbg}), if 
 we take a  two point target $\V=\{Z_1\}\cup\{Z_2\}$ and consider $H(x)=\min[H(x, a_1, Z_1), H(x, a_2, Z_2)]$ such that these hyperboloids $H(\cdot, a_i, Z_i), i=1,2$ have non empty intersection over $\U$. Then 
 approximating $\V$ by smooth R-convex sets $\V_t$ we obtain  a sequence of admissible functions $H_t$, solving the 
 refractor problem with target $\V_t$, and converging  to $H$ as $t\to 0$.
But if $t$ is sufficiently close to 0 then $H_t$ cannot remain $C^1$ smooth because otherwise the limit $H$ would also be $C^1$ which is impossible., see \cite{KW-1} for more discussion on such constructions.  
{We would also like to point out 
a recent paper of Guti\'erred and Tournier \cite{GutTC1} where the authors study the local $C^{1, \alpha}$ regularity of reflector/refractor problems without using the explicit form of the equation. 
There one can find a detailed account of the case when the supporting functions are 
ellipsoids of revolution. Our work contributes in this direction only by deriving the explicit equation for 
general receiver $\Sigma$ and establishing a simple form  of the corresponding regularity condition 
(involving the second fundamental form of $\Sigma$) for the existence of 
$C^2$ smooth solutions, see Section \ref{ssec-example}. In this paper, however, we mainly focus on the case when the 
supporting functions are the hyperboloids of revolution. }

\smallskip

The rest of the paper is organized as follows: in the next section we derive the
main formulae. Then we prove Theorem A in Section 4. The main result there
is Proposition \ref{th-eq} from which the proof of Theorem A easily follows.
Section 5 contains some preliminary discussion on the condition (\ref{A3-cond00}) and after that in Section 6 we give the proof of Theorem B. The admissible functions are introduced in Section 7 where we also exhibit some interesting properties of hyperboloids of revolution, notably the dual admissibility and confocal expansion.  
 Employing the polyhedral approximation technique and weak convergence of measures $\beta_{u,f}$ we prove Theorem C1 in Section 8.  
%%% Sec 9 approximation lemma 
The first direct application of (\ref{A3-cond00}) is given in Lemma \ref{loc-glob}, which is G. Loeper's geometric interpretation of the A3 condition  from \cite{MTW}. A direct consequence  of this is Lemma \ref{lem-apprx} stating that 
a suitable dilation of an admissible function by a paraboloid  of revolution can be approximated via 
 smooth subsolutions of (\ref{eq-main-th}). This is a crucial ingredient in the proof of Theorem D.
%%% Sec 10 Legendre + Theorem C2
Next we introduce the Legendre-like transformation of an admissible $u$ and conclude Theorem C2. 
%%% Sec 11 Comparing A and B Theorem C3-4.
The proofs of Theorem C3-4 follow from a comparison of $A$ and $B$ type weak solutions
by extending the results of Luis Caffarelli \cite{C92} and John Urbas \cite{Urbas} for the classical Monge-Amp\`ere equation to (\ref{eq-main-th}). This is done in Section 11.
%%% Sec 12 Weak Dirichlet + Perron
The last two sections are devoted to the study of the higher regularity of  $A$-type weak solutions. 
We follow the classical approach developed by A. Pogorelov for the classical Monge-Amp\`ere equation, see \cite{Pog-2D}, \cite{P-Mink}. Therefore, we first 
prove the solvability of weak Dirichlet's problem when the boundary data is given as the trace of an $A$-type weak subsolution. That done, the uniqueness  follows from comparison principle stated in  Proposition \ref{prop-compar}.
%%% Section 13
Finally in Section 13 we give the proof of our main regularity result, Theorem D.

\section{Notations}
\begin{tabbing}
$C, C_0, C_n, \cdots$ \hspace{1.55cm}       \=\hbox{generic constants,}\\
$\Pi$     \>     $\Pi=\R^n\times \{0\}$,\\
$\overline \Omega$       \>\hbox{closure of a set} $\Omega$,\\
$\partial \Omega$        \>\hbox{boundary of a set }  $\Omega$, \\
$\widehat \Omega$\> the projection of $\Omega\subset \R^{n+1}$ on $\Pi$,\\
$\widehat X$\> $(x_1, x_2, \dots, x_n, 0)$ projection of $X=(x_1, x_2, \dots, x_n, x_{n+1}),$\\
$\ee$\>eccentricity, \\
$\kap$\> $\kap=\frac{\eps^2-1}{\eps^2}$,\\
$\H$\> $n$-dimensional Hausdorff measure on $\Sigma$,\\
$\p_{i}$\> partial derivate with respect to $x_i$ variable,\\
$Du$\> the gradient of a function $u$,\\
$\rho(z)$  \> $\inf\{R>0 : \U\subset B_{z}(R)\}$ is the maximal visibility radius from $z\eqdef\widehat Z\in\widehat\V$,\\
$q$\> see (\ref{def-q}),\\
$\mathbb H (\U, \Sigma)$\> the class of hyperboloids of revolution with focal axis parallel to $Z^{n+1}$ and upper focus on $\Sigma$,\\
$\mathbb H^+_{a_0}(\U, \V)$\> hyperboloids from $\mathbb H(\U, \V)$ which are nonnegative in  $\U$ and $a> a_0$ for some fixed $a_0$,\\
$\Wup, \Wlo$\>  upper and lower admissible functions, see Lemma \ref{lem-lips},\\
$\Wup^0(\U, \V)$\> polyhedral admissible functions.
\end{tabbing}
%%%%%%%%%%%%%%%%%%%%%%%%%%%%%%%%%%%%%%%%%%%%%%%%%%%%%%%%%%%%
%                                                          %
%                SECTION                                   %
%                                                          %
%%%%%%%%%%%%%%%%%%%%%%%%%%%%%%%%%%%%%%%%%%%%%%%%%%%%%%%%%%%%
\section{Main formulae}
In this section we derive the Monge-Amp\`ere type equation \eqref{eq-main-th}
manifesting  the energy balance condition (\ref{eq-energ-bal})
in the refractor problem (\ref{problem}), see Introduction.

\subsection{Computing $Y$}
We first compute the unit direction of the refracted ray.
Denote by $\gamma$ the unit normal to the graph of $u$, that is
\begin{eqnarray}\label{eq-gamma}
  \gamma=\frac{(-D_1u,\dots,-D_nu,1)}{\sqrt{1+|Du|^2}}.
\end{eqnarray}
Since $\ell_x, Y$ and $\gamma$
lie in the same hyperplane we have
\begin{eqnarray}\label{eq-Y-comb}
  Y=\a e_{n+1}+\b\gamma,
\end{eqnarray}
for some coefficients $\a$ and $\b$.
Computing the scalar products $Y\cdot \gamma$ and $Y\cdot e_{n+1}$ we obtain the following equations (cf. (\ref{angles}))
\begin{eqnarray*}
  \left\{
  \begin{array}{ll}
    \cos\theta_2=\a\cos\theta_1+\b,\\
    \cos{(\theta_1-\theta_2)}=\a+\b\cos\theta_1.
  \end{array}
  \right.
\end{eqnarray*}
Multiplying the first equation by $\cos\theta_1$ and subtracting from the second one we conclude
\begin{eqnarray*}
  \a=\frac{\sin\theta_2}{\sin\theta_1}, \hspace{2cm} \b=\cos\theta_2-\a\cos\theta_1.
\end{eqnarray*}

\smallskip

Recalling our notations
\begin{equation}\label{def-sig}
 \kap=\frac{\eps^2-1}{\eps^2}, \qquad \eps=\frac{n_1}{n_2},
\end{equation}
we see that  $\mathpzc A=\eps$. Furthermore
\begin{eqnarray*}
n_2^2-n_2^2\cos^2\theta_2=n_2^2\sin^2\theta_2=n_1^2\sin^2\theta_1=n_1^2-n_1^2\cos^2\theta_1.
\end{eqnarray*}
Dividing both sides  of this identity by $n_2^2$ we obtain
$$\cos^2\theta_2=\eps^2\cos^2\theta_1-(\eps^2-1)=\eps^2(\cos^2\theta_1-\kap).$$
Therefore from  $\a=\eps$ we conclude that $\b=\eps(\sqrt{\cos^2\theta_1-\kap}-\cos\theta_1)$. Returning to
(\ref{eq-Y-comb})
we infer that the unit direction of the refracted ray is
\begin{eqnarray}\label{eq-Y}
 Y=\eps\left(e_{n+1}+(\sqrt{\cos^2\theta_1-\kap}-\cos\theta_1)\gamma\right).
\end{eqnarray}
Notice that (\ref{eq-gamma}) implies
$$\cos\theta_1=\gamma\cdot e_{n+1}=\frac1{\sqrt{1+|Du|^2}}.$$
Consequently, denoting  $Y=(Y^1,Y^2,\dots,Y^n,Y^{n+1})$
and $y\in\R^n$, the projection of $Y$ onto $\Pi=\{X\in\R^{n+1} : X^{n+1}=0\}$, (i.e. $y=(Y^1,Y^2,\dots,Y^n,0)$) we get
\begin{eqnarray}\label{y-par}
  y&=&-\eps\frac{Du}{\sqrt{1+|Du|^2}}\left(\sqrt{\cos^2\theta_1-\kap}-\cos\theta_1\right)\\\nonumber
  &=&\frac{\eps\kap Du}{\sqrt{1+|Du|^2}}\frac{1}{\sqrt{\cos^2\theta_1-\kap}+\cos\theta_1}\\\nonumber
  &=&\eps\kap\frac{Du}{\sqrt{1-\kap(1+|Du|^2)}+1}.
\end{eqnarray}
From this computation it follows that
\begin{equation}\label{y-npls1}
 Y^{n+1}=\eps\left(1-\frac{\kap}{1+ \sqrt{1-\kap(1+|Du|^2)}}\right).
\end{equation}

Combining \eqref{y-par} and \eqref{y-npls1} we obtain
\begin{equation}\label{unit-Y}
 Y=\eps\left(\frac{\kap Du}{1+\sqrt{1-\kap(1+|Du|^2)}}, 1-\frac{\kap}{1+ \sqrt{1-\kap(1+|Du|^2)}}\right).
\end{equation}
If we use the notation $q(x)=\sqrt{1-\kap(1+|Du|^2)}$ (see (\ref{def-q})) then (\ref{unit-Y}) takes the form
\begin{equation}\label{unit-Y-q}
 Y=\eps\left(\frac{\kap Du}{1+q}, 1-\frac{\kap}{1+q}\right).
\end{equation}
Notice that by \eqref{unit-Y} $Y^{n+1}>0$ for all values of $\kap$.
\smallskip

\subsection{Stretch function}
Assume that $\psi $ is a smooth function $\psi:\R^{n+1}\to\R$, and the
receiver $\Sigma$ is given as the zero set of $\psi$
\begin{eqnarray}
  \Sigma=\{Z\in \R^{n+1} : \psi(Z)=0\}.
\end{eqnarray}
Let us represent the mapping $Z: \U\to \Sigma$ in the following form
\begin{equation}\label{def-Z}
Z=x+e_{n+1}u(x)+Yt,
\end{equation}
where $t=t(x, u(x), Du(x))$ is determined from the equation
$\psi(Z)=0$ and  is called  the \textit{stretch function}.
It is worthwhile to point out that
the stretch function $t$ can be explicitly computed for a wide class of elementary surfaces.
For instance, if $\Sigma$ is the horizontal plane $Z^{n+1}=m>0$ then from simple geometric considerations one finds that
$$t=\frac{m-u}{Y^{n+1}}$$
where $Y^{n+1}$ is given by \eqref{y-npls1}.

In lemma to follow
we denote by $z$ the projection of $Z$ onto $\Pi$, that is
$z=x+ty$.

\begin{lemma}\label{lem-1}
Let $dS_\U$ and $dS_\V$ be  the area elements on  $\U$ and $Z(\U)=\V\subset\Sigma$ respectively
and $z$ being the projection of $Z$ onto $\Pi=\{Z\in \R^{n+1} :Z^{n+1}=0\}.$
Then  we have
\setlength\arraycolsep{2pt}
\begin{eqnarray}\label{eq-Jac}
  J & =&\frac{dS_\V }{dS_\U}
   =\left|\begin{array}{cccc}
        Z_{1}^1, &\cdots, & Z^{1}_n, & \nu^1\\
         \vdots  &\ddots &\vdots &\vdots \\
        Z^{n}_1, &\cdots, & Z_{n}^n, & \nu^n\\
        Z^{n+1}_1, &\cdots, & Z^{n+1}_n, & \nu^{n+1}\\
         \end{array} \right|\\\nonumber
&=&-\frac{|\nabla \psi|}{\psi_{n+1}}\det Dz,\\\nonumber
\end{eqnarray}
where $\nu$ is the unit normal of $\Sigma$.
\end{lemma}
\pr The first equality in (\ref{eq-Jac}) follows from the change of variables formula.
Differentiating the equality $\psi(Z)=0$ by $x_i$ we have that
$$\p_iZ^{n+1}=-\frac1{\p_{n+1}\psi}\sum\limits_{k=1}^n\p_iz^k\p_{z_k}\psi.$$
Using this identity we multiply $j$-th row of matrix in (\ref{eq-Jac}) by $\p_{z_j}\psi$ and subtract it
from the $(n+1)$st row in order to get
\setlength\arraycolsep{2pt}
\begin{eqnarray*}
 \det\left|\begin{array}{cccc}
        Z_{1}^1, &\cdots, & Z_{n}^1, & \nu_1\\
         \vdots  &\ddots &\vdots &\vdots \\
        Z^{n}_1, &\cdots, & Z^{n}_n, & \nu_n\\
        Z^{n+1}_1, &\cdots, & Z^{n+1}_n, & \nu_{n+1}\\
         \end{array} \right|&=&
%%%%%%%%%%%%%%%%%%%%%%%%%%%%%%%%%%%%%%%%%%%%%%%%%%%%%%%%%%%%%
-\frac1{\psi_{n+1}}\det \left|\begin{array}{cccc}
        Z_{1}^1, &\cdots, & Z_{n}^1, & \nu_1\\
         \vdots  &\ddots &\vdots &\vdots \\
        Z^{n}_1, &\cdots, & Z^{n}_n, & \nu_n\\
        \sum\limits_{k=1}^n\p_1z^k\p_{z_k}\psi, &\cdots, & \sum\limits_{k=1}^n\p_nz^k\p_{z_k}\psi, & -\psi_{n+1}\nu_{n+1}\\
         \end{array} \right|\\\nonumber
%%%%%%%%%%%%%%%%%%%%%%%%%%%%%%%%%%%%%%%%%%%%%%%%%%%%%%%%%%%%%
&=&-\frac1{\psi_{n+1}}\det \left|\begin{array}{cccc}
        Z_{1}^1, &\cdots, & Z_{n}^1, & \nu_1\\
         \vdots  &\ddots &\vdots &\vdots \\
        Z^{n}_1, &\cdots, & Z^{n}_n, & \nu_n\\
        0, &\cdots, & 0, & -\sum\limits_{k=1}^{n+1}\psi_{k}\nu_{k}\\
         \end{array} \right|.
\end{eqnarray*}
Finally noting that $\nu=\frac{\nabla\psi}{|\nabla \psi|}$ the desired identity follows.\qed

\begin{lemma}\label{lem-prelim-comp}
Let $C\in \R$ and $\xi, \eta\in \R^n$. Consider the matrix $\mu=\Id+C\xi\otimes \eta=\delta_{ij}+C\xi^i\eta^j$ where $\Id=\delta_{ij}$ is the identity matrix.
Then  the inverse matrix of $\mu $ is
\begin{eqnarray*}
\det\mu&=&1+C\xi\cdot\eta,\\
 \mu^{-1}&=&\Id-\frac{C\xi\otimes\eta}{1+C(\xi\cdot\eta)}.
\end{eqnarray*}
Here and henceforth $\Id$ is the identity matrix.
\end{lemma}
\pr\ Without loss of generality we assume that $\xi=e_1$ then $\det\mu=1+C\eta^1$.
As for the second identity, it is a partucal case of Sherman-Morrison formula. \qed

\smallskip

Finally, we derive a formula for the first order derivatives
of the stretch function $t$.
Let us differentiate the equation $\psi(Z)=0$ with respect to $x_j$ to get
$$ \sum_{k=1}^n\psi_k(\delta_{kj}+t_jy^k+ty^k_j)+\psi_{n+1}(u_j+t_jY^{n+1}+tY^{n+1}_j)=0.$$
From here we find
\begin{eqnarray}\label{eq-Dt}
  t_j=-\frac{1}{\nabla\psi\cdot Y}[\psi_j+\psi_{n+1}u_j+t(\nabla\psi\cdot Y_j)].
\end{eqnarray}

\medskip

%%%%%%%%%%%%%%%%%%%%%%%%%%%%%%%%%%%%%%%%%%%%%%%%%%%%%%%%%%%%
%                                                          %
%                SECTION                                   %
%                                                          %
%%%%%%%%%%%%%%%%%%%%%%%%%%%%%%%%%%%%%%%%%%%%%%%%%%%%%%%%%%%%
\section{Proof of Theorem A}\label{sec-eq-derv}
In this section we prove Theorem A. We begin with a computation for the matrix $Dz$, where $z$ is the projection of
$Z$ on to $\Pi$.

\begin{proposition}\label{th-eq}
 Let $u\in C^2(\U)$ and $Z$ be the corresponding refractor map, then with the same notations as in Lemma \ref{lem-1}
 we have
  \begin{eqnarray}
    Dz&=&\mu_1\mu_2\left[\Id-\kap\eps^2 Du\otimes Du+\frac{t\kap\eps}{1+h}D^2u\right],
    \\\nonumber
  \end{eqnarray}
where
\begin{eqnarray}\label{mtrx-def}
       \mu_1=\Id-\frac{y\otimes(\widehat\nabla\psi-y\frac{\psi_{n+1}}{Y^{n+1}})}{\nabla\psi\cdot Y}, \qquad
    \mu_2=\Id+\kap\frac{Du\otimes Du}{q(q+1)},
      \end{eqnarray}
$q=\sqrt{1-\kap(1+|Du|^2)}$ and
\begin{equation}\label{D-psi}
 \widehat \nabla \psi=(\psi_1, \dots, \psi_n, 0).
\end{equation}
\end{proposition}

\smallskip

In order to prove Proposition \ref{th-eq} we will need the following
\begin{lemma}\label{lem-Dz-comp}
Let $z(x), x\in \U$ be the projection of the mapping $Z(x)$ onto $\Pi=\{X\in \R^{n+1} : X^{n+1}=0\}$. Then
  \begin{eqnarray}
    Dz=\mu_1\left(\Id -y\otimes[y+DuY^{n+1}]+tDy\right)
  \end{eqnarray}
where $\mu_1$ is defined by (\ref{mtrx-def}).
\end{lemma}
\pr Introduce the matrix
\begin{equation}\label{mu-0}
 \mu_0=\delta_{ij}-y^i\frac{\psi_j+u_j\psi_{n+1}}{\nabla\psi\cdot Y}.
\end{equation}
Using (\ref{eq-Dt}) and recalling $z=x+ty$ we compute
\begin{eqnarray}\label{Dz-1}
  z^i_j&=&\delta_{ij}+t_jy^i+ty^i_j\\\nonumber
  &=&\delta_{ij}+ty^i_j-y^i\frac{1}{\nabla\psi\cdot Y}[\psi_j+\psi_{n+1}u_j+t(\nabla\psi\cdot Y_j)]\\\nonumber
  &=&\underbrace{\delta_{ij}-y^i\frac{[\psi_j+u_j\psi_{n+1}]}{\nabla\psi\cdot Y}}_{\mu_0}
  +t\left[y^i_j -\frac{y^i(\nabla\psi\cdot Y_j)}{\nabla\psi\cdot Y}\right]\\\nonumber
  &=&\mu_0+t\left[y^i_j -\frac{y^i(\nabla\psi\cdot Y_j)}{\nabla\psi\cdot Y}\right]
.\\\nonumber
\end{eqnarray}
In order to deal with the remaining matrix we recall that
$(Y^{n+1})^2=1-|y|^2$ and hence $Y^{n+1}_j=-\frac{yy_j}{Y^{n+1}}$.
Consequently, setting $\widehat\nabla \psi=(\psi_1, \dots,\psi_n,0)$ (see (\ref{D-psi})) we infer
\begin{eqnarray}\label{Dz-2}
  y^i_j-\frac{y^i(\nabla\psi\cdot Y_j)}{\nabla\psi\cdot Y}&=&y^i_j-\frac{y^i}{\nabla\psi\cdot Y}
  \left(\widehat\nabla \psi\cdot y_j-\psi_{n+1}\frac{y\cdot y_j}{Y^{n+1}}\right)\\\nonumber
  &=&y^i_j-\frac{y^i}{\nabla\psi\cdot Y}\left[\left(\widehat\nabla \psi-\frac{\psi_{n+1}}{Y^{n+1}}y\right)y_j\right].
\end{eqnarray}
Combining (\ref{Dz-1}) and \eqref{Dz-2} we obtain the following formula for $Dz$, written in intrinsic form
\begin{eqnarray}\label{Dz-122}
 Dz&=&\mu_0+t\left[\Id-\frac{y\otimes\left(\widehat\nabla \psi-\frac{\psi_{n+1}}{Y^{n+1}}y\right)}{\nabla \psi \cdot Y}\right]Dy\\\nonumber
 &=&\mu_0+t\mu_1Dy\\\nonumber
 &=&\mu_1(\mu_1^{-1}\mu_0+tDy)
\end{eqnarray}
where the second equality follows from the definition of matrix $\mu_1$, see (\ref{mtrx-def}).

Next, we compute $\mu_1^{-1}$. From Lemma \ref{lem-prelim-comp} and the identity $[Y^{n+1}]^2=1-|y|^2$ we get
\begin{eqnarray}
 \mu_1^{-1}&=& \Id+y\otimes\frac{\widehat\nabla \psi-\frac{\psi_{n+1}}{Y^{n+1}}y}{\nabla\psi\cdot Y-
 \left(\widehat\nabla \psi\cdot y-|y|^2\frac{\psi_{n+1}}{Y^{n+1}}\right)}\\\nonumber
 &=&\Id+\frac{Y^{n+1}}{\psi_{n+1}}y\otimes\left[\widehat\nabla\psi-\frac{\psi_{n+1}}{Y^{n+1}}y\right],
\end{eqnarray}
where the last equality follows from the observation
\begin{eqnarray}
\nabla\psi\cdot Y-  \left(\widehat\nabla \psi\cdot y-|y|^2\frac{\psi_{n+1}}{Y^{n+1}}\right) &=&\psi_{n+1}Y^{n+1}+(1-(Y^{n+1})^2)\frac{\psi_{n+1}}{Y^{n+1}}
\\\nonumber
&=&\frac{\psi_{n+1}}{Y^{n+1}}.
\end{eqnarray}
It is convenient to rewrite this identity in the following form
\begin{equation}\label{obvs}
 \left[\widehat\nabla \psi\cdot y-|y|^2\frac{\psi_{n+1}}{Y^{n+1}}\right]\frac{Y^{n+1}}{\psi_{n+1}}\frac{1}{\nabla\psi\cdot Y}=
\frac{Y^{n+1}}{\psi_{n+1}}-\frac{1}{\nabla\psi\cdot Y}.
\end{equation}

Consequently, we obtain
\begin{eqnarray*}
 \mu_1^{-1}\mu_0&=&\left(\Id+\frac{Y^{n+1}}{\psi_{n+1}}y\otimes \left[\widehat \nabla\psi-\frac{\psi_{n+1}}{Y^{n+1}}y\right]\right)
 \left(\Id -y\otimes\frac{\widehat\nabla\psi+Du\psi_{n+1}}{\nabla\psi\cdot Y}\right)\\\nonumber
 &=&\Id-y\otimes\frac{\widehat\nabla\psi+Du\psi_{n+1}}{\nabla\psi\cdot Y}
 +\frac{Y^{n+1}}{\psi_{n+1}}y\otimes \left[\widehat \nabla\psi-\frac{\psi_{n+1}}{Y^{n+1}}y\right]-\\\nonumber
 &&-\left[\widehat\nabla \psi\cdot y-|y|^2\frac{\psi_{n+1}}{Y^{n+1}}\right]
 \frac{Y^{n+1}}{\psi_{n+1}}\frac{1}{\nabla\psi\cdot Y}\left\{y\otimes {\widehat\nabla\psi+Du\psi_{n+1}}\right\}.
\end{eqnarray*}
Applying (\ref{obvs}) to the last term in this computation we get
\begin{eqnarray*}
 \mu_1^{-1}\mu_0
  &=&\Id-y\otimes\frac{\widehat\nabla\psi+Du\psi_{n+1}}{\nabla\psi\cdot Y}
 +\frac{Y^{n+1}}{\psi_{n+1}}y\otimes \left[\widehat \nabla\psi-\frac{\psi_{n+1}}{Y^{n+1}}y\right]-\\\nonumber
 &&-\left[\frac{Y^{n+1}}{\psi_{n+1}}-\frac{1}{\nabla\psi\cdot Y}\right]
 \left\{y\otimes {\widehat\nabla\psi+Du\psi_{n+1}}\right\}\\\nonumber
 %%%%
 &=&\Id +\frac{Y^{n+1}}{\psi_{n+1}}y\otimes \left[\widehat \nabla\psi-\frac{\psi_{n+1}}{Y^{n+1}}y\right]-\\\nonumber
 &&-\frac{Y^{n+1}}{\psi_{n+1}}
 \left\{y\otimes {\widehat\nabla\psi+Du\psi_{n+1}} \right\}\\\nonumber
 %%%%
 &=&\Id-y\otimes	[y+Du\psi_{n+1}].
\end{eqnarray*}
Plugging in the computed form of $\mu_1^{-1}\mu_0$ into \eqref{Dz-122} the result follows.
\qed

\medskip

\subsection{Proof of Proposition \ref{th-eq}}

To finish the proof of Proposition \ref{th-eq}, it remains to express
$Dz$ through the Hessian $D^2 u$.
We have from (\ref{unit-Y-q})
\begin{eqnarray}\label{eq-y-sig}
  y&=&\eps\kap\frac{Du}{q+1},\\
  Y^{n+1}&=&\eps\left(1-\frac{\kap}{q+1}\right),
\end{eqnarray}
where $q=\sqrt{1-\kap(1+|Du|^2)}$, see (\ref{def-q}).
From the definition of $q$
we have $Dq=-\kap DuD^2u/q$, thus
\begin{eqnarray*}
Dy&=&\eps\kap\left[\Id+\kap\frac{Du\otimes Du}{q(q+1)}\right]\frac{D^2u}{q+1}\\\nonumber
&=&\eps\kap\mu_2 \frac{D^2u}{q+1},
\end{eqnarray*}
where $\mu_2$ is the matrix in (\ref{mtrx-def}). Now Lemma \ref{lem-Dz-comp} yields
\begin{eqnarray}\label{443}
 Dz&=&\mu_1\left(\Id-y\otimes[y+DuY^{n+1}]+t\eps\kap\mu_2 \frac{D^2u}{q+1}\right)\\\nonumber
 &=& \mu_1\mu_2\left(\mu_2^{-1}\left\{\Id-y\otimes[y+DuY^{n+1}]\right\}+t\eps\kap\frac{D^2u}{q+1}\right)\\\nonumber
  &=& \mu_1\mu_2\left(\mu_2^{-1}\mathcal M+t\eps\kap\frac{D^2u}{q+1}\right)
\end{eqnarray}
where $\mathcal M=\Id-y\otimes[y+DuY^{n+1}]$.

\smallskip

Using (\ref{eq-y-sig}) we can further simplify the matrix $\mathcal M=\Id-y\otimes[y+DuY^{n+1}]$ to get
\begin{eqnarray}\label{mu-2}
\mathcal M&=&\Id - y\otimes(y+DuY^{n+1})\\\nonumber
&=&\Id -\frac{\eps^2\kap^2}{(1+q)^2}Du\otimes Du
-\frac{\eps^2\kap (1-\frac\kap{1+q})}{1+q}Du\otimes Du\\\nonumber
&=&\Id-\frac{\eps^2\kap}{1+q}Du\otimes Du.
\end{eqnarray}
By  Lemma \ref{lem-prelim-comp} we have for the inverse of $\mu_2$ (see \eqref{mtrx-def})
\begin{eqnarray}\label{mu-3-inv}
 \mu_2^{-1}&=&\Id-\frac{\kap}{q^2+q+\kap|Du|^2}Du\otimes Du\\\nonumber
&=&\Id -\frac{\kap}{1-\kap+q}Du\otimes Du,
\end{eqnarray}
where the last equality  follows from the definition of $q$, see (\ref{def-q}).
\smallskip
It remains to compute
$\mu_2^{-1}\cal M$. From  (\ref{mu-3-inv}) and (\ref{mu-2}) we obtain

\begin{eqnarray*}
 \mu_2^{-1}\mathcal M&=& \left[\Id -\frac{\kap}{1-\kap+q}Du\otimes Du\right]\left[\Id-\frac{\eps^2\kap}{1+q}
Du\otimes Du\right]\\\nonumber
&=&\Id +\left[I+II+III\right]Du\otimes Du
\end{eqnarray*}
where
\begin{eqnarray*}
 I&=&-\frac{\kap}{1-\kap+q},\\\nonumber
II&=& -\frac{\eps^2\kap}{1+q},\\\nonumber
III&=&\frac{\eps^2\kap^2|Du|^2}{(1+q)(1-\kap+q)}.
\end{eqnarray*}
It follows from (\ref{def-q}) that $-\kap|Du|^2=q^2-1+\kap$, therefore
$$III=\frac{\eps^2\kap(-q^2+1-\kap)}{(1+q)(1-\kap+q)}.$$
Adding this to $II$ we have
\begin{eqnarray*}
 II+III&=&\frac{\eps^2\kap}{1+q}\left[-1+\frac{-q^2+1-\kap}{1-\kap+q}\right]\\\nonumber
&=&-\frac{q\eps^2\kap}{1-\kap+q}.
\end{eqnarray*}

Finally we compute  the total sum
\begin{eqnarray*}
 I+II+III&=&-\frac{\kap}{1-\kap+q}-\frac{q\eps^2\kap}{1-\kap+q}\\\nonumber
&=& -\frac{\kap}{1-\kap+q}\left[q\eps^2 +1\right]\\\nonumber
&=&  -\frac{\kap}{1-\kap+q}\left[\frac{q}{1-k} +1\right]\\\nonumber
&=&-\frac{\kap}{1-\kap}\\\nonumber
&=&-\kap\eps^2,
\end{eqnarray*}
where the last line follows from the definition of $\kap$,
see (\ref{def-sig}).

\smallskip

Returning to (\ref{443}) and utilising these computations
we get

\begin{eqnarray*}
  Dz&=&\mu_1\mu_2\left[ \mu_2^{-1}\mathcal M  +t\eps\kap\frac{D^2u}{q+1} \right]\\\nonumber
&=&\mu_1\mu_2\left[\Id -\kap\eps^2Du\otimes Du+t\eps\kap\frac{D^2u}{q+1}\right].\\\nonumber
\end{eqnarray*}
This finishes the proof of Proposition  \ref{th-eq}. \qed

\smallskip

\subsection{Proof of Theorem A}\label{cor-with-det}
Now we are ready to finish the proof of
Theorem A. Let $u\in C^2(\U)$ be a solution to the refractor problem (\ref{problem})
then from Proposition \ref{th-eq} we obtain
\begin{equation}
\det Dz=\det\mu_1\det\mu_2\left[\frac{t\eps\kap}{q+1}\right]^n\det\left[\frac{q+1}{t\eps\kap}\left\{\Id-\kap\eps^2Du\otimes Du\right\}+D^2u\right].
\end{equation}
 By Lemma \ref{lem-prelim-comp} and (\ref{def-q}) we have
$$\det \mu_2= 1+\frac{\kap|Du|^2}{q(q+1)}=\frac{1-\kap+q}{q(q+1)}.$$
Similarly, we get
$$\det\mu_1=\frac{\psi_{n+1}}{Y^{n+1}}\frac1{\nabla \psi\cdot Y}.$$
These in conjunction with (\ref{eq-Jac}) gives

\begin{eqnarray*}
 \det\left[\frac{q+1}{t\eps\kap}\left\{\Id-\kap\eps^2Du\otimes Du\right\}+D^2u\right] &=&\left[\frac{q+1}{t\eps\kap}\right]^n\frac{\det Dz}{\det\mu_1\det\mu_2}\\\nonumber
 &=&-\frac{f}g \frac{\psi_{n+1}}{|\nabla\psi|}\left[\frac{q+1}{t\eps\kap}\right]^n\frac{1}{\det\mu_1\det\mu_2}\\\nonumber
 &=&-
(\nabla \psi\cdot Y)\frac{Y^{n+1}}{|\nabla \psi|}\frac{q(q+1)}{1-\kap+q}\left[\frac{q+1}{t\eps\kap}\right]^n\frac fg.
 \end{eqnarray*}
 Finally, recalling (\ref{unit-Y-q}) and substituting the value of $Y^{n+1}$ we see that
\begin{eqnarray}\label{eq-with-det}
 \det\left[\frac{q+1}{t\eps\kap}\left\{\Id-\kap\eps^2Du\otimes Du\right\}+D^2u\right]&=&-
(\nabla \psi\cdot Y)\frac{Y^{n+1}}{|\nabla \psi|}\frac{q(q+1)}{1-\kap+q}\left[\frac{q+1}{t\eps\kap}\right]^n\frac fg\\\nonumber
&=&-\eps q \left[\frac{q+1}{t\eps\kap}\right]^n
\frac{\nabla \psi\cdot Y}{|\nabla \psi|}\frac fg
\end{eqnarray}
and the proof of Theorem A is now complete.
\qed

\medskip

%%%%%%%%%%%%%%%%%%%%%%%%%%%%%%%%%%%%%%%%%%%%%%%%%%%%%%%%%%%%
%                                                          %
%                SECTION                          %
%                                                          %
%%%%%%%%%%%%%%%%%%%%%%%%%%%%%%%%%%%%%%%%%%%%%%%%%%%%%%%%%%%%
\section{Existence of smooth solutions}\label{sec-weak-sol}

In this section we will have a provisional discussion on the   existence of smooth solutions to (\ref{eq-main-th}).  
Our main objective is to apply the available regularity theory for the Monge-Amp\`ere type equations, stemming from seminal paper
\cite{MTW},  in order to establish the regularity of
weak solutions of the refractor problem.

We first rewrite the equation (\ref{eq-with-det}) in a more concise form.
Let us introduce the following matrix
\begin{equation}\label{def-G-M}
 G^{ij}=\frac1t %M^{ij}, \qquad M^{ij}=
 (q+1)[\delta_{ij}-\kap\eps^2u_iu_j].
\end{equation}
Here $q=\sqrt{1-\kap(1+|Du|^2)}$, see (\ref{def-q}) and $t$ is the stretch function determined from implicit equation $\psi(x+e_{n+1}u +tY)=0$
as in Theorem A.
Then the equation (\ref{eq-with-det}) transforms into
\begin{eqnarray}\label{M-A-short+}
 \det\left[-\frac G{\eps\kap }-D^2u\right]&=&|h(x, u, Du)|, \quad {\rm if} \ \kap>0,\eps>1\quad u\in C^2(\U)\ {\rm and}\ -\frac G{\eps\kap }-D^2u\ge0,\ %(\mbox{supporting hyperboloids})
 \\\label{M-A-short-}
 \det\left[D^2u +\frac G{\eps\kap }\right]&=&|h(x, u, Du)|,\quad  {\rm if} \ \kap<0,\eps<1\quad u\in C^2(\U)\ {\rm and}\ \frac G{\eps\kap }+D^2u\ge0%(\mbox{supporting ellipsoids})
\end{eqnarray}
with
\begin{equation}\label{def-rhs}
h(x, u, Du)=
-\eps q \left[\frac{q+1}{t\eps\kap}\right]^n
\frac{\nabla \psi\cdot Y}{|\nabla \psi|}\frac fg.
\end{equation}

\smallskip

The existence of $C^2$ smooth solutions of (\ref{M-A-short+}) or (\ref{M-A-short-}) depend on the
properties of the matrix $G$.
Namely, it is shown in \cite{MTW} that if we regard
$G$ as a function of variable $p=Du$ then the condition
\begin{equation}\label{eq-A3}
- D^2_{p_kp_l}G^{ij}\xi_i\xi_j\eta_k\eta_l\
\begin{array}{lll}
\leq -c_0|\xi|^2|\eta|^2 \ \ \ \ \ \ &\textrm{if}\ \kap>0\\
\geq c_0|\xi|^2|\eta|^2  &\textrm{if}\ \kap<0
\end{array}
\qquad \forall \xi, \eta\in \R^n, \xi\perp \eta,
\end{equation}
with $c_0$ being a positive constant, is sufficient to obtain a priori $C^{1,1}$ bounds
for the smooth solutions.

It is noteworthy to point out that the condition (\ref{eq-A3})  and the $C^2$ estimates
 were derived in \cite{MTW} for the Monge-Amp\`ere type equations with variational structure emerging in optimal transport theory.
 The method used there is based on comparing the weak solution with the smooth one
 in a small ball. To employ this method successfully  in the outset of refractor problem we need to establish a comparison principle,
 suitable mollification of the weak solution and a priori estimated for the smooth solutions of Dirichlet's problem in small balls.

 The method outlined above gives the
 $C^2$ estimates for non-variational case as well, see \cite{KW-1,KW-2}.
 Therefore the local regularity result for the solutions to (\ref{M-A-short+})-(\ref{M-A-short-}) with
smooth $w$ will follow once the matrix $G$ verifies the condition (\ref{eq-A3}).
That done,  the regularity of weak solutions reduces to the verification of the inequality (\ref{eq-A3})
with some positive constant $c_0$.

\smallskip

The conditions imposed on the matrix in (\ref{M-A-short+})-(\ref{M-A-short-}) involving the Hessian
implies that the Monge-Amp\`ere equation is degenerate elliptic. The weak formulation of
degenerate ellipticity will be discussed in Section \ref{sec-A-type}. Postponing the precise definition of weak solutions until then
we would like to point out how the ellipticity of equation follows if we
consider those $C^2$ solutions of  (\ref{M-A-short+}) (resp.  (\ref{M-A-short-})) for
which at every point $x\in \U$ there is a hyperboloid (resp. ellipsoid) of revolution $H(\cdot, a, Z)$ touching
$u$ from above (reps. below) at $x$.
Indeed, for $H(x)=\ell_0-\frac ab\sqrt{b^2+|x-x_0|^2}$ the
matrix $\cal W_H=-\frac G{\eps\kap}+D^2H$ is identically zero. To see this we consider
the case of planar receiver $\Sigma$ given as $X^{n+1}=m$ with $m>0$. Without loss of generality
we take $x_0=0$.
Then $H(0)=\ell_0-a$. On the other hand it follows from the definition of
eccentricity $\eps=\frac{\sqrt{a^2+b^2}}{a}$ that $\ell_0=m-a\eps$, see Section \ref{sec-ell-hyp}. Next, a simple geometric
reasoning yields  the following explicit formula for the stretch function
\begin{equation}\label{t-comp-0}
 t=\frac{m-H}{Y^{n+1}}=\frac{c+\frac ab\sqrt{b^2+|x|^2}}{\eps(1-\frac{\kap}{q+1})}.
\end{equation}

We have $DH=-\frac ab\frac{x}{\sqrt{b^2+|x|^2}}$. Consequently
\begin{equation}\label{eq-Hess-H}
 D^2H=-\frac{a}{b\sqrt{b^2+|x|^2}}\left(\Id-\frac{x\otimes x}{b^2+|x|^2}\right).
\end{equation}

Moreover, recalling (\ref{def-sig})
we obtain $\kap=1-\frac1{\ee^2}=\frac{b^2}{c^2}$ where $c=\sqrt{a^2+b^2}$. This gives %
\begin{equation}\label{q-comp}
 q(x)=\frac1\ee\frac b{\sqrt{b^2+|x|^2}},
\end{equation}
in lieu of  (\ref{def-q}).

Thus combining these formulae for $t$ and $q$
we get from (\ref{def-G-M}), (\ref{t-comp-0}) and (\ref{eq-Hess-H})

\begin{eqnarray*}
\cal W_H&=&-\frac{G}{\eps\kap}-D^2H\\\nonumber
&=&-\frac{q+1}{\kap \eps  t}\left\{[\delta_{ij}-\kap \eps^2 H_iH_j]+ \frac{\kap\eps t D^2_{ij}H}{q+1}\right\} \\\nonumber %%%
&=&-\frac{q+1}{\kap \eps  t}\left\{  \Id-\bigg[\kap\eps^2\underbrace{\frac{a^2}{b^2}\bigg]\frac{x\otimes x}{b^2+|x|^2}}_{DH\otimes DH}+
  \frac{\kap\eps t}{q+1}   \underbrace{\left[-\frac{a}{b\sqrt{b^2+|x|^2}}\left(\Id-\frac{x\otimes x}{b^2+|x|^2}\right)\right]}_{D^2 H}\right\}.
\end{eqnarray*}
From the definition of $\kap$ (\ref{def-kap}) it follows
that $\kap\eps^2\frac{a^2}{b^2}=\frac{b^2}{c^2}\ee^2\frac{a^2}{b^2}=1$ implying

\begin{eqnarray*}
 \cal W_H=-\frac{q+1}{\kap \eps  t}\left(\Id-\frac{x\otimes x}{b^2+|x|^2}\right)
 \left\{1-
  \frac{\kap\eps t}{q+1} 
  \left[\frac{a}{b\sqrt{b^2+|x|^2}}\right]\right\}.
\end{eqnarray*}

Therefore,  recalling (\ref{t-comp-0}) and (\ref{q-comp}) we easily compute
\begin{eqnarray}\label{t-comp-1}
 t&=&\frac{c+\frac ab\sqrt{b^2+|x|^2}}{\eps(1-\frac{\kap}{h+1})}=
 (q+1)\frac{c+\frac ab\sqrt{b^2+|x|^2}}{\eps(q+1-\kap)}\\\nonumber
 &=&(q+1)\frac{c+\frac ab\sqrt{b^2+|x|^2}}{\eps(q+\frac1{\eps^2})}\\\nonumber
 &=&\eps(q+1)\frac{c+\frac ab\sqrt{b^2+|x|^2}}{\eps^2 q+1}\\\nonumber
  &=&\eps(q+1)\sqrt{b^2+|x|^2}\frac{c+\frac ab\sqrt{b^2+|x|^2}}{\eps b +\sqrt{b^2+|x|^2}}.\\\nonumber
\end{eqnarray}

Returning to $\cal W_H$ and utilizing  (\ref{t-comp-1}) we obtain
\begin{eqnarray*}
 1-  \frac{\kap\eps t}{q+1}
  \left[\frac{a}{b\sqrt{b^2+|x|^2}}\right]
  &=&
   1-{\kap}\ee^2 \frac{c+\frac ab\sqrt{b^2+|x|^2}}{\sqrt{b^2+|x|^2}+b\ee}\\
   &=&1-\eps\kap\frac{c}{a}\frac{a^2}{b^2}\\\nonumber
   &=&1-\eps^2\kap\frac{a^2}{b^2}\\\nonumber
   &=&0.
\end{eqnarray*}
A  similar  computation  for the matrix $\displaystyle\mathcal W_E=\frac{G}{t\eps|\kap|}+D^2 E$
can be carried out for the ellipsoids of revolution $E$ (i.e. for  $\eps<1, \kap<0$).
\smallskip

Since $-D^2u\geq -D^2H_{x_0}$ at $x_0$ and
$ \mathcal W_H=-\frac G{\eps\kap}-D^2H_{x_0}\equiv 0$
it follows that the equation $\det\left[-\frac G{\eps\kap}-D^2u\right]=h$ is degenerate elliptic.
\smallskip

Notice that for $\eps<1$ the weak solution
has a supporting ellipsoid of revolution $E_{x_0}$ at each point $x_0\in \overline{ \U}$ touching $\Gamma_u$
from below.  In particular we see that if $u\in C^2$ then
$Du=DE_{x_0}, -D^2u\leq -D^2 E_{x_0}$ at $x_0$.  Thus $\displaystyle  \frac G{\eps\kap}+D^2u\geq 0$
and we infer that (\ref{M-A-short+}) is degenerate elliptic. Analogously, using the hyperboloids as supporting functions, one can check that (\ref{M-A-short-}) is also degenerate elliptic.

%%%%%%%%%%%%%%%%%%%%%%%%%%%%%%%%%%%%%%%%%%%%%%%%%%%%%%%%%%%%
%                                                          %
%                SECTION                                   %
%                                                          %
%%%%%%%%%%%%%%%%%%%%%%%%%%%%%%%%%%%%%%%%%%%%%%%%%%%%%%%%%%%%
\section{Proof of Theorem B: Verifying the A3 condition}
In this section we explicitly compute the second derivatives in $p$ variable of the matrix $G^{ij}(x, u, p)$ introduced in 
\eqref{def-G-M} where $p$ is the dummy variable for $D u$.
We will find a concise representation of the form $D_{p_kp_l}G^{ij}\xi^i\xi^j\eta^k\eta^l$ for $ \xi, \eta\in \R^n, \xi\perp\eta$ and relate it with the
second fundamental form of the receiver $\Sigma=\{Z\in \R^{n+1}\  :\ \psi(Z)=0\}$ where $\psi:\R^{n+1}\to \R$ is a smooth function
such that  (\ref{vis-cond}) holds. Recall that the existence of smooth solutions 
depends on the sign of the form $D_{p_kp_l}G^{ij}\xi^i\xi^j\eta^k\eta^l$, 
see \eqref{eq-A3} and the discussion in the previous section.
%

%%%%%%%%%%%%%%%%%%%%%%%%%%%%%%%%%%%%%%%%%%%%%%%%%%%%%%%%%%%%%
%%%%%%%%%%%%%%%%%%%%%%%%%%%%%%%%%%%%%%%%%%%%%%%%%%%%%%%%%%%%%

\subsection{Computing the derivatives of stretch function $t$}
Recall that by (\ref{def-Z}) $Z(x)=x+e_{n+1}u(x)+tY$. Differentiating
$\psi(Z(x))=0$ with respect to $p_k$ we get
\begin{eqnarray}\label{t_k}
  \frac{t_{p_k}}{t}=-\frac{\sum_m\psi_m Y^m_{p_k}}{\sum_m\psi_m Y^m}, \quad k=1, \dots n.
\end{eqnarray}
After differentiating again by $p_l$ we get
\begin{eqnarray}
  \frac{t_{p_kp_l}}{t}-\frac{t_{p_k}t_{p_l}}{t^2}
  &=&-\Bigg[\frac{\sum_{ms}\psi_{ms}({}Y^s_{p_l}t+{}Y^st_{p_l}){}Y^m_{p_k}+\sum\psi_m{{}Y^m_{p_kp_l}}}{(\nabla\psi\cdot{}Y)}\\\nonumber
  &&-\frac{\sum_m\psi_m{}Y^m_{p_k}}{(\nabla\psi\cdot{}Y)^2}
  \left(\sum_{m,s}\psi_{ms}({}Y^s_{p_l}t+{}Y^st_{p_l}){}Y^m+\sum_m\psi_m{}Y^m_{p_l}\right)\Bigg]\\\nonumber
  %%%%
  &=&-\frac{1}{(\nabla\psi\cdot{}Y )}\bigg[
  \left(\nabla^2\psi{}Y_{p_k}{}Y_{p_l}-\frac{\nabla\psi{}Y_{p_k}}{(\nabla\psi\cdot{}Y)}\nabla^2\psi{}Y{}Y_{p_l}\right)t\\\nonumber
  &&+\left(\nabla^2\psi{}Y_{p_k}{}Y-\frac{\nabla\psi{}Y_{p_k}}{(\nabla\psi\cdot{}Y)}\nabla^2\psi{}Y{}Y\right)t_{p_l}\\\nonumber
  &&+\nabla\psi{}Y_{p_kp_l}-\frac{\nabla\psi{}Y_{p_k}}{(\nabla\psi\cdot{}Y)}\nabla\psi{}Y_{p_l}\bigg]=\\\nonumber
%%%%%
  &=&-\frac{1}{(\nabla\psi\cdot{}Y) }\bigg[
  \left(\nabla^2\psi{}Y_{p_k}{}Y_{p_l}-\frac{\nabla\psi{}Y_{p_k}}{(\nabla\psi\cdot{}Y)}\nabla^2\psi{}Y{}Y_{p_l}\right)t\\\nonumber
  &&+\left(\nabla^2\psi{}Y_{p_k}{}Y-\frac{\nabla\psi{}Y_{p_k}}{(\nabla\psi \cdot {}Y)}\nabla^2\psi{}Y{}Y\right)t_{p_l}
  +\nabla\psi{}Y_{p_kp_l}\bigg]\\\nonumber
  &&+\frac{t_{p_k}t_{p_l}}{t^2}.
\end{eqnarray}
Rearranging the terms we infer
\begin{eqnarray}\label{t_2nd_der}
D^2_{p_kp_l}\left(\frac1t\right)=-\frac{t_{p_kp_l}}{t^2}+\frac{2t_{p_k}t_{p_l}}{t^3}&=&\frac{1}{t(\nabla\psi\cdot{}Y)
}\bigg[
  \left(\nabla^2\psi{}Y_{p_k}{}Y_{p_l}-\frac{\nabla\psi{}Y_{p_k}}{(\nabla\psi\cdot{}Y)}\nabla^2\psi{}Y{}Y_{p_l}\right)t\\\nonumber
  &&+\left(\nabla^2\psi{}Y_{p_k}{}Y-\frac{\nabla\psi{}Y_{p_k}}{\nabla\psi{}\cdot Y}\nabla^2\psi{}Y{}Y\right)t_{p_l}
  +\nabla\psi{}Y_{p_kp_l}\bigg]\\\nonumber
  &=&\frac{1}{t(\nabla\psi\cdot{}Y)
}\bigg[\frac{1}{t}\nabla^2\psi
Z_{p_k}Z_{p_l}+\nabla\psi{}Y_{p_kp_l}\bigg]
\end{eqnarray}
where the last line follows from  (\ref{t_k}).
Thus the second derivatives of $\frac1t$ can be computed from (\ref{t_2nd_der}), while for the first order
derivatives we have the formula (\ref{t_k}).

\smallskip

Next,  we want to compute the derivatives of $M_{ij}=(q+1)[\delta_{ij}-\kap\eps^2p_ip_j]$
with respect to $p$.
We have
\begin{eqnarray*}
 D_{p_k}M_{ij}=q_{p_k}(\delta_{ij}-\kap\eps^2 p_ip_j)-(q+1)\kap\eps^2[\delta_{kj}p_i+\delta_{ki}p_j].
\end{eqnarray*}
The condition $\xi\perp \eta $ implies that the contribution of the terms involving $\delta_{kj}$ and $\delta_{ki}$ is zero.
Thus we infer
$$D_{p_kp_l}^2M_{ij}\xi^i\xi^j\eta^k\eta^l=q_{p_kp_l}\eta^k\eta^l\left[|\xi|^2-\kap\eps^2(p\cdot \xi)^2\right].
$$
Recall that by definition $G=\frac Mt$ hence from the product rule we have
\begin{eqnarray}\label{A3-1}
 D^2_{p_kp_l}G^{ij}\xi^i\xi^j\eta^k\eta^l&=& D_{p_kp_l}^2\left(\frac1t\right)\eta^k\eta^l(q+1)\left[|\xi|^2-\kap\eps^2(p\cdot \xi)^2\right]
\\\nonumber&&+2\left(D_p\frac1t\cdot\eta\right)(D_p q\cdot\eta)\left[|\xi|^2-\kap\eps^2(p\cdot \xi)^2\right]\\\nonumber
&&+\frac1t(D^2_{p_kp_l}q )\eta^k\eta^l\left[|\xi|^2-\kap\eps^2(p\cdot \xi)^2\right]\\\nonumber
&=& S_{kl}\eta^k\eta^l\left[|\xi|^2-\kap\eps^2(p\cdot \xi)^2\right],
\end{eqnarray}
where
\begin{eqnarray}
 S_{kl}=
 (q+1)D_{p_kp_l}^2\left(\frac1t\right)%\eta^k\eta^l
%\\\nonumber&&
+ D_{p_k}\left(\frac1t \right) D_{p_l} q + D_{p_l}\left(\frac1t \right) D_{p_k} q %\\\nonumber&&
+\frac{D^2_{p_kp_l}q }{t}.%\eta^k\eta^l.\\\nonumber
\end{eqnarray}
%and the first and second order derivatives of $\frac1t$ are given by (\ref{t_k}) and (\ref{t_2nd_der}).

It follows from (\ref{unit-Y-q}) that
\begin{equation}
 Y_{p_kp_l}(q+1)+Y_{p_k}q_{p_l}+ Y_{p_l}q_{p_k} +Yq_{p_kp_l}=e_{n+1} \eps q_{p_kp_l}.
\end{equation}
which after taking the inner product with $\nabla \psi$ and  dividing the
by $\nabla \psi \cdot Y$ yields
\begin{eqnarray}
 \frac{(q+1)\psi_{n+1}q_{p_kp_l}}{\nabla \psi\cdot Y}&=&\frac{(q+1)\nabla \psi \cdot Y_{p_kp_l}}{\nabla \psi \cdot Y}
 +\frac{q_{p_k}\nabla\psi\cdot Y_{p_l}}{\nabla\psi\cdot Y}+\frac{q_{p_l}\nabla \psi \cdot Y_{p_k}}{\nabla \psi \cdot Y}
 +q_{p_kp_l}=\\\nonumber
 &=&\frac{(q+1)\nabla \psi \cdot Y_{p_kp_l}}{\nabla \psi \cdot Y}
 +\frac{q_{p_k}\nabla\psi\cdot Y_{p_l}}{\nabla\psi\cdot Y}+\frac{q_{p_l}\nabla \psi \cdot Y_{p_k}}{\nabla \psi \cdot Y}
 +q_{p_kp_l}\\\nonumber
 &=&\frac{(q+1)\nabla \psi \cdot Y_{p_kp_l}}{\nabla \psi \cdot Y}+tq_{p_k}D_{p_l}\left(\frac 1t \right)+
 tq_{p_l}D_{p_k}\left(\frac 1t \right)+q_{p_kp_l}\\\nonumber
 &=&\frac{(q+1)\nabla \psi \cdot Y_{p_kp_l}}{\nabla \psi \cdot Y}+t\left(S_{kl}-
 (q+1)D_{p_kp_l}^2\left(\frac1t\right)\right).
\end{eqnarray}

Consequently, with the aid of (\ref{t_2nd_der}) we find  that
\begin{eqnarray*}
 S_{kl}&=&\frac{q+1}{t(\nabla \psi \cdot Y)}
 \left( {\psi_{n+1}q_{p_kp_l}} -
 {\nabla \psi \cdot Y_{p_kp_l}}\right)+ (q+1)D_{p_kp_l}^2\left(\frac1t\right)\\\nonumber
 &=&\frac{q+1}{t(\nabla \psi \cdot Y)}
 \left( {\psi_{n+1}q_{p_kp_l}} -
 {\nabla \psi \cdot Y_{p_kp_l}}\right)+ \frac{q+1}{t(\nabla\psi\cdot{}Y)
}\bigg[\frac{1}{t}\nabla^2\psi
Z_{p_k}Z_{p_l}+\nabla\psi{}Y_{p_kp_l}\bigg]\\\nonumber
&=&\frac{q+1}{t(\nabla \psi \cdot Y)}
 \left[\frac{1}{t}\nabla^2\psi
Z_{p_k}Z_{p_l}+{\psi_{n+1}q_{p_kp_l}}\right].
\end{eqnarray*}
It remains to recall that by (\ref{def-q})
\begin{equation}\label{hess-q}
q_{p_kp_l}=-\frac{\kap}{q}\left[\delta_{kl}+\kap\frac{p_kp_l}{q^2}\right]
\end{equation}
and we conclude

\begin{eqnarray}\label{A3-111}
 D^2_{p_kp_l}G^{ij}\xi^i\xi^j\eta^k\eta^l=\frac{q+1}{t(\nabla \psi \cdot Y)}
 \left[\frac{1}{t}\nabla^2\psi
Z_{p_k}Z_{p_l}+{\psi_{n+1}q_{p_kp_l}}\right]\eta^k\eta^l\left[|\xi|^2-\kap\eps^2(p\cdot \xi)^2\right].
\end{eqnarray}

It is worth noting that $|\xi|^2-\kap\eps^2(p\cdot \xi)^2$ is always positive. This is obvious if  $\kap<0$. As for
$\kap >0$ then we note that $|\xi|^2-\kap\eps^2(p\cdot \xi)^2=|\xi|^2\left(1-\kap\eps^2\left(p\cdot\frac\xi{|\xi|}\right)^2\right)>0$
in view of the estimate $|p|< \frac1{\sqrt{\eps^2-1}}$, see Lemma \ref{lem-lips}.
Furthermore, from (\ref{vis-cond}) it follows that $D^{2}_{p_kp_l}G^{ij}$ and $\widehat S_{kl}$ defined by
\begin{equation}\label{S-form}
\widehat S_{kl}=\frac{1}{t}\nabla^2\psi
Z_{p_k}Z_{p_l}+{\psi_{n+1}q_{p_kp_l}}\
\end{equation}
have the same signs. Thus it is enough to explore the form $\widehat S_{kl}$ instead.
%%%
%%%
\subsection{Refining condition (\ref{eq-A3})}
Let $Z_0$ be a fixed point on $\Sigma$.
Introduce a new coordinate
system $\hat x_1, \dots, \hat x_n,
\hat x_{n+1}$ near $Z_0$, with $\hat x_{n+1}$  having direction $Y$.
Since (\ref{vis-cond}) and (\ref{dist-cond}) implies  $\nabla \psi\ne 0$, without loss of
generality we assume that near $Z_0$, in $\hat x_1, \dots, \hat x_n,
\hat x_{n+1}$ coordinate system $\Sigma$ has a representation
$\hat x_{n+1}=\phi(\hat x_1, \dots, \hat x_n).$
Recall that the second fundamental form of $\Sigma $ is
\begin{eqnarray}\label{second-form}
\mbox{II}=\frac{\p^2_{\hat x_i, \hat x_j}\phi}{\sqrt{1+|D\phi|^2}} , \qquad i,j=1,\dots, n
\end{eqnarray}
if we choose the normal of $\Sigma$ at $Z_0$ to be
$\frac{(-D_{\hat x_1}\phi, \dots, -D_{\hat x_n}\phi, 1)}{\sqrt{1+|D\phi|^2}}, D\phi=(D_{\hat x_1}\phi, \dots, D_{\hat x_n}\phi,0)$.

Denote  $\widetilde \psi(Z)=Z^{n+1}-\phi(z)$ and assume that
near $Z_0$, $\Sigma $ is given by the equation $\widetilde \psi=0$. It follows that
\begin{eqnarray}\label{Hess-psi}
 \nabla^2\widetilde \psi= -\left|\begin{array}{cccc}
        \phi_{11} &\cdots & \phi_{1n} & 0\\
         \vdots  &\ddots &\vdots &\vdots \\
        \phi_{n1} &\cdots & \phi_{nn} & 0\\
        0 &\cdots & 0 & 0\\
         \end{array} \right|.
\end{eqnarray}

Therefore for  $Z=x+u e_{n+1}+tY$ we have $\nabla^2\widetilde \psi  Y=0$ and hence
\begin{eqnarray}\label{2nd-f}
  \nabla^2\widetilde \psi Z_{p_k}Z_{p_l}&=&\nabla^2\widetilde\psi(tY_{p_k}+t_{p_k}Y)(tY_{p_l}+t_{p_l}Y)\\\nonumber
  &=&t^2\nabla^2\widetilde\psi Y_{p_k}Y_{p_l}\\\nonumber
  &=&-t^2\nabla^2\phi  Y_{p_k}Y_{p_l}.\\\nonumber
\end{eqnarray}

By (\ref{eq-y-sig}) $Y(q+1)=(\eps \kap p, \eps q+\eps -\kap)$  where $y(q+1)=\eps \kap p$.
Differentiating this equality with respect to  $p_k$ we infer

\begin{equation}
 Y_{p_k}(q+1)+Yq_{p_k}=\eps(\kap \hat e_k+q_{p_k}\hat e_{n+1})
\end{equation}
hence 
\begin{equation}
 Y_{p_k}=\frac1{q+1}\left[-Yq_{p_k}+\eps(\kap \hat e_k+q_{p_k}\hat e_{n+1})\right].
\end{equation}
On the other hand (\ref{Hess-psi}) and $\hat e_{n+1}= Y$ yield
\begin{eqnarray}
 \nabla^2\widetilde \psi Y_{p_k}&=&\frac1{q+1}
 \nabla^2\widetilde \psi \left[-Yq_{p_k}+\eps(\kap \hat e_k+q_{p_k}\hat e_{n+1})\right]\\\nonumber
 &=&\frac\eps{q+1}
 \nabla^2\widetilde \psi (\kap \hat e_k+q_{p_k}\hat e_{n+1})\\\nonumber
 &=&\frac{\eps\kap}{q+1}
 \nabla^2\widetilde \psi \hat e_k.\\\nonumber
\end{eqnarray}
Since $\nabla^2\widetilde \psi$ is symmetric we infer
\begin{eqnarray}\label{Y-sclr}
 \nabla^2\widetilde \psi Y_{p_k}Y_{p_l}&=&\frac{\eps^2\kap^2}{(q+1)^2}
 \nabla^2\widetilde \psi \hat e_k\hat e_l.\\\nonumber
\end{eqnarray}
Plugging (\ref{Y-sclr}) into (\ref{2nd-f}) we finally obtain
\begin{eqnarray}
 \nabla^2\widetilde \psi Z_{p_k}Z_{p_l}&=&-t^2\frac{\eps^2\kap^2}{(q+1)^2}
 \nabla^2 \phi \hat e_k\hat e_l.\\\nonumber
 &=&-t^2\frac{\eps^2\kap^2}{(q+1)^2}\sqrt{1+|D\phi|^2}\hbox{II}
\end{eqnarray}
where $\hbox{II}$ is the second fundamental form of $\Sigma$ at $Z_0$, see (\ref{second-form}).
%%% \td{double check this condition}
This in conjunction with (\ref{hess-q}) yields
\begin{equation}
\widehat S_{kl}=-\left[ \frac1t\left(\frac{t\eps\kap}{q+1}\right)^2 \sqrt{1+|D\phi|^2}{\rm II}+\frac{%\widetilde\psi_{n+1}
\kap}q (\Id+\kap \frac{p\otimes p}{q^2})
\right].
\end{equation}
Summarizing 
\begin{eqnarray}\label{A3-1111}
\\\nonumber
 D^2_{p_kp_l}G^{ij}\xi^i\xi^j\eta^k\eta^l=-\frac{q+1}{t(\nabla \psi \cdot Y)}
 \left[ \frac1t\left(\frac{t\eps\kap}{q+1}\right)^2 \sqrt{1+|D\phi|^2}{\rm II}+\frac{%\widetilde\psi_{n+1}
\kap}q (\Id+\kap \frac{p\otimes p}{q^2})
\right]
 \eta^k\eta^l\left[|\xi|^2-\kap\eps^2(p\cdot \xi)^2\right].
\end{eqnarray}
\subsection{Examples of receivers}\label{ssec-example}
Let us consider the case of horizontal receiver $Z^{n+1}=m>0$ for some positive number $m$.
Then ${\rm II}=0$ implying that 
$$D^2_{p_kp_l}G^{ij}\xi^i\xi^j\eta^k\eta^l=-\frac{q+1}{t(\nabla \psi \cdot Y)}\frac{%\widetilde\psi_{n+1}
\kap}q \left(\Id+\kap \frac{p\otimes p}{q^2}\right)
\eta^k\eta^l\left[|\xi|^2-\kap\eps^2(p\cdot \xi)^2\right].$$ 
If $\kap>0$ then in view of \eqref{blya-1} we have $|\xi|^2-\kap\eps^2(p\cdot \xi)^2>0$ and 
clearly $D^2_{p_kp_l}G^{ij}\xi^i\xi^j\eta^k\eta^l<0$. Hence
\eqref{eq-A3} is not satisfied for this case. 

\smallskip 

As for $\kap<0$ we compute
\begin{equation}\label{blya-7}
-\frac{\kappa}{q}(|\eta|^2+\kappa\frac{(p\cdot\eta)^2}{q^2})=\frac{|\kap| |\eta|^2}q(1-\frac{|\kap| (p\cdot \eta)^2}{q^2|\eta|^2})\ge
\frac{|\kap| |\eta|^2}q (1-\frac{|\kap |p|^2}{q^2})=\frac{|\kap|(1+|\kap|)}{q^3}|\eta|^2 \geq c_0|\eta|^2
\end{equation}
where $c_0>0$ depends only on $\sup |p|$ and $\eps$. Consequently, $D^2_{p_kp_l}G^{ij}\xi^i\xi^j\eta^k\eta^l>0$ and 
(\ref{eq-A3}) is not true for
horizontal receivers $Z^{n+1}=m>0$.

Next we construct examples of $\psi$ {\it satisfying} \eqref{eq-A3}.

{\bf Case\ 1:} $\eps>1, \kappa>0$.

As we have seen above $|\xi|^2-\kap\eps^2(p\cdot \xi)^2>0$. Therefore in order to get 
$D^2_{p_kp_l}G^{ij}\xi^i\xi^j\eta^k\eta^l>0$ we must assume that ${\rm II}<-\delta\Id$
for some constant $\delta>0$.  It is enough to show that 
$$J:=\frac1t\left(\frac{t\eps\kap}{q+1}\right)^2 \sqrt{1+|D\phi|^2}{\rm II}+\frac{\kap}q (\Id+\kap \frac{p\otimes p}{q^2})<0.$$
We have 
\begin{eqnarray}
\frac{\kap}q (|\eta|^2+\kap \frac{(p\cdot\eta)^2}{q^2})\leq |\eta|^2\frac{\kap}q (1+\kap\frac{|p|^2}{q^2} )=|\eta|^2\frac{\kap}{q^3}(1-\kappa).
\end{eqnarray}
Therefore using the fact that $q<1$ for $\kap>0$ we get 
\begin{eqnarray}
J&\leq& -\delta t\left(\frac{\eps\kap}{q+1}\right)^2\Id + \frac{\kap}{q^3}(1-\kappa)\Id\\\nonumber
&=& -\kap\Id \left(t\delta \frac{\eps^2\kap}{(1+q)^2}-\frac1{\eps^2}\right)\\\nonumber
&\le& -\kap\Id\left(t\delta \frac{\eps^2\kap}{4}-\frac1{\eps^2}\right)\\\nonumber
&<0&
\end{eqnarray}
if $t>\frac4{\delta\eps^4\kap}$. Thus if $\Sigma$ is strictly concave and is sufficiently far from 
$\U$ then \eqref{eq-A3} is satisfied for $\eps>1, \kap>0$.

\smallskip

{\bf Case\ 2:} $\eps<1, \kappa<0$.

Again, we obviously have  $|\xi|^2-\kap\eps^2(p\cdot \xi)^2>0$. Therefore we demand $J>0$.
For  this we suppose that $\rm{II}>\delta\Id$ for some $\delta>0$. Then by \eqref{blya-7}
\begin{eqnarray}
J&\ge&\delta t \left(\frac{\eps\kap}{q+1}\right)^2\Id-\frac{|\kap|}q \Id=\\\nonumber
&=&|\kap|\Id \left(\delta t\frac{\eps^2|\kap|}{(1+q)^2}-\frac1q\right) \quad (\text{recall that}\ q=\sqrt{1+|\kap|(1+|p|^2)}>1)\\\nonumber
&=&|\kap|\Id \left(\delta t\frac{\eps^2|\kap|}{(1+q)^2}-1\right)\\\nonumber
&\ge& |\kap|\Id \left(\delta t\frac{\eps^2|\kap|}{(1+\sqrt{1+|\kap|(1+\mu^2)})^2}-1\right)\\\nonumber
&>0&
\end{eqnarray}
if $t>\frac{(1+\sqrt{1+|\kap|(1+\mu^2)})^2}{\delta|\kap|\eps^2}$ and $|p|\leq \mu$, where $\mu$ depends on the 
Lipschitz constant of the solution, see \eqref{blya-1}. 
Thus if $\Sigma$ is strictly convex and is sufficiently far from 
$\U$ then \eqref{eq-A3} is satisfied for $\eps<1, \kap<0$.

\begin{remark}\label{rem-extns}
We summarize the following conclusions based on the discussion in this section:  
\begin{itemize}
\item 
The computation above shows that (\ref{eq-A3}) is true if $\kap>0, {\rm II}>0$ or if $\kap<0, {\rm II}<0$
and $\dist(\U, \V)$ is large enough.
\item
If $k>0$ and $\Sigma$ is a graph of strictly concave function $Z^{n+1}=\phi(z)$ then \eqref{vis-cond} holds.
Hence the refracted rays intersect $\Sigma$ at a unique point.
\item
We can extend $\Sigma$ to entire space such that the resultd surface is still concave if say $\kap>0$, 
hence without loss of generality we can assume that $\Sigma$ is entire concave 
surface and so is $\Sigma+Me_{n+1}$, for $M\gg1$. We will take advantage of this in 
Lemmas \ref{lem-conf} and \ref{lem-apprx}
\end{itemize}
 \end{remark}

%%%%%%%%%%%%%%%%%%%%%%%%%%%%%%%%%%%%%%%%%%%%%%%%%%%%%%%%%%%
%%%%%%%%%%%%%%%%%%%%%%%%%%%%%%%%%%%%%%%%%%%%%%%%%%%%%%%%%%%

%%%%%%%%%%%%%%%%%%%%%%%%%%%%%%%%%%%%%%%%%%%%%%%%%%%%%%%%%%%%
%                                                          %
%                SECTION                                   %
%                                                          %
%%%%%%%%%%%%%%%%%%%%%%%%%%%%%%%%%%%%%%%%%%%%%%%%%%%%%%%%%%%%
\section{Admissible functions}\label{sec-ell-hyp}

The refractive properties of ellipses and hyperbolas have been known since ancient times \cite{Mountford}.
Furthermore, hyperboloids and ellipsoids of revolution share the same properties. This section is devoted to
the class of functions obtained as envelopes of halves of ellipsoids and hyperboloids of revolution.

\smallskip
%%%%%%%%%%%%%%%%%%%%%%%%%%%%%%%%%%%%%%%%%%%%%%%%%%%%%%%%%%%%%%%%%%%%%%%%%%%%%%%%%%%%%
%%%%%%%%%%%%%%%%%%%%%%%%%%%%%%%%%%%%%%%%%%%%%%%%%%%%%%%%%%%%%%%%%%%%%%%%%%%%%%%%%%%%%
\begin{figure}
  \centering
  \subfloat[Upper admissible]{\includegraphics[width=0.4\textwidth]{hyper.eps}}
  \subfloat[Lower  admissible]{
\includegraphics[width=0.4\textwidth]{ellipse.eps}}
  \caption{Supporting surfaces}
  \label{fig:0-pic}
\end{figure}

\subsection{Ellipsoids}
Throughout this paper by ellipsoid  we mean the lower half of an ellipsoid of revolution
with focal axis parallel to $e_{n+1}$. Such surface can be regarded as the graph of
\begin{equation}\label{def-ell}
 E(x, a, Z)=Z^{n+1}-a\eps -a\sqrt{1-\frac{(x-z)^2}{a^2(1-\eps^2)}}
\end{equation}
where $a$ is the larger semiaxis, $\eps$- the eccentricity,  and $Z$ the higher focus, see Figure 2.
Moreover we have that
\begin{equation}\label{ell-grad}
 DE=\frac1{a(1-\eps^2)}\frac{x-z}{\sqrt{1-\frac{(x-z)^2}{a^2(1-\eps^2)}}}.
\end{equation}
Notice that at the points $x$ where $|x-z|=a\sqrt{1-\eps^2}$ the gradient
$|DE|$ is unbounded.

\smallskip

\subsection{Hyperboloids}
It is convenient to introduce the lower sheet of hyperboloids of revolution
\begin{equation}\label{def-hypr}
 H(x, a, Z)=Z^{n+1}-a\eps -a \sqrt{1+\frac{(x-z)^2}{a^2(\eps^2-1)}}
\end{equation}
where $a$ is the larger semiaxis, $\eps$ the eccentricity, and  $Z$ the upper focus, see Figure 2.
Differentiating $H$ we obtain
\begin{equation}\label{hypr-grad}
 DH=-\frac1{(\eps^2-1)}\frac{x-z}{\sqrt{a^2+\frac{(x-z)^2}{\eps^2-1}}}.
\end{equation}

\smallskip

\subsection{Supporting hyperboloids}
\begin{definition}\label{def-adms}
A function $u:\U\to \R$ is said to be upper (resp. lower) admissible if for any
$x_0\in \U$ there is $Z\in \V$ and $a>0$ such that $H(x_0, a, Z)=u(x_0)$ (resp. $E(x_0, a, Z)=u(x_0)$) and $H(x, a, Z)\ge u(x), x\in \U$
(resp. $E(x, a, Z)=u(x)$).
$H$ (resp. $E$) is called a supporting function of $u$ at $x_0$. The class of all upper admissible functions is denoted by $\Wup(\U, \V)$ (resp. $\Wlo(\U, \V)$).
\end{definition}
In what follows we focus on upper admissible functions, the lower admissible functions can be studied in similar fashion.
If the generalisation is not straightforward then we will outline the proof.
\smallskip

Formula  (\ref{hypr-grad}) yields uniform Lipschitz estimates for $\Wup(\U, \V)$.  
\begin{lemma}\label{lem-lips}
Let $\mathbb H^+_{a_0}(\U, \V)$ be the set of all hyperboloids $H(x, a, Z)\ge 0, x\in \U$ for some $Z\in \V$ 
such that  $a>a_0$ for some fixed $a_0\ge 0$.
Then 
\begin{equation}\label{blya-1}
\sup\limits_{\mathbb H^+_{a_0}(\U, \V)}\|D H\|_{\infty}\le \frac1{\sqrt{\eps^2-1}}\frac{d_0}{\sqrt{a_0^2(\eps^2-1)+d_0^2}} \quad (if  \ a_0>0)
\end{equation}
where $d_0=\sup_{x\in\U, z\in\widehat\V}|x-z|.$
Furthermore, 
$$\|D u\|_{L^\infty(\U')}<\frac1{\sqrt{\eps^2-1}}, \quad \forall \U'\subset\subset \U, u\in {\Wup (\U, \V)}.$$
\end{lemma}
\pr From (\ref{hypr-grad}) we have for $x\not=z$
\begin{eqnarray}
|DH|&=&\frac1{(\eps^2-1)}\frac{|x-z|}{\sqrt{a^2+\frac{|x-z|^2}{\eps^2-1}}}\\\nonumber
&=&\frac1{\sqrt{\eps^2-1}}\frac{|x-z|}{\sqrt{a^2(\eps^2-1)+|x-z|^2}}=\frac1{\sqrt{\eps^2-1}}\frac{1}{\sqrt{\frac{a^2(\eps^2-1)}{|x-z|^2}+1}}\\\nonumber
&\leq&\frac1{\sqrt{\eps^2-1}}\frac{1}{\sqrt{\frac{a^2(\eps^2-1)}{\sup_{x\in\U, z\in\widehat\V}|x-z|^2}+1}}.
\end{eqnarray}
If  $a\geq a_0>0$  then the first inequality immediately follows. 

In order to prove the second inequality let us suppose that for some fixed 
subdomain $\U'\subset \subset \U$ there are $p_k\in \p u(x_k)$ such that $x_k \in \U'$, 
$p_k\to p_0, x_k \to x_0$ and $|p_0|=\frac1{\sqrt{\eps^2-1}}$. Here $\p u(x)$ is the subdifferential of $u$ at $x$. It is clear that 
$p_0\in \p u(x_0)$ and hence there is a supporting hyperplane for $u$ at $x_0$ with slope $p_0$. If 
$u$ is strictly concave at $x_0\in U$ then near $x_0$ one can find $\bar x_0$ such that there is $\bar p_0\in \p u(\bar x_0)$
with $|\bar p_0|>\frac1{\sqrt{\eps^2-1}}$ which is in contradiction with the first inequality. Thus suppose that 
there is a straight segment in the graph of $u$ passing through $(x_0, u(x_0))$. But this is impossible 
because $u$ is admissible and therefore $\Gamma_u$ cannot contain straight segments. 
\qed

\smallskip
\begin{lemma}\label{lem-support}
Let $\{u_k\}$ be a sequence of upper admissible function such that
$u_k\to u_0$ uniformly in $\U$. If $x_k\in \U$, $x_k\to x_0$ and $H_k$ are supporting functions of $u_k$ at $x_k$ then
$u_0$ has an upper supporting function $H_0$ at $x_0$ and $H_k\to H_0$ uniformly in $\U$.
\end{lemma}

\pr 
One way to check the claim is to use some well known fact from convex analysis.
Consider the convex sets $G_k=\{X\in \R^{n+1} : x\in \U, 0< u_k(x)<X^{n+1}\}$ and 
$\mathscr H_k=\{X\in \R^{n+1} : x\in \U, 0< H_k(x)<X^{n+1}\}$ where $x=\widehat X$. 
Then $\mathscr H_k\subset G_k$ and $(x_k, u_k(x_k))\in \p\mathscr H_k\cap  \p G_k$. 
Thus, from uniform convergence $u_k\to u_0$ we infer that 
 the limit set $G_0=\{X\in \R^{n+1} : x\in \U, 0< u_0(x)<X^{n+1}\}$ is a subset of 
$\mathscr H_0=\{X\in \R^{n+1} : x\in \U, 0< H_0(x)<X^{n+1}\}$, see \cite{Alek-Book} Chapter 5.2.
Furthermore, from $x_k\to x_0\in \ol \U$ it follows that there is $X_0\in \p G_0\cap \p \mathscr H_0$ such that 
$\widehat X_0=x_0$.
Therefore we conclude that $H_0$ is a supporting hyperboloid  of $u_0$ at $x_0$.
\qed

\subsection{Continuous expansion of hyperboloids}\label{ssec-confocal}
If $u\in \Wup(\U, \Sigma)$ then it turns out that $u$ is also admissible with respect with 
$\widetilde \Sigma$, the receiver moved vertically upwards in $e_{n+1}$ direction.
In other words, the same admissible $u$ will be $R-$convex with respect to a family of 
surfaces obtained from $\Sigma$ by translation is $e_{n+1}$ direction. We will need this 
observation in order to construct smooth solutions of our problem in small balls, see Section \ref{sec-pr-Th2}.
\begin{lemma}\label{lem-conf}
Let $\widetilde \Sigma=\Sigma +M e_{n+1}$ for some $M>0$.  
\begin{itemize}
 \item[(i)] For any fixed $x_0$ and $H_1(x)=H(x, a_1, Z_1)\in \mathbb H (\U, \Sigma)$ there is
 $H_2(x)=H(x, a_2, Z_2)$ with $Z_2\in \widetilde \Sigma$ and touching $H_1$ from above at $x_0$.
 \item[(ii)] In particular if $u\in \Wup(\U, \Sigma)$ then also $u\in \Wup(\U, \widetilde \Sigma)$.
 \end{itemize}
\end{lemma}
\pr (i) Let $\xi_1=H_1(x_0)$ and $X_0=(x_0, \xi_1)$. For $s>1$ we consider $Z_2=X_0+s(Z_1-X_0)$.
By construction $X_0, Z_1$ and $Z_2$ lie on the same line. To determine $a_2$ we utilize two geometric properties of hyperbola, namely that the difference of distances of $X_0$ from 
$Z_2$ and the lower focus $Z_2'$ is $2a_2$ and $|X_0Z_2'|=\eps |X_0D|$ where 
$|X_0D|$ is the distance of $X_0$ from the lower directrix $X^{n+1}=Z^{n+1}-a_2\eps-a_2/\eps$.
Therefore if $P$ is on the graph of $H_2$ we get that $|PZ_2|=2a_2+|PZ_2'|=-a_2(\eps^2-1)+s\eps(Z_1^{n+1}-\xi_1)$. Taking 
$P=X_0$ in this equation $|PZ_2|=\eps|Z_1-X_0|$ one finds that

\begin{equation}\label{a-2-eq}
a_2=\frac{1}{\eps^2-1}[s\eps(Z_1^{n+1}-\xi_1)-|s(Z_1-X_0)|].
\end{equation}

\smallskip

As for (ii), we choose $s_0>1$ so that $X_0+s(Z_1-X_0)\in \widetilde \Sigma$. Consequently from 
(i) it follows that $Z_2=X_0+s_0(Z_1-X_0)$ is the focus of supporting hyperboloid $H(\cdot, a_2, Z_2)$
at $x_0$ where $a_2$ is given by (\ref{a-2-eq}). Therefore $u\in \Wup(\U, \widetilde\Sigma)$.
\qed

\smallskip

%%%%%%%%%%%%%%%%%%%%%%%%%%%%%%%%%%%%%%%%%%%%%%%%%%%%%%%%%%%
%%%%%%%%%%%%%%%%%%%%%%%%%%%%%%%%%%%%%%%%%%%%%%%%%%%%%%%%%%%

%%%%%%%%%%%%%%%%%%%%%%%%%%%%%%%%%%%%%%%%%%%%%%%%%%%%%%%%%%%%
%                                                          %
%                SECTION                                   %
%                                                          %
%%%%%%%%%%%%%%%%%%%%%%%%%%%%%%%%%%%%%%%%%%%%%%%%%%%%%%%%%%%%
\section{$B$-type weak solutions: Proof of Theorem C1}\label{sec-weak-B}

In this section we introduce our first notion of weak solution for the refractor problem (\ref{problem}).
For any upper admissible function  $u\in \Wup(\U, \V)$ we define the mapping $\s_u: \V\to \U$ as follows
\begin{eqnarray*}
 \s_u(Z)&=&\{x\in \U:\ \exists\ \text{a supporting hyperboloid of $u$
 at $x$ with focus at $Z\in\V$}\}.
\end{eqnarray*}
For any Borel set $\omega\subset \V$ we put
\begin{eqnarray}
 \s_u(\omega)=\bigcup_{Z\in \omega}\s_u(Z).
\end{eqnarray}
We will write $\s(E)$ instead of
$\s_u(E)$ if there is no confusion.

\smallskip

\begin{proposition}\label{prop-S}
For  $u\in \Wup(\U, \V)$ the corresponding mapping $\s$ enjoys the following properties:
\begin{itemize}
 \item[a)]  $\s: \V\longrightarrow \Pi$ maps the closed sets to closed sets.
\item[b)] The mapping $\s$ is one-to-one modulo a set of vanishing measure, i.e.
$$\big|\big\{x\in \Pi: x\in\s(Z_1)\cap\s(Z_2)\  \textrm{for}\ Z_1\not= Z_2, \ \ Z_i\in \V, i=1,2\big\}\big|=0.$$
\item[c)] The family $\mathscr F=\{E\subset \V \ \textrm{such that} \ \s(E) \ \textrm{is measurable}\}$
is $\sigma-$algebra.
\end{itemize}
\end{proposition}

\pr
 The first claim a) follows directly from Lemma \ref{lem-support}.

In order to prove b) we set $A=\big\{x\in \Pi: x\in\s(Z_1)\cap\s(Z_2)\  \textrm{for}\ Z_1\not= Z_2, \ \ Z_i\in \V, i=1,2\big\}$. If
$x\in A$ then $u$ cannot be differentiable at $x$ thanks to \eqref{vis-cond}. Notice that 
if $\Sigma$ is a strictly concave graph over the plane $\{x_{n+1}=0\}$ then 
\eqref{vis-cond} is satisfied, see Remark \ref{rem-extns}.  
By Aleksandrov's theorem the concave function $u$ is twice
differentiable a.e. Hence $|A|=0$.

As for  c) we must check that the following three conditions hold, see e.g. \cite{Ash}
\begin{itemize}
 \item[1)] $\V\in \mathscr F$,
\item[2)] if $A\in\mathscr F$ then $\V\setminus A\in \mathscr F$,
\item[3)] if $A_i\in\mathscr F$ then $\bigcup_{i=1}^\infty A_i\in \mathscr F$.
\end{itemize}
We first prove 1). If $A_i\in \V$ is any sequence of subsets of $\V$ then clearly
$\s(\cup_{i=1}^\infty A_i)=\cup_{i=1}^\infty \s(A_i)$.
Writing $\V=\cup_{i=1}^\infty E_i$, where $E_i\subset \V$ are closed subsets we conclude that
$\s(\V)=\s(\cup_{i=1}^\infty E_i)=\cup_{i=1}^\infty\s(E_i)$. From a) it follows
that $\s(E_i)$ is closed for any $i$, and hence measurable, implying that $\s(\V)$ is measurable.

\smallskip

2) Let $A\in \mathscr F$. We use the following elementary identity
\begin{equation}\label{compl-ident}
 \s(\V\setminus A)=[\s(\V)\setminus \s(A)]\bigcup [\s(\V\setminus A)\cap \s(A)].
\end{equation}

From b) it follows that $|\s(\V\setminus A)\cap \s(A)|=0$. Therefore
$|\s(\V\setminus A)|=|\s(\V)\setminus \s(A)|$ and 2) is proven.

\smallskip
It remains to check 3). Without loss of generality we assume that $A_i$'s are disjoint, see \cite{Ash}.
Thus, letting $A_i\in\mathscr F, A_i\cap A_j=\emptyset, i\not=j$ we get
\begin{eqnarray*}
 \sum_{i=1}^\infty|\s(A_i)|&\geq& |\s(\cup_{i=1}^\infty A_i)|\ge\\\nonumber
&\geq &\sum_{i=1}^\infty|\s(A_i)|-
\sum_{ij=1}^\infty|\s(A_i)\cap\s(A_j)|\geq\\\nonumber
&\geq &\sum_{i=1}^\infty|\s(A_i)|.
\end{eqnarray*}
 \qed

\smallskip

For a given function $u\in \Wu(\U,\V)$ we consider  the set function
\begin{eqnarray}\label{def-beta}
\beta_{u} (\omega)= \int_{\s(\omega)}f
\end{eqnarray}
where $\omega\subset \V$ is a Borel subset.
Since $\mathscr F$ contains the closed sets (see part a) above) we infer that $\beta_{u,f}$ is a Borel measure.
Moreover, from the proof of Proposition \ref{prop-S} b) it follows  that $\beta_{u,f}$ is countably additive.

\begin{definition}\label{def-beta-weak}
A function $u$  (or its graph $\Ga_u$) is said to be a $B$-type weak solution to
(\ref{problem}) if $u\in \Wup(\U,\V)$ and
the following two identities holds
\begin{eqnarray}\label{eq-beta}
  \left\{
  \begin{array}{ccc}
  \displaystyle\beta_{u,f}(\omega)=\int_{\omega}gd\H, \ \hbox{for\ any\ Borel\ set}\ \omega\subset \V \qquad \hbox{and}\\
   \s_u(\V)=\U.
   \end{array}
   \right.
\end{eqnarray}
\end{definition}

\subsection{Existence of weak solutions of $B$-type}\label{sec-B-ex}
The measure $\beta$, defined in (\ref{def-beta}) is weakly continuous. We have

\begin{lemma}\label{lem-weakconv-B}
  Let $u_k$ be a sequence of $B$-type weak solutions in the sense of Definition \ref{def-beta} and $\beta_k$
  is the associated measure, defined by (\ref{def-beta}). If $u_k\rightarrow u$ uniformly on
  compact subsets of ${\U}$ then $u$ is $R-$concave and $\beta_k$ weakly
  converges to $\beta_{u,f}$.
\end{lemma}
\pr Once the $\sigma$ additivity is established then the proof of lemma is standard,  
see for the classical case \cite{P-88} pp 14-18,\cite{Gut-b} and \cite{GutH}, \cite{K-ref}  for refractor and reflector problems respectively. 
The proof for supporting ellipsoids is carried out also \cite{GutT}. 
That $u$ is admissible follows from Lemma \ref{lem-support}. Recall that the weak
convergence is equivalent to the following two inequalities (see \cite{Ash} Theorem 4.5.1)

\begin{itemize}
  \item[1)] $\limsup\limits_{k\rightarrow \infty} \beta_k(E)\leq \beta(E)$ for any compact $E\subset \V,$
  \item[2)] $\liminf\limits_{k\rightarrow \infty} \beta_k(J)\geq \beta(J)$ for any open $J\subset \V$.
\end{itemize}
Take a closed set $E$ and let $E^*_\delta$ be an $\delta-$neighbourhood of the closed set $E^*=\s(E)$, see Lemma \ref{prop-S} a).
We claim that for any $\delta>0$ there is $i_0\in \mathbb N$ such that
$\s_i(E)\subset E^*_\delta$ whenever $i>i_0$, where
$\s_i$ is the mapping corresponding to $u_i$. If this fails then
there is $\delta>0$ and a sequence of points $x_i\in\s_i(E)$ such that $x_i\in \complement E_\delta^* $.
By definition there is $Z_i\in E$ such that $x_i\in \s_i(Z_i)$. Suppose that
$x_i\rightarrow x_0$, for some $x_0$,   and $Z_i\rightarrow Z_0\in E$ at least for a subsequence.
Thus, $x_0\in \complement E^*_\delta, x_0\in \s(Z_0)$ and $Z_0\in E$ which is a contradiction.

\smallskip

To prove the second inequality we let $J\subset \V$ be an open subset and denote $J^*=\s(H)$.
By Lemma \ref{prop-S} c) $J^*$ is measurable, hence  for any 
small $\delta>0$ there is a closed set $J^*_\delta$ such that
$J^*_\delta\subset J^*$ and $|J^*|-\delta\leq |J^*_\delta|\leq |J^*|$. 
This is possible because by Proposition \ref{prop-S} b)
$\s$ is one-to-one modulo a set of measure zero.
Let $N_\delta$ be an open set, $|N_\delta|<\delta$ containing the points
where the inverse of $\s$ is not defined. We claim that there is $k_0$ such that 
\begin{equation}\label{eq-H-open}
J^*_\delta\setminus N_\delta\subset J^*_k\eqdef
\s_k(J), \quad {\rm for\ any}\ k\geq k_0.
\end{equation}
Here $\s_k$ is the mapping generated by $u_k$.
Proof of (\ref{eq-H-open}) is by  contradiction. If (\ref{eq-H-open}) fails then
there is $x_k\in J^*_\delta\setminus N_\delta$ and $x_k\not \in J^*_k$.
We can assume that
$x_k\rightarrow x_0$. Since $J^*_\delta\setminus N_\delta$ is closed it follows that
$x_0\in J^*_\delta\setminus N_\delta$. By definition of $N_\delta$ the inverse of $\s$ is
one-to-one on $J^*_\delta\setminus N_\delta$. Thus there is a unique $Z_0\in H$ such that
$x_0=\s(Z_0)$. Furthermore, there is an open neighborhood of $Z_0$ contained in $J$ because $J$ is open.
If $H(x,\sigma_k,Z_k)$ is a supporting hyperboloid of $u_k$ at $x_k$ it follows from
Lemma \ref{lem-support} that $x_k\in \s_k(Z_k), Z_k\rightarrow Z_0$. Thus for large $k$, $\{Z_k\}$ is in some neighborhood of
$Z_0\in J$  implying that $x_k\in J^*_k$ which contradicts our supposition.\qed

\begin{proposition}\label{prop-B-exst}
 Let $f:\U\to \R$ and $g: \V\to \R$ be two nonnegative integrable functions. If $\U\subset \Pi$ and
$\V\subset \Sigma$ are bounded domains such that the energy balance condition (\ref{glob-balance}) and (\ref{vis-cond-B}) hold then there exists  a $B-$type weak solution  to the problem (\ref{problem}).
\end{proposition}
Notice that we do not exclude the case $\U\cap \widehat \V\not=\emptyset.$

\pr The proof of Proposition \ref{prop-B-exst} is by approximation 
argument, see \cite{GutT}, \cite{K-ref}, \cite{W-96}. Let  
$g_N=\sum_{i=1}^NC_i\delta_{Z_i}$ with $C_i\geq 0$ such that
$\sum_{i=1}^NC_i=\int_\U f(x)dx, Z_i\in \Sigma$ and $\delta_{Z_i}$ are atomic measures 
supported at $Z_i$. For each $g_N$ we construct a $B-$type solution $u_N$. Then sending $N\to \infty$ and using
the compactness argument together with 
weak convergence of $g_N$ to
$g$, Lemma \ref{lem-weakconv-B}, one will arrive at desired result.

%%%%%%%%%%%%%%%%%%%%%%%%%%%%%%%%%%%%
% case of N points

First, for each $Z\in \Sigma$ we  define 
\begin{equation}\label{max-a}
\bar a(Z)=\frac{\eps Z^{n+1}-\sqrt{(Z^{n+1})^2+\rho^2}}{\eps^2-1}
\end{equation}
where 
\begin{equation}\label{max-rad}
\rho(z)=\inf\{R>0 : \U\subset B_{R}(z)\}.
\end{equation}

Clearly $\bar a(Z)$ is the maximal value of larger semiaxis of hyperboloid  $H$
such that $\Gamma_H$ is visible from $\U$ in the $e_{n+1}$ direction. In other words $H(x, \bar a(Z), Z )$
 is the lowest possible hyperboloid with focus $Z\in \Sigma$ such that $H(x, \bar a(Z), Z)\geq 0$.
 Thus for $a\in (0, \bar a(Z)]$ we have $H(\cdot, a, Z)\in \mathbb H^+_{0}(\U, \V)$. 
To check (\ref{max-a}) we fix $Z$ and pick $x_0$ such that 
$\rho(z)=|x_0-z|$. Since the ratio of distances of $x_0$ from lower focus 
$Z'$ and the plane $\Pi_d=\{X\in \R^{n+1} : X^{n+1}=Z^{n+1}-\bar a\eps-\bar a/\eps\}$ is $\eps$, it follows that 
$|x_0Z'|=\eps(Z^{n+1}-\bar a\eps-\bar a/\eps)$. On the other hand 
$|x_0Z|-|x_0Z'|=2\bar a$. Consequently, we find that 
$\sqrt{(Z^{n+1})^2+\rho^2(z)}=2\bar a+\eps(Z^{n+1}-\bar a\eps-\bar a/\eps)$ which gives 
(\ref{max-a}).

Next we define the maximal level $L_0=\sup\limits_{\U\times \V}H(x, \bar a(Z), Z)$. Since 
\begin{eqnarray*}
\max_{x\in \U} H(x, \bar a(Z), Z)&=&Z^{n+1}-\eps\bar a(Z)-\bar a(Z)=\frac{\sqrt{(Z^{n+1})^2+\rho^2(z)}-Z^{n+1}}{\eps-1}\\\nonumber
&=&\frac{\rho^2(z)}{(\eps-1)(\sqrt{(Z^{n+1})^2+\rho^2(z)}+Z^{n+1})}
\end{eqnarray*}
it follows that 
\begin{equation}\label{def-L-0}
L_0=\sup_{\V}\frac{\rho^2(z)}{(\eps-1)(\sqrt{(Z^{n+1})^2+\rho^2(z)}+Z^{n+1})}\leq \frac1{\eps-1}\sup_\V\frac{\rho^2(z)}{2Z^{n+1}}.
\end{equation}

Next, we  bound $H(\cdot, a, Z)$ by below for $a>0$ close to zero.
By definition (\ref{def-hypr}) we have that for this case $H(x, a, Z)\sim Z^{n+1}-\frac{\rho(z)}{\sqrt{\eps^2-1}}$.
We demand $Z^{n+1}-\frac{\rho(z)}{\sqrt{\eps^2-1}}\ge 2L_0$ or equivalently in lieu of 
(\ref{def-L-0})
$$Z^{n+1}\geq \rho(z)\left[\frac{1}{\sqrt{\eps^2-1}}+\frac{2\rho(z)}{(\eps-1)(\sqrt{(Z^{n+1})^2+\rho^2(z)}+Z^{n+1})}\right].$$
But clearly $\frac{2\rho(z)}{(\eps-1)(\sqrt{(Z^{n+1})^2+\rho^2(z)}+Z^{n+1})}\leq 2/{(\eps-1)}$. Therefore
 it is enough to 
assume that 
$Z^{n+1}\geq [\frac2{\eps-1}+\frac1{\sqrt{\eps^2-1}}]\rho(z)$ which is exactly (\ref{vis-cond-B}).
It follows that if $\Sigma$ satisfies (\ref{vis-cond-B}) then $\widetilde \Sigma=\Sigma+Me_{n+1}, M\gg1$ 
also does.

\smallskip 
%%%%
Let ${\bf a}=(a_1, \dots, a_N), a_i\in(0, \bar a(Z_i)], i=1, \dots, N$ and set 
$$H({\bf a}, x)=\min\big[H(x,a_1,Z_1),\dots,H(x,a_N, Z_N)\big].$$
We also let $\cal E_i({\bf a})=\{x\in \U : H({\bf a}, x)=H(x, a_i, Z_i)\}$ be the $i$-th visibility sets  and 
$$\cal A^N=\left\{{\bf a}\in\prod\limits_{i=1}^N(0, \bar a_i(Z_i)] : 
\int_{\cal E_i({\bf a})}f\leq C_i, \int_{\cal E_N({\bf a})}f\geq C_N,\quad  i=1,\dots,
N-1
\right\}.$$

From (\ref{vis-cond-B}) it follows that $\cal A^N$
is not empty for taking $a_i, 1\leq i \le N-1$ close to zero and $a_N=\bar a_N(Z_N)$ one readily gets that 
such ${\bf a}$ is in $\cal A^N.$

The visibility sets $\cal E_i({\bf a})$ enjoy the following property:
if for some $a_k<\bar a(Z_k)$ we set ${\bf a}^k_\delta=(a_1, \dots, a_{k}+\delta, \dots, a_N)$
and ${\bf a}=(a_1,\dots, a_N)$ for $\delta >0$ small, then 

\begin{equation}\label{mon-incl}
\cal E_{k}({\bf a})\subset \cal E_{k}({\bf a}_\delta^k)
\quad {\rm while}\quad  \cal E_{i}({\bf a}_\delta^k)\subset \cal E_{i}({\bf a}), i\not =k.
\end{equation} 
This can be seen for $N=2$
by simple geometric considerations, and general case is by induction.

Let $\mathfrak a=\sup\limits_{{\bf a}\in \cal A^N}\sum\limits_{i=1}^N a_i$ and $\hat{\bf a}\in \cal A^N$ be such that 
the supremum is realised, i.e. $\mathfrak a=\sum\limits_{i=1}^N \hat a_i$. We claim that
$H(\hat {\bf a},x)$ solves the refractor problem with measure $g_N$. If not, then 
there is $i_0$, say $i_0=1$, such that $\int_{\cal E_{1}(\hat{\bf a})}f<C_1$. Then in view of the 
energy balance condition this implies $\int_{\cal E_{N}(\hat{\bf a})}f>C_N$.
For $\delta>0$ small $F_N(\delta)=\int_{\cal E_N (\hat {\bf a}^1_\delta) }f(x)dx\ge C_N$ because 
$F_k(\delta)$ is continuous function of $\delta$. Furthermore, using (\ref{mon-incl}) it follows that 
${\bf a}^1_\delta\in \cal A^N$ which is a contradiction. 
Now the proof of Theorem C1 follows from the above polyhedral approximation $H({\bf \hat a}, x)$ as $N\to \infty$ 
and the weak convergence of measures $\beta_{H,f}$, Lemma \ref{lem-weakconv-B}.
\qed

%%%%%%%%%%%%%%%%%%%%%%%%%%%%%%%%%%%%%%%%%%%%%%%%%%%%%%%%%%%%
%                                                          %
%                SECTION                                   %
%                                                          %
%%%%%%%%%%%%%%%%%%%%%%%%%%%%%%%%%%%%%%%%%%%%%%%%%%%%%%%%%%%%
\section{An approximation lemma}

\subsection{Refraction cone}
Recall that for smooth refractors the unit direction of the refracted ray is
$$Y=\eps\left(e_{n+1}+\gamma\left[\sqrt{(\gamma\cdot e_{n+1})^2-\kap}-\gamma\cdot e_{n+1}\right]\right),$$
see (\ref{eq-Y}).
This formula may be generalized for non smooth refractors as follows: let $\gamma_1, \gamma_2$ be the
normals of two supporting planes of $u$ at $x$. Then for any two constants $c_1, c_2$ the unit vector
$\gamma_{c_1c_2}=\frac{c_1\gamma_1+c_2\gamma_2}{|c_1\gamma_1+c_2\gamma_2|}$ generates a mapping to the 
unit sphere $\mathbb S^{n+1}$ given by
$$\gamma_{c_1c_2}\mapsto \eps\left(e_{n+1}+\gamma_{c_1c_2}\left[\sqrt{(\gamma_{c_1c_2}\cdot e_{n+1})^2-\kap}-\gamma_{c_1c_2}\cdot e_{n+1}\right]\right).$$

\smallskip

\begin{definition}
  For $\gamma_1, \gamma_2 \in \mathbb S^{n+1}$ the refractor cone at $\xi\in \R^{n+1}$ is defined
  as $$C_{\xi, \gamma_1, \gamma_2}=\left\{Z\in \R^{n+1} : \frac{Z-\xi}{|Z-\xi|}=\eps\left(e_{n+1}+\gamma_{c_1c_2}\left[\sqrt{(\gamma_{c_1c_2}\cdot e_{n+1})^2-\kap}-\gamma_{c_1c_2}\cdot e_{n+1}\right]\right)\right\}.$$
\end{definition}

\smallskip

One can easily verify that $C_{\xi, \gamma_1, \gamma_2}$ is a convex cone. Indeed, for any $\gamma_0\perp \hbox{Span}\{ \gamma_1, \gamma_2 \}$ we have that $\frac{Z-\xi}{|Z-\xi|}\cdot \gamma_0=\eps ( e_{n+1}\cdot\gamma_0)$.
Thus $C_{\xi, \gamma_1, \gamma_2}$ is a cone.

In view of Lemma \ref{lem-lips} $\|Du\|_{L^\infty(\U)}\le \frac1{\sqrt{\eps^2-1}}$ for any admissible $u\in \Wup (\U, \V)$, and $\sqrt{(\gamma\cdot e_{n+1})^2-\kap}$ is well defined thanks to this gradient estimate.

%%%%
\subsection{Contact set}\label{subsec-contact}
In this subsection we study the contact set of two hyperboloids
$$\Lambda=\{x\in \R^n\ :\ H(x, a_1, Z_1)=H(x, a_2, Z_2)\}$$
where $Z_i\in \Sigma, a_i>0, i=1.2$. We show that 
$\Lambda$ is a conic section.
To check this we simplify the equation 
\begin{equation*}
Z_1^{n+1}-Z_2^{n+1}-\eps(a_1-a_2)-\sqrt{a_1^2+\frac{|x-z_1|^2}{\eps^2-1}}=-\sqrt{a_2^2+\frac{|x-z_2|^2}{\eps^2-1}}
\end{equation*} 
by squaring both sides of it.
Then denoting $C=Z_1^{n+1}-Z_2^{n+1}-\eps(a_1-a_2)$
we infer
\begin{eqnarray*}
C^2-2C\sqrt{a_1^2+\frac{|x-z_1|^2}{\eps^2-1}}+a_1^2+\frac{|x-z_1|^2}{\eps^2-1}=a_2^2+\frac{|x-z_2|^2}{\eps^2-1}.
\end{eqnarray*}
Recognizing the terms and denoting $D=\frac12\left[\frac{|z_1|^2-|z_2|^2}{\eps^2-1}+C^2+a_1^2-a_2^2\right]$ we get 
\begin{equation*}
D+\frac{x\cdot(z_2-z_1)}{\eps^2-1}=C\sqrt{a_1^2+\frac{|x-z_1|^2}{\eps^2-1}}.
\end{equation*}
By choosing a suitable coordinate system we can assume that $z_2-z_1$ is collinear to the unit direction of $x_1$ axis.
Thus 
\begin{equation*}
D+\frac{x_1|z_2-z_1|}{\eps^2-1}=C\sqrt{a_1^2+\frac{|x-z_1|^2}{\eps^2-1}}.
\end{equation*}
Finally squaring  both sides of the last identity and assuming that $E:=\frac{|z_1-z_2|^2}{\eps^2-1}-C^2\not=0$
 we infer 
\begin{eqnarray}\nonumber
C^2\frac{|x'-z'_1|^2}{\eps^2-1}&=&D^2-C^2 a_1^2+\frac1{\eps^2-1}\left(2Dx_1|z_1-z_2|+\frac{x_1^2|z_1-z_2|^2}{\eps^2-1}
-C^2\left[x_1^2-2x_1z_1^1+(z_1^1)^2\right]\right)\\\nonumber
&=& D^2-C^2 a_1^2+\frac{E}{\eps^2-1}\left(x_1^2+2x_1\frac{D|z_1-z_2|+C^2z_1^1}{E}-\frac{C^2(z_1^1)^2}{E}\right)\\\nonumber
&=& D^2-C^2 a_1^2+\frac{E}{\eps^2-1}\left(\left[x_1+\frac{D|z_1-z_2|+C^2z_1^1}{E}\right]^2 
-\left[\frac{D|z_1-z_2|+C^2z_1^1}{E}\right]^2-\frac{C^2(z_1^1)^2}{E}\right)\\\nonumber
&=& F+\frac{E}{\eps^2-1}\left[x_1+\frac{D|z_1-z_2|+C^2z_1^1}{E}\right]^2
\end{eqnarray}
where 
\begin{equation}
E= \frac{|z_1-z_2|^2}{\eps^2-1}-C^2, \quad F= D^2-C^2 a_1^2-\frac{E}{\eps^2-1}\left( 
\left[\frac{D|z_1-z_2|+C^2z_1^1}{E}\right]^2+\frac{C^2(z_1^1)^2}{E}\right).
\end{equation}
Note that if $E=0$ then $\Lambda$ is a paraboloid. Otherwise 
\[
\Lambda\ \text{is}\left\{
\begin{array}{lll}
\text{the sheet of  hyperboloid}\  \frac{|x'-z_1'|^2}{a^2}=F+\frac{|x_1-{x_1^0}|^2 }{b^2}\ \text{passing through}\ x^0,& \text{if}\  E>0,\\
\text{the ellipsoid}\  \frac{|x'-z_1'|^2}{A^2}+\frac{|x_1-{x_1^0}|^2 }{B^2}=F\quad \text{with}\quad  \frac{A^2}{B^2}=\frac{C^2-|z_1-z_2|^2}{C^2},
& \text{if}\  E<0.
\end{array}
\right.
\] 
We see that the rotational axis of $\Lambda$ for both cases $E<0$ and $E>0$ is parallel to the direction of $z_2-z_1$.
Moreover, if $\Lambda$ is an ellipsoid then this direction corresponds to the larger semiaxis.
This observation will be used in the proof of Lemma \ref{loc-glob} below.
%%%%%%
\subsection{R-convexity of $\V$}
\begin{definition}\label{def-R-conv}
  We say that $\V\subset \Sigma$ is $R-$convex with respect to a point $\xi\in [0,\infty)\times \U$
  if for any two unit vectors $\gamma_1, \gamma_2$ the intersection $C_{\xi, \gamma_1, \gamma_2}\cap \V$
  is connected. If $\V$ is $R-$convex with respect to any $\xi\in [0, \infty)\times \U$ then we simply say that
  $\V$ is $R-$convex.
\end{definition}
In particular a  geodesic ball on the  convex  surface $\Sigma$ is an example of $R-$convex $\V$. 

\subsection{Local supporting function is also global}
In the Definition \ref{def-adms} of admissibility the supporting hyperboloid $H$ is staying above 
$u$ in whole $\U$. Consequently, one may wonder if the locally admissible functions (i.e. $H$ stays above $u$ only in a vicinity of the contact point) are still in $\Wup(\U, \V)$.   This issue was addressed by G. Loeper in \cite{Loep} for the optimal transfer problems. We have  
\begin{lemma}\label{loc-glob}
Under the condition (\ref{A3-cond00}) a local supporting hyperboloid is also global.
\end{lemma}
\pr 
The proof is very similar to that of in \cite{Loep}, \cite{Turd-Wang-scrt}.
Let $H_i(x)=H(x, a_i, Z_i), i=1, 2$ be two global supporting hyperboloids of $u$ at $x_0$ such that 
the contact set $\Lambda\not=\{x_0\}$. Thus $u$ is not differentiable at $x_0$. To fix the ideas take $x_0=0$.
%Suppose that there is a local support hyperboloid $H(x)=H(x, a, Z)$ at $x_0$.
If $\gamma_i$ is the normal of the graph of $H_i, i=1,2$ at $x_0$
then for any   $\theta\in(0, 1)$ there is $Z_{\theta}\in \Sigma\cap \mathcal C_{0, \gamma_1,\gamma_2}$ and $a_{\theta}>0$
such that  $H(x)=H(x, a_{\theta}, Z_{\theta})$ is a local supporting hyperboloid 
of $u$ at $0$ and 
\begin{equation}\label{blya-2}
DH_{\theta}(0)=\theta DH_1(0)+(1-\theta)DH_2(0).
\end{equation} 
Observe that the correspondence $\theta\mapsto Z_\theta$ is one-to-one thanks to the 
assumption \eqref{vis-cond}, see Remark \ref{rem-extns}. 
By choosing a suitable coordinate system we can assume that $DH_1(0)-DH_2(0)=(0, \dots, 0, \alpha)$ for some 
$\alpha>0$. 
%Moreover, if $\mathscr  L$ is a linear function then 
%$\theta(H_1-\mathscr L)+(1-\theta)(H_2-\mathscr L)-(H-\mathscr L)=\theta H_1+(1-\theta)H_2-H$. Thus without loss of 
%generality we assume that $DH_2(0)=0$. 
Then  we have that for all $0<\theta<1$
\begin{eqnarray}\label{blya-3}
\min[H_1(x), H_2(x)]&\le &\theta H_1(x)+(1-\theta)H_2(x)\\\nonumber
&=&u(0) + \left[DH_2(0)+\alpha\theta\right]x_n+\frac12\left[\theta D^2H_1(0)+(1-\theta)D^2 H_2(0)\right]{x\otimes x}+o(|x|^2)
\end{eqnarray}
where the last line follows from Taylor's expansion.
 
Using the notations of Section \ref{sec-weak-sol} we have that 
\begin{equation*}
D^2H_\theta(0)=-\frac{G(x_0, u(0), p_1+\theta(p_2-p_1))}{\eps\kappa}
\end{equation*}
where we set $p_i=DH_i(0), i=1,2$ and used \eqref{blya-2}. For all unit vectors $\tau$ perpendicular to $x_n$ axis we have 
\begin{eqnarray}
\frac{d^2}{d\theta^2}D^2_{\tau\tau}H_\theta(0)&=&
-\frac{d^2}{d\theta^2}\frac{G^{ij}(0, u(0), p_1+\theta(p_2-p_1))\tau_i\tau_j}{\eps\kappa}\\\nonumber
&=&-\alpha^2\frac{\partial^2}{\partial p_n^2} \frac{G^{jj}(0, u(0), p_1+\theta(p_2-p_1))\tau_i\tau_j}{\eps\kappa}\\\nonumber
&\le&-\alpha^2 c_0
\end{eqnarray}
where the last line follows from \eqref{A3-cond00} with $c_0>0$, see also \eqref{eq-A3} and  Section \ref{ssec-example}.

Therefore
\begin{equation}\label{blya-4}
D^2_{\tau\tau}H_\theta(0)\ge\theta D^2_{\tau\tau}H_1(0)+(1-\theta)D^2_{\tau\tau}H_2(0)+\widehat{ c_0}\theta(1-\theta)|p_1-p_2|^2
\end{equation}
where $\widehat{c_0}$ depends on $c_0$.

Observe that at $0$ we have 
\begin{equation}
\theta p_1+(1-\theta)p_2=\frac1{\eps^2-1}\frac{z_\theta}{\sqrt{a_\theta^2+\frac{|z_\theta|^2}{\eps^2-1}}}, 
\quad \phi(z_\theta)-u(0)-\eps a_\theta=\sqrt{a_\theta^2+\frac{|z_\theta|^2}{\eps^2-1}}
\end{equation}
where we assume that $\Sigma$ is the graph of a function $\phi$  such that $\psi(Z)=Z^{n+1}-\phi(z)$ satisfies \eqref{A3-cond00}.
From these $n+1$ equations we see that $z_\theta$ and $a_\theta$ are smooth functions of $\theta p_1+(1-\theta)p_2$.
This yields  the following crude estimate for the remaining second order derivatives
\begin{eqnarray}\label{blya-5}
\left|\theta D_{jn}^2H_1(0)+(1-\theta)D_{jn}^2H_2(0)-D^2_{jn}H(0)\right|\le C\theta(1-\theta)|p_1-p_2|^2,\quad  j=1, \dots, n%|x_j-x^0_j||x_k-x^0_k|
\end{eqnarray}
where $C$ depends of $C^2$ form of $\phi$.
Consequently, after plugging \eqref{blya-4} and \eqref{blya-5} into \eqref{blya-3} and recalling that 
$|p_1-p_2|=\alpha$ we conclude

\begin{eqnarray}\nonumber
\min[H_1(x), H_2(x)]&\le & u(0)+\left[p_2+\alpha\theta\right] x_n+\frac12\left[\theta D^2H_1(0)+(1-\theta)D^2 H_2(0)\right]{x\otimes x}+o(|x|^2)\\\nonumber
&=&u(0)+\left[p_2+\alpha\theta\right] x_n+
\frac12 D^2H_\theta (0){x\otimes x}-\widehat{ c_0}\theta(1-\theta)\alpha^2\sum_{j=1}^{n-1}x_j^2+ C\theta(1-\theta)\alpha^2|x_n||x|+ o(|x|^2)\\\nonumber
&=&u(0)+\left[p_2+\alpha\theta\right] x_n+
\frac12 D^2H_\theta (0){x\otimes x}-\\\nonumber
&&-\widehat{ c_0}\theta(1-\theta)\alpha^2|x|^2+ \theta(1-\theta)\alpha^2(C|x_n||x|+\widehat{ c_0}x_n^2)+ o(|x|^2)\\\nonumber 
&=&u(0)+\left[p_2+\alpha\theta\right] x_n+
\frac12 D^2H_\theta (0){x\otimes x}-\\\nonumber
&&-\frac{\widehat{ c_0}}2\theta(1-\theta)\alpha^2|x|^2+ \theta(1-\theta)\alpha^2\left(\frac{2C^2}{\widehat{c_0}}+\widehat{ c_0}\right)x_n^2+ o(|x|^2)
\end{eqnarray}
where the last line follows from H\"older's inequality.
Now fixing $\theta_0$ as in  \eqref{blya-2} and using the estimate $|D^2H_\theta (0)-D^2H_{\theta_0} (0)|\le C|\theta-\theta_0|$
(with $C$ depending on $\phi$)
we obtain 
\begin{eqnarray}\nonumber
\min[H_1(x), H_2(x)]&\le &u(0)+\left[p_2+\alpha\theta_0\right] x_n+
\frac12 D^2H_{\theta_0} (0){x\otimes x}+\\\nonumber
&&+\alpha\left[\theta-\theta_0\right] x_n+\frac12 D^2H_{\theta} (0){x\otimes x}-\frac12 D^2H_{\theta_0} (0){x\otimes x}\\\nonumber
&&-\frac{\widehat{ c_0}}2\theta(1-\theta)\alpha^2|x|^2+ \theta(1-\theta)\alpha^2\left(\frac{2C^2}{\widehat{c_0}}+\widehat{ c_0}\right)x_n^2+ o(|x|^2)\\\nonumber
&\le &H_{\theta_0} (x)+\alpha\left[\theta-\theta_0\right] x_n+C|\theta_0-\theta||x|^2-\\\nonumber
&&-\frac{\widehat{ c_0}}2\theta(1-\theta)\alpha^2|x|^2+ \alpha^2\left(\frac{2C^2}{\widehat{c_0}}+\widehat{ c_0}\right)x_n^2+ o(|x|^2).\\\nonumber
\end{eqnarray}
Choosing $\theta=\theta_0-x_n\alpha\left(\frac{2C^2}{\widehat{c_0}}+\widehat{ c_0}\right)$ such that 
$|x_n|<\delta$ with sufficiently small $\delta$
we finally obtain 
\begin{eqnarray}\label{Leop-000}
\min[H_1(x), H_2(x)]
&\le &H_{\theta_0} (x)+C|\theta_0-\theta||x|^2-\frac{\widehat{c_0}}2\theta(1-\theta)\alpha^2|x|^2+ o(|x|^2)=\\\nonumber
&=& H_{\theta_0} (x)-\frac{\widehat{c_0}}2\theta(1-\theta)\alpha^2|x|^2+ o(|x|^2)\\\nonumber
&\le& H_{\theta_0} (x), \quad \forall x\in B_\delta.
\end{eqnarray}

This, in particular, implies that $H_{\theta_0}$ is a local supporting hyperboloid near $x=0$. 

It remains to check that $H_{\theta_0}$ is also a global supporting hyperboloid. 
The set $\Lambda_{1,2}=\{x\in \R^n\ :\  H_1(x)=H_2(x)\}$ passes through $0$ and splits $\U$ into two parts $\U^+$ and $\U^-$
(recall that $\Lambda_{12}$ is a conic section, see Section \ref{subsec-contact}).
It follows from \eqref{Leop-000} that the contact sets $\Lambda_{i,\theta_0}=\{x\in \R^n\ :\  H_i(x)=H_{\theta_0}(x)\}, i=1,2$
are tangent to $\Lambda_{12}$ from one side in some vicinity of $0$, say in $\U^+$. 
If there is $\bar x_0\not =0$ such that, say,  
$\bar x_0\in \Lambda_{1,2}\cap \Lambda_{1, \theta_0}$ then $H_1(\bar x_0)=H_2(\bar x_0)=H_{\theta_0}(\bar x_0)$
and $DH_{\theta_0}(\bar x_0)=\bar \theta_0 DH_1 (\bar x_0)+(1-\bar\theta_0)DH_2(\bar x_0)$ with possibly 
different $\bar\theta_0$. 
Observe that by construction the ray emitted from $\bar x_0$ in the direction of $e_{n+1}$ after refraction 
from $H_1, H_2$ and $H_{\theta_0}$ hits the point $Z_1, Z_2$ and $Z_{\theta_0}$, respectively.
Then repeating the argument above with $0$ replaced by $\bar x_0$ and $\theta_0$ by $\bar\theta_0$ (but keeping $H_{\theta_0}$ fixed),  we can see that 
\eqref{Leop-000} is satisfied in 
$B_\delta(\bar x_0)$ implying that $ \Lambda_{1, \theta_0}$ is tangent with $ \Lambda_{1, 2}$ at $\bar x_0$ and lies in $\U^+$.
Thus $H_{\theta_0}$ is a global supporting hyperboloid.
\qed

\smallskip 

As an application of Lemma \ref{loc-glob} we have the following approximation result.

\begin{lemma}\label{lem-apprx}
 If $u\in \Wup(B_r, \Sigma)$ then
 \begin{itemize}
  \item[(i)] $u_\eps(x) +K(r^2-|x|^2)\in \Wup(B_r, \widetilde \Sigma)$
  where $u_\eps$ is the standard mollification of $u$, $K>0$ and $\widetilde \Sigma=\Sigma +Me_{n+1}$
  for some large constant $M>0$,
  \item[(ii)] $u_\eps(x) +K(r^2-|x|^2)$ is a classical subsolution of (\ref{M-A-short+}).
 \end{itemize}
\end{lemma}
\pr
(i) It is well known that $u_\eps$ is concave and $\|Du_\eps\|_{L^\infty(B_r)} \leq \|Du\|_{L^\infty(B_r)}<\frac1{\sqrt{\eps^2-1}}$. Therefore if $K$ is fixed then we can choose $r$ so small that 
\begin{equation}\label{sml-grad}
\|D\bar u_\eps\|_{L^\infty(B_r)} \leq \|Du\|_{L^\infty(B_r)}+2Kr<\frac1{\sqrt{\eps^2-1}}.
\end{equation}
Moreover
$K(r^2-|x|^2)$ is concave, hence $\bar u_\eps=u_\eps(x) +K(r^2-|x|^2)$ is concave as well. Notice that 
$D^2\bar u_\eps =D^2u_\eps-2K\Id\leq -2K\Id<0$ implying that $\bar u_\eps$ is strictly concave.
In order to bound the curvature of $\Ga_{\bar u_\eps}$ from below we recall that 
for fixed $Z$, $H(\cdot, a, Z)$ becomes flatter as $a\to \infty$ because 
$$D^2H=-\frac{1}{(\eps^2-1)\sqrt{a^2+\frac{|x-z|^2}{\eps^2-1}}}\left[\Id-\frac{(x-z)\otimes(x-z)}{(\eps^2-1)a^2+|x-z|^2}\right].$$  
In particular, for large $K$ and $a$ we will have $-D^2\bar u_\eps\geq 2K\Id \ge -D^2H$. 
Consequently, for each $x\in \U$ there is $Z$ and $a>0$ such that 
$H(\cdot, a, Z)$ touches $\bar u_\eps$ from above at $x$, in some neighbourhood of $x$. Furthermore,
from Lemma \ref{lem-conf} on confocal expansion we can choose $a, \widetilde Z$ so that 
$\widetilde Z\in \widetilde \Sigma=\Sigma+Me_{n+1}, M\gg1.$
Finally applying Lemma \ref{loc-glob} we infer that $H(\cdot, a, \widetilde Z)$ is a global supporting hyperboloid 
of $u$ at $x$ and thus $\bar u_\eps\in \Wup(\U, \widetilde\Sigma)$.
%%%%%%%%%%%%%%%%%%%%%%%%%%%%%%%%%%%%%%%%%%%%%%%%%%%%%%%%%%%

(ii) By direct computation we have 
\begin{eqnarray}
\cal M&=&-D^2\bar u_\eps-\frac{G(x, \bar u_\eps, D\bar u_\eps)}{\eps\kap}=-D^2 u_\eps +2K\Id-\frac{G(x, \bar u_\eps, D\bar u_\eps)}{\eps\kap}\ge \\\nonumber
&\ge&2K\Id -\frac{G(x, \bar u_\eps, D\bar u_\eps)}{\eps\kap}.
\end{eqnarray}
By definition, (\ref{def-G-M}) we have $$\frac{G(x, \bar u_\eps, D\bar u_\eps)}{\eps\kap}=\frac{[q+1](\Id-\eps^2\kap D\bar u_\eps\otimes D\bar u_\eps)}{\eps\kap t(x, \bar u_\eps, D\bar u_\eps)}\leq \frac C{t(x, \bar u_\eps, D\bar u_\eps)}$$
with some tame constant $C>0$ depending only on $\eps$. Recall that by (\ref{sml-cone}) $t=\frac{(M+[Z^{n+1}-\bar u_\eps])}{Y^{n+1}}\sim \frac{M}{c(\eps)}, Z\in \Sigma$. Therefore 
choosing $M$ large enough, one sees that  $\cal M\geq \left[2K-\frac{Cc(\eps)}M\right]\Id \geq K\Id$ if $K>\frac{Cc(\eps)}M$.
Fixing $K\ge \max\left[\frac{Cc(\eps)}M, \sup |h|^{\frac1n}\right]$, where $h$ is defined by (\ref{def-rhs}) and 
choosing $r$ small enough such that 
(\ref{sml-grad}) holds we finally arrive at $\det\left[-D^2\bar u_\eps-\frac{G(x, \bar u_\eps, D\bar u_\eps)}{\eps\kap}\right]\geq |h|$ and the proof is complete.
\qed

%%%%%%%%%%%%%%%%%%%%%%%%%%%%%%%%%%%%%%%%%%%%%%%%%%%%%%%%%%%

%%%%%%%%%%%%%%%%%%%%%%%%%%%%%%%%%%%%%%%%%%%%%%%%%%%%%%%%%%%%
%                                                          %
%                SECTION                                   %
%                                                          %
%%%%%%%%%%%%%%%%%%%%%%%%%%%%%%%%%%%%%%%%%%%%%%%%%%%%%%%%%%%%
\section{$A$-type weak solutions and the Legendre-like transform}\label{sec-A-type}
In this section we are concerned with the second notion of weak solution to (\ref{problem}).
For $u\in \Wup(\U, \V)$ let us consider the mapping $\r_u: \U\to \Sigma$ defined as
\begin{equation*}
 \r_u(x)=\{Z\in \Sigma : \mbox{there \ is\ a \ supporting \ hyperboloid}\  H(\cdot,a, Z ) \ \mbox{ of} \ u\  \mbox{at} \ x\}.
\end{equation*}
Let  $E\subset \U$ be a Borel set and  put
\begin{equation*}
 \r_u(E)=\bigcup_{x\in E}\r_u(x).
\end{equation*}
Our primary goal is to prove that  $\r_u(E)$ is measurable with respect to the restriction of
$\H$ on $\Sigma$ for any Borel set $E\subset \U$.
That done, we can proceed as in \cite{K-ref} and establish that the set function $\alpha_{u, g}$
is $\sigma$-additive measure.

To take advantage of the geometric intuition coming from 
supporting hyperboloids of $u\in \Wup(\U, \V)$  it is convenient 
to define the Legendre-like transformation of $u$. We use the construction of smallest focal parameter
introduced by Xu-Jia Wang in \cite{W-96} (equation (1.15)).
Let $u\in \Wup(\U, \Sigma)$ and $Z\in \Sigma$ be a fixed point. Then the
smallest semi-axis among all hyperboloids $H(\cdot, a, Z)$  that stay above $u$ is 
$$a_0=\sup_{a\in  I(Z)}a, \quad
I(Z)=\{a>0: H(x, a, Z)\geq u(x)\ \mbox{in }\ \U \}.$$
Suppose that $H(\cdot, a_0, Z)$ touches $u$ at $x_0\in \U$ then
\begin{eqnarray*}
 u(x_0)=\psi(z)-a_0\eps-\sqrt{a^2_0+\frac{(x_0-z)^2}{\eps^2-1}}.\\\nonumber
 \end{eqnarray*}
 From here we can easily find that
 \begin{equation}\label{def-a-0}
 a_0=\frac{1}{\eps^2-1}\left\{\eps[u(x_0)-\psi(z)]-\sqrt{[u(x_0)-\psi(z)]^2+(x_0-z)^2}\right\}.
\end{equation}
  Alternatively, one can use the property that 
 the distance of a point $P$ on hyperboloid from lower focus $Z'$ is $\eps$ times the distance
  of $P$ from the hyperplane $\Pi_d=\{X\in \R^{n+1} :X^{n+1}=Z^{n+1}-a\eps-\frac a\eps\}$ (which in one dimensional case 
  is the directrix). Since by definition of hyperboloid $|PZ|-|PZ'|=2a$ and $|PZ'|=\eps\dist(P, \Pi_d)$ we infer 
$|PZ|=2a+\eps([\psi(z)-u(x_0)]-a\eps-\frac a\eps)=-a(\eps^2-1)+\eps([\psi(z)-u(x_0)]$ and 
(\ref{def-a-0}) follows.

%%%%%%%%%%%%%%%%%%%%%%%%%%%%%%%%%%%%%%%%%%%%%%%%%%%%%%%%%%%
%%%%%%%%%%%%%%%%%%%%%%%%%%%%%%%%%%%%%%%%%%%%%%%%%%%%%%%%%%%
\subsection{$A$-type weak solutions}

\begin{definition}
Let $u\in \Wup(\U, \Sigma)$ then
\begin{equation}\label{def-Legendre}
 v(z)=\inf_{x\in \U}\left\{\eps[\psi(z)-u(x)]-\sqrt{[u(x)-\psi(z)]^2+(x-z)^2}\right\}
\end{equation}
is called the Legendre-like transformation of $u$.
\end{definition}

If $\dist(\U,\V)>0$ and $\psi\in C^2$ then the function
$\mathscr L_x(z)=\eps[\psi(z)-u(x)]-\sqrt{[u(x)-\psi(z)]^2+(x-z)^2}$ is $C^2$-smooth for any fixed $x\in \U$.
Since $v$ is the upper envelope of $C^2$ smooth functions $\mathscr L_x, x\in \U$ ($x$ being the parameter) then $v(z)$ is semi-convex.
Next lemma gives an important characterization  of $v(z)$.

\begin{lemma}\label{lem-Leg-prop}
  Let $v$ be the Legendre-like transformation  of $u\in\Wup(\U,\Sigma)$. Then
\begin{itemize}
 \item[(i)] $v(z)=\eps[\psi(z)-u(x_0)]-\delta_u(x_0,z)$ if $Z=(z,\psi(z))\in\r_u(x_0)$ where $\delta_u(x, z)=\sqrt{[u(x)-\psi(z)]^2+(x-z)^2}$,
\item[(ii)] $v(z)$ is semi-concave.
\end{itemize}

\end{lemma}
\pr By definition $v(z)$ is locally bounded, non-negative, lower semi-continuous function.
Let $\delta_u(x, z)$ denote the distance between the points of
graph $\Ga_u$ and $\Sigma$. To check (i) we first observe that by definition of $v(z)$, see (\ref{def-Legendre}), we have
$v(z)\leq \eps[\psi(z)-u(x_0)]-\delta_u(x_0, z)$. If 
$v(z) < \eps[\psi(z)-u(x_0)]-\delta_u(x_0, z)$ it follows from (\ref{def-a-0}) and the discussion 
above that $H(\cdot, a_0, Z)$ is a supporting hyperboloid of $u$ at $x_0$, where $a_0=(\eps^2-1)^{-1}(\eps[\psi(z)-u(x_0)]-\delta_u(x_0, z))$
because  $Z\in\r_u(x_0)\subset\Sigma$.
On the other hand, there is a sequence $\{x_k\}$ in $\U$ 
such that $x_k\to  \bar x_0\in\ol{\U}$ and $\lim\limits_{x_k\to \bar x_0}\left(\eps[\psi(z)-u(x_k)]-\delta_u(x_k, z)\right)=v(z)$. 
Setting $\bar{a}_0=(\eps^2-1)^{-1}v(z)$ we conclude that $H(\cdot, \bar{a}_0, Z)$ is touching 
$\Ga_u$ at $\bar{x}_0$. By construction  $\bar {a}_0<a_0$ and it follows from 
confocal expansion of hyperboloids \ref{ssec-confocal} that $H(\cdot, \bar a_0, Z)>H(\cdot,a_0, Z)$ in $\U$.
But this inequality is in contradiction with the fact that $H(\cdot, a_0, Z)$ is a supporting 
hyperboloid of $u$ at $x_0$ and $H(\cdot, \bar a_0, Z)$ touches $\Gamma_u$ at $\bar x_0$ whilst staying above $\Ga_u$.

\smallskip

To prove (ii) we let $\mathscr L_{x_0}(y)=\eps[\psi(z)-u(x_0)]-\delta_u(x_0,y)$. Then
$$v(y)=\inf\limits_{x\in \U}\big\{\eps[\psi(z)-u(x)]-\delta_{u}(x,y)\big\}\leq \eps[\psi(z)-u(x_0)]-\delta_u(x_0,y)$$
which implies that $v(y)\leq\mathscr L_{x_0}(y)$ and $v(z)=\mathscr L_{x_0}(z)$,
where  $Z\in\mathscr R_u(x_0)$. We can regard $\mathscr L_{x_0}(y)$ as an upper
supporting function
of $v$ at $z$. Differentiating $\mathscr L_{x_0}$ twice in $z$ variable we 
see that 
$|D^2 \mathscr L_{x_0}(z)| \leq \frac{C}{(\dist(\U, \V))^3}$ 
for some tame constant $C>0$, consequently $v(z)-C|z|^2$ is concave for  large $C>0$.
\qed

The main result of this section is contained in the following
\begin{lemma}\label{Aleksandrov}
 Let $\cal S=\{Z\in \V : \mbox{such \ that}\ Z\in \r_u(x_1)\cap \r_u(x_2), x_1\ne x_2\}$.
 Then $\cal S$ has vanishing   surface measure on $\Sigma$.
\end{lemma}

\pr Let us show that if $Z\in \cal S $ then the Legendre-like transformation of $u$
is not differentiable at $Z$. This will suffice to
conclude the proof because by definition $v$ is semiconcave and hence by
Aleksandrov's theorem  twice differentiable almost everywhere. Let $v$ be the Legendre-like transformation of $u$, then by  Lemma \ref{lem-Leg-prop} for any $Z\in \r_u(x_0)$ at which $v(z)$ is differentiable 
 there holds
\begin{equation}\label{drv-v}
Dv(z)=\eps D\psi(z)-(y(x)+D\psi(z)y^{n+1}(x)).
\end{equation}
Indeed, $Dv(z)=\eps D\psi(z)-{\delta_u(x, z)}^{-1}\left[(z-x)+D\psi(z)(\psi(z)-u(x))\right]$. 
From the definition of stretch function $t$ it follows that
$(z-x, \psi(z)-u(x))=Y\delta_u(x, z)$ where $Y=(y, y^{n+1})$ is the unit direction of the refracted ray and 
(\ref{drv-v}) follows.
Consequently, if $x_1\ne x_2$ such that $\r_u(x_1)\cap\r_u(x_2)\ni Z$ then we must have

$$Dv(z)=-\frac{z-x_i+D\psi(z)(\psi(z)-x_i)}{\delta_u(x_i, z)}+\varepsilon D\psi(z), \quad i=1,2.
%Dv(z)=y(x_i)+D\psi(z)y^{n+1}(x_i)-\eps D\psi(z), \quad i=1,2.
$$
Equating the right hand sides for $i=1$ and $i=2$ we obtain
$$\frac{z-x_1+D\psi(z)(\psi(z)-x_1)}{\delta_u(x_1, z)}=\frac{z-x_2+D\psi(z)(\psi(z)-x_2)}{\delta_u(x_2, z)}$$
 With the aid of
this observation and (\ref{drv-v}) we can rewrite
the last line as follows 
\begin{eqnarray*}
 y_1+D\psi(z)y^{n+1}_1=y_2+D\psi(z)y_2^{n+1}\ {\rm in}\ \R^n\quad \Rightarrow\quad  Y_1+(D\psi(z), -1)y_1^{n+1}=Y_2+(D\psi(z), -1)y_2^{n+1},\ {\rm in}\ \R^{n+1}.
\end{eqnarray*}
The last identity implies that $Y_1-Y_2$ is collinear to the normal of $\Sigma$ at $Z$. Consequently, from
the assumption (\ref{vis-cond}) (see also (\ref{sml-cone})) we obtain that this is possible if and only if $Y_1=Y_2$.
%From this equality we can infer that $x_1=x_2$ which will be a contradiction.
Next, from $Y_1=Y_2$ we have $y_1=y_2$ and consequently we conclude that
\begin{equation}\label{fig-reciproc}
 \frac{z-x_1}{\delta_u(x_1, z)}=\frac{z-x_2}{\delta_u(x_2, z)}.
\end{equation}
 Taking the reciprocal of both sides
in the last identity and recalling the definition of the distance $\delta_u(x, z)$ one gets
\begin{eqnarray*}
 \frac{u(x_1)-\psi(z)}{|x_1-z|}&=&\frac{u(x_2)-\psi(z)}{|x_2-z|}
 \end{eqnarray*}
 yielding
 \begin{eqnarray*}
 %\\\nonumber
 u(x_1)&=&\psi(z)+\frac{|z-x_1|}{|z-x_2|}(u(x_2)-\psi(z))\\\nonumber
 &=&\psi(z)+\frac{\delta_u(x_1, z)}{\delta_u(x_2, z)}(u(x_2)-\psi(z)).
\end{eqnarray*}
On the other hand $Y^{n+1}_1=Y^{n+1}_2$ gives $u(x_1)-u(x_2)=\delta_u(x_2,z)-\delta_u(x_1,z)$ and hence 
combining this with the last equation yields
\begin{eqnarray*}
 \psi(z)\left[1-\frac{\delta_u(x_1,z)}{\delta_u(x_2, z)}\right]- u(x_2)\left[1-\frac{\delta_u(x_1,z)}{\delta_u(x_2, z)}\right]=\delta_u(x_1,z)-\delta_u(x_2,z).
\end{eqnarray*}
If $\delta_u(x_2,z)\not =\delta_u(x_1,z)$ then the last equality implies  $u(x_2)-\psi(z)=\sqrt{(u(x_2)-\psi(z))^2 +(z-x_2)^2}$. Hence
$x_2=z$ and by (\ref{fig-reciproc}) $x_1=x_2$, which is contradiction.
Thus we must have $\delta_u(x_2,z)=\delta_u(x_1,z)$ and  in view of (\ref{fig-reciproc}) this implies that
$x_1=x_2$, again contradicting our supposition. 
Therefore we infer that $v(z)$ cannot be differentiable at $z$.
By Rademacher's  theorem $v(z)$ is differentiable a.e. in $z$. Thus $\cal S$ has vanishing surface measure.
\qed

\smallskip

\begin{corollary}\label{cor-1}
For any $u\in \Wup(\U, \V)$ and any Borel subset $E\subset \U$ the set function 
\begin{equation}\label{alpha-mes}
\alpha_{u, g}(E)=\int_{\r_u(E)}gd\H
\end{equation} is a Radon measure.
\end{corollary}
\pr In order to show that $\alpha_{u,g}$ is Radon measure
it suffices to check that
$\widetilde{\mathscr  F}=\{E\subset\U : \r_u(E)\ \mbox{is measurable} \}$ is a
$\sigma-$algebra. This can be done  exactly in the same way as in the proof of Proposition
\ref{prop-S} c). It remains to recall that by Lemma \ref{lem-support},
$\widetilde{ \mathscr F}$
contains the closed sets.
\qed
%%%%

\begin{definition}\label{bndr-A-def}
 A function $u \in \Wup({\U, \V})$ is said to be $A$-type weak solution of (\ref{problem}) if $\int_{E}f(x)dx=\alpha_{u,g}(E)$ or any Borel 
set $E\subset \U$ and
\begin{equation}\label{bndr-A}
 \overline{\V}\subset\overline{\r_u(\U)}, \qquad |\{x\in \U : \r_u(x)\not \subset \V \}|=0
\end{equation}
\end{definition}
This definition is natural, stating that  the target domain $\V$ is covered by the refracted
rays and the endpoints of those rays that after refraction do not strike
$\V$ constitute a null set on $\U$. We shall establish the 
existence of $A$-type weak solution in the next section.

In closing this section we state the weak convergence result for the
$\alpha$-measures, see Corollary \ref{cor-1}.
\begin{lemma}\label{lem-weakconv-A}
 Let $u_k$ be a sequence of $A$-type weak solutions and $\alpha_k$
  is the corresponding measure, defined by (\ref{alpha-mes}). If $u_k\rightarrow u$ uniformly on
  compact subsets of ${\U}$ then $u$ is $R-$concave and $\alpha_k$ weakly
  converges to $\alpha_{u,g}$.
\end{lemma}
The proof is very similar to that of Lemma \ref{lem-weakconv-B} (modulo minor adjustments) and hence omitted.

\begin{remark}
The Legendre like transformation (\ref{def-Legendre}) can be used to reformulate (\ref{problem}) as a nonlinear optimization
problem, following the method set out in \cite{Liu}. 
\end{remark}
%%%%%%%%%%%%%%%%%%%%%%%%%%%%%%%%%%%%%%%%%%%%%%%%%%%%%%%%%%%
%%%%%%%%%%%%%%%%%%%%%%%%%%%%%%%%%%%%%%%%%%%%%%%%%%%%%%%%%%%

%%%%%%%%%%%%%%%%%%%%%%%%%%%%%%%%%%%%%%%%%%%%%%%%%%%%%%%%%%%%
%                                                          %
%                SECTION                                   %
%                                                          %
%%%%%%%%%%%%%%%%%%%%%%%%%%%%%%%%%%%%%%%%%%%%%%%%%%%%%%%%%%%%
\section{Comparing $A $  and $ B$ type weak solutions: Proof of Theorem C3-4}
In this section we prove the equivalence of 
$A$ and $B$ type weak solutions under some conditions. These 
results are known for in $\R^n$ for the sub-differential \cite{C92}, \cite{Urbas}. Let $\phi:\R^{N}\to\R^n$ be a Borel mapping and $\mu(\R^N)=\nu(\R^n)<\infty$ with
$\mu, \nu$ being two Radon measure on  $\R^N$ and $\R^n$, respectively. Then
$\phi$ induces a (push-forward) measure on $\R^n$ defined by $\phi_{\#}\mu(E)=\mu(\phi^{-1}(E))$ for Borel subsets
$E\subset \R^n$. We say that a Borel mapping
$\phi$ {\it measure preserving} if
\begin{equation}\label{push-f}
 \phi_\#\mu(E)=\nu(E)\qquad \mbox{for any Borel set}\ E\subset \R^n.
\end{equation}
By the change of variables formula (\ref{push-f}) can be rewritten
in the following equivalent form
\begin{equation}\label{push-int}
 \int h(\phi(x))d\mu=\int h(y)d\nu, \qquad \forall h\in C(\R^n),
\end{equation}
see \cite{Brenier}.

\begin{remark}\label{rem-push}
If $u\in\W^+(\U,\Sigma)$ is the $B$-type solution of (\ref{problem}), the existence of
which is established in Section \ref{sec-weak-B},
then taking  $\phi(Z)=\s_u(Z), N=n+1$, $d\mu=gd\H$, and $\nu$ being the Lebesgue measure
one immediately observes that $\s_u$ is measure preserving in the sense of (\ref{push-f}) or
(\ref{push-int}).
\end{remark}
\begin{lemma}\label{lem-hull}
 If $\r_u(x)\subset\V$ for a.e. $x\in \U$ then $\r_u(E)\subset\mbox{\rm Hull}(\V)$, where
$\mbox{\rm Hull}(\V)$ is the R-convex hull of $\V$ defined as the smallest
R-convex subset of $\Sigma$ containing $\V$.
\end{lemma}

\pr We only have to consider the points where $u$ is non-differentiable.
Let $u$ be non-differentiable at $x_0\in \U$ and suppose that
$\gamma_1, \gamma_2$ are the normals of two supporting planes of $u$
at $x_0$. The ray with endpoint $x_0$ after reflection will lie in the
reflector cone $\cal C_{\xi_0, \gamma_1, \gamma_2}$, with $\xi_0=(x_0, u(x_0))$ and
the reflected ray will strike $\mbox{Hull}(\V)$, because
$\cal C_{\xi_0, \gamma_1, \gamma_2}\cap \mbox{Hull}(\V)$ is connected.
Considering all normals of supporting planes at $x_0$ we obtain the desired result.
\qed

\begin{proposition}\label{lem-A=B}
  Let $\Sigma$ be $R-$convex with respect to $Q_m=\U\times(0,m),m>0$
  and the densities $f, g$ are positive. Then $B$-type weak solution  is also of $A$-type.
\end{proposition} 
% \td{in this prop $f$ and $g$ must be positive}
\pr We split the proof into three parts.

%\noindent
\mbox{1)} First we show that for  any compact $K_1\subset \U$ there holds
$\int_{K_2}gd\H\geq \int_{K_1}f(x)dx$ with $K_2=\r_u(K_1)$.
In other words the $B$-type solution is $A$-type subsolution. It is worthwhile to point out that
for the proof of
this inequality we don't need $\V$ to be $R-$convex. Take $\eta\in C(\Sigma)$ such that
$\eta\equiv1$ on $K_2\subset \Sigma$ and $0\leq \eta\leq 1$. From
(\ref{push-int}) we see that
\begin{equation*}
 \int_{\V}\eta gd\H=\int_\U\eta(\r_u(x))f(x)dx\geq \int_{K_1}f(x)dx.
\end{equation*}
Letting $\eta$ to decrease to the characteristic function of
$K_2$, $h\downarrow \chi_{K_2}$ we infer
\begin{equation}\label{ineq-1}
 \int_{K_2}gd\H\geq \int_{K_1}f(x)dx.
\end{equation}
Notice
by Corollary  \ref{cor-1} the measure $\alpha_{u,g}$ is Borel regular, therefore in the last inequality
$K_1$ can be replaced by any Borel subset of $\U$.
As a result we conclude from (\ref{ineq-1}) that
\begin{equation}\label{meas-cont-H}
 \mbox{if}\ \H(\r_u(E))=0\ \mbox{then}\ |E|=0.
\end{equation}
%\td{this is where we need to have $f>0$}

\smallskip

%\noindent
\mbox{2)} Next, we prove the converse estimate of (\ref{ineq-1}). Here we will utilize the $R-$convexity of $\V$.
Take any compact $K_1\in\U$ and apply Lemma  \ref{Aleksandrov} to conclude
$\H\big(\r_u(K_1)\cap\r_u(\U\setminus K_1)\big)=0$. Let us show that
\begin{equation}\label{xarn}
|\r_u^{-1}(\r_u(K_1))\setminus K_1|=0
\end{equation}
where $\r_u^{-1}(\r_u(K_1))$ is the pre-image of $\r_u(K_1)$.
Denote $E=\r_u^{-1}(\r_u(K_1))$ and $G=K_1$.
If $\H(E\setminus G)=0$ then in view of (\ref{meas-cont-H})  we obtain (\ref{xarn}).
Indeed, form the identity (\ref{compl-ident}) it follows that
\begin{eqnarray}
 |\r_u(E\setminus G)|&=&\left|[\r_u(E)\setminus \r_u(G)]\bigcup[\r_u(E\setminus G)\cap\r_u(G)]\right|\\\nonumber
&=&|\r_u(E\setminus G)\cap\r_u(G)|\\\nonumber
&=&0
\end{eqnarray}
where to get the last line we used the definitions of $E$ and $G$ in order to obtain
$\r_u(E)\setminus \r_u(G)=\r_u(K_1)\setminus \r_u(K_1)=\emptyset$ and
Lemma \ref{Aleksandrov}. Thus (\ref{meas-cont-H}) implies
$0=|E\setminus G|=|\r_u^{-1}(\r_u(K_1))\setminus K_1|$.

\smallskip

%Now we are ready to finish the proof and establish the converse of the inequality (\ref{ineq-1}).
Let $h\in C(\Sigma)$ such that $0\leq h\leq 1$ and $h\geq \chi_{\r_u(K_1)}$.
Since $\V$ is $R-$convex it follows that
$\r_u(K_1)\subset \hbox{Hull}\V$, see Lemma \ref{lem-hull}. If $u$ is a
$B$-type weak solution then (\ref{push-int}) holds, see Remark \ref{rem-push}. Therefore

\begin{eqnarray*}
 \int_\U\eta(\r_u(x))f(x)dx&=&\int_\V\eta gd\H\\\nonumber
&=&\int_{\textrm{Hull}(\V)}\eta gd\H\\\nonumber
&\geq&\int_{\r_u(K_1)}gd\H.
\end{eqnarray*}

Letting  $\eta\rightarrow 0$ on compact subsets of $\V\setminus\r_u(K_1)$,
it follows that $\eta(\r_u(x))$ uniformly converges to zero one the
compact subsets of $\U\setminus \r_u^{-1}(\r_u(K_1))$. Consequently
$$\int_{\r_u(K_1)}gd\H\leq \int_\U\eta(\r_u(x))f(x)dx\longrightarrow \int_{\r^{-1}(\r_u(K_1))}f(x)dx=\int_{K_1}f(x)dx$$
where the last line follows from (\ref{xarn}). This implies that $u$ is a supersolution.

\smallskip

\mbox{3)} It remains to check that $u$ verifies the boundary condition (\ref{bndr-A}).
Suppose that there is  $Z_0\in \overline{\V}$ such that $Z_0\not \in \overline{\r_u(\U)}$.
Since $u$ is of $B$-type, it follows that $\s_u(\V)=\U$ implying
$x_0\in \s_u(Z_0)$ in other words, there is a supporting hyperboloid
$H(x,a_0, Z_0)$ at $x_0$. Thus $Z_0\in\r_u(x_0)$which yields
$\overline{\V}\subset \overline{\r_u(\U)}$.
From energy balance condition we have
$$\int_{\overline{\r_u(\U)}}gd\H=\int_\U f(x)dx=\int_{\V}g d\H\quad \Rightarrow\quad  \int_{\overline{\r_u(\U)}\setminus \V}g=0.$$
This yields $|\{x\in \U : \r_u(x)\not \subset \V\}|=0$ for $f, g>0$.
\qed 
%% \td{This is where we use $g>0$}

\begin{remark}\label{rem-set-eqv}
We always have
$\ol{\V}\subset \r_u(\U)$, however if in addition $\Sigma$ is $R-$convex then it follows that
$\r_u(\U)\subset \V$. Thus we get the equality $\overline{\r_u(\U)}=\overline{\V}$
for R-convex $\V$.
\end{remark}
%%%%%%%%%%%%%%%%%%%%%%%%%%%%%%%%%%%%%%%%%%%%%%%%%%%%%%%%%%%
%%%%%%%%%%%%%%%%%%%%%%%%%%%%%%%%%%%%%%%%%%%%%%%%%%%%%%%%%%%

%
%%
%%%
%%%%
%%%%%
%%%%%%
%%%%%%%
%%%%%%%%  Existence of weak solutions
%%%%%%%%%
%%%%%%%%%%
%%%%%%%%%%%
%%%%%%%%%%%%
%%%%%%%%%%%%%

\subsection{Existence of $A$-type weak solutions: Proof of Theorem C4}

Suppose that $\V\subset \Sigma$ and let  $\mbox{\rm Hull}(\V)$ be the
$R-$convex hull of $\V$. For small $\delta, \delta'>0$ we consider
\begin{equation}
 g_\delta(Z)=
\left
\{\begin{array}{lll}
g(Z)-\delta & \mbox{if}\ Z\in \V\\
\delta' &\mbox{if} \ Z\in\mbox{Hull}(\V)\setminus \V
\end{array}
\right.
\end{equation}
where we choose $\delta, \delta'$ so that $g_\delta$ satisfies the energy balance
condition (\ref{glob-balance}).
By Proposition \ref{prop-B-exst} for each $g_\delta$ there is a $B$-type weak solution which according to Proposition
\ref{lem-A=B} is also of $A$-type. Moreover, from Remark \ref{rem-set-eqv} we infer
\begin{equation}
 \overline{\r_{u_\delta}(\U)}=\overline{\V}.
\end{equation}
 Sending $\delta\rightarrow 0$ we obtain from Lemma \ref{lem-weakconv-A}
that $u_\delta\rightarrow u$ and $u$ is an
$A$-type solution, i.e. (\ref{alpha-mes}) is satisfied, and

\begin{equation}
 \overline\V\subset \overline{\r_u(\U)}.
\end{equation}
Since $u$ is second order differentiable a.e. in $\U$
it follows that $\r_u$ is defined for a.e. $x\in \U$.
Finally we want to show that   $|S|=0$ where
$S=\{x\in \U : \exists Z\in \r_u(x)\ \mbox{such\ that}\ Z\in \r_u(\U)\setminus \V\}$.
Indeed, from energy balance condition (\ref{glob-balance}) we have
\begin{eqnarray*}
 \int_Sf(x)dx&=&\int_{\U}f(x)dx-\int_{\U\setminus S}f(x)dx=\\\nonumber
&=&\int_\U f(x)dx-\int_\V g d\H=0.
\end{eqnarray*}
Since $f>0$ we conclude that $|S|=0$ and hence (\ref{bndr-A}) holds and
$u$ is a weak $A$-type weak solution of (\ref{problem}). \qed

\begin{remark}\label{rem-zzbg}
As the proof of Proposition \ref{lem-A=B} exhibited if $\V$ is $R-$convex then $S=\emptyset$. If
$S\not=\emptyset$ then $u$ is only Lipschitz continuous. Therefore 
if $\V$ is not R-convex
then $u$ may not be $C^1$ smooth, see Introduction. It is worthwhile to point out that
even if $S=\emptyset$ then $u$ may not be $C^1$, and hence further
assumptions must be imposed to assure the smoothness of $u$.
\end{remark}

%%%%%%%%%%%%%%%%%%%%%%%%%%%%%%%%%%%%%%%%%%%%%%%%%%%%%%%%%%%%
%%%%%%%%%%%%%%%%%%%%%%%%%%%%%%%%%%%%%%%%%%%%%%%%%%%%%%%%%%%%

%%%%%%%%%%%%%%%%%%%%%%%%%%%%%%%%%%%%%%%%%%%%%%%%%%%%%%%%%%%%%%%%%%%%
%                                                                  %
%                                                                  %
%                    SECTION                              %
%                                                                  %
%                                                                  %
%%%%%%%%%%%%%%%%%%%%%%%%%%%%%%%%%%%%%%%%%%%%%%%%%%%%%%%%%%%%%%%%%%%%
\medskip

\section{Dirichlet's problem}\label{sec-Dirichlet}

This section concerns the Dirichlet problem for $A$-type weak solutions. We formally rewrite the equation 
(\ref{M-A-short+}) below
\begin{equation}\label{diff-op-A}
 \cal F[u](x)=\frac{f(x)}{g\circ \r_u(x)}, \qquad x\in \U,
\end{equation}
where for $u\in C^{2}(\U)$, $\cal F[u](x)$ is the determinant of the Jacobian matrix of $\r_u$. For non-smooth solutions we give the following definition. 

\smallskip

\begin{definition}\label{def-AS+}
A function $u \in \Wup(\U,\Sigma)$ is said to be a
weak A-subsolution of (\ref{diff-op-A})
if for any Borel set $E$
\begin{equation}
 \int_{\r_u(E)}gd\H\geq \int_Ef(x)dx.
\end{equation}
 If $\alpha_{u,g}(E)=\int_Ef(x)dx$ then we say that $u$ is a weak  A-solution.
The class of all generalized A-subsolutions is denoted by ${\cal {AS}}^{+}(\U)$.
\end{definition}

For smooth and bounded $D\subset\Sigma$ and smooth function $\phi:\ol D\to \R$  let us consdier the Dirichlet problem
\begin{eqnarray}\label{BVP-A-1}
 \left\{
\begin{array}{lll}
 \displaystyle\cal F[u](x)=\frac{f(x)}{g\circ \r_u(x)} &\qquad x\in D,\\
u=\phi &\qquad x\in \p D.
\end{array}
\right.
\end{eqnarray}
Our main objective here is to prove  the existence and uniqueness of $A$-type weak solution to (\ref{BVP-A-1}) for a smooth boundary data. In fact, for our purposes it suffices to consider the case where $D$ is a ball of small radius.
At this point we first we establish the following comparison principle. 

%%%%%%%%%%%%%%%%%%%%%%%%%%%%%%%%%%%%%%%%%%%%%%%%%%%%%%%%%%%
%%%%%%%%%%%%%%%%%%%%%%%%%%%%%%%%%%%%%%%%%%%%%%%%%%%%%%%%%%%
\smallskip

\begin{proposition}\label{prop-compar}%[Comparison principle]
 Let $u_i$ be weak solutions of (\ref{diff-op-A}) in $\U$ with
$f=f_i, i=1,2$, where $\Omega\subset \Pi$ is a smooth, bounded domain and conditions in Theorem C hold. 
Suppose that $\r_{u_1}(\Omega)\subset\overline\Sigma$, $f_1
< f_2$ in $\Omega$ and $u_1\le u_2$ on $\p \Omega$. If  $\Gamma_{1}$, the
graph of $u_1$, lies in the region $\cal D$ then we have $u_1\le u_2$ in $\Omega$.
\end{proposition}

\pr Suppose that $\Omega_1=\{x\in \Omega : u_1(x)>u_2(x)\}$ is not empty.
Let $x_0\in \Omega_1$ and $H(x, a_0, Z_0), Z_0\in \Sigma$ is a supporting
hyperboloid of $u_2$ at $x_0$, i.e. $Z_0\in \r_{u_2}(x)$. Let us show that there is $\bar x$ such that 
$Z_0\in \r_{u_1}(\bar x)$. Observe that by \eqref{vis-cond-B} the hyperboloids $H(x, a_0+s, Z_0)$ 
stay above $\U$ for all $s<0$ such that $\bar a_0>a_0+s$, see  \eqref{max-a} for the definition of $\bar a_0$.
Notice that if $a_0+s>0$ is very small then the corresponding hyperboloid $H(x, a_0+s, Z_0)$ is very close 
to the asymptotic cone of the hyperboloid, with vertex very close to the fixed focus $Z_0$.

\smallskip 

Consequently, from the confocal expansion of hyperboloids
(see subsection \ref{ssec-confocal}) we infer that there is $s_0<0$ such that $H(x, a_0+s, Z_0)$ is a
supporting hyperboloid of $u_1$ at an interior point $x_1\in\Omega_1$. Thus
$H(x, a_0+s_0, Z_0)$ is a local supporting hyperboloid of $u_1$. Since
$\Gamma_{u_1}$ is in the regularity domain $\cal D$, where (\ref{vis-cond})-(\ref{A3-cond00}) are fulfilled, we can apply
Lemma \ref{loc-glob} to conclude  that
$H(x, a_0+s, Z_0)$ is also a global supporting hyperboloid  of $u_1$. Hence $Z_0\in \r_{u_1}(x_1)$. Therefore
$$\r_{u_2}(\Omega_1)\subset \r_{u_1}(\Omega_1)$$
implying
\begin{eqnarray*}
 \int_{\Omega_1}f_1dx<\int_{\Omega_1}f_2dx=\int_{\r_{u_2}(\Omega_1)}gd\H\leq
\int_{\r_{u_1}(\Omega_1)}gd\H=\int_{\Omega_1}f_1dx
\end{eqnarray*}
which gives a contradiction. Thus $\Omega_1=\emptyset.$
\qed

\smallskip

\subsection{Discrete Dirichlet problem}
To outline our next two steps, we note that for the classical Monge-Amp\`ere equation 
the standard way of proving the existence of globally smooth solutions to Dirichlet's problem with, say, $\phi\in C^4(\ol D)$
is to employ  the continuity method
combined with standard mollification argument, see \cite{Pog-2D}. 
Moreover, in this argument $\phi$ must be a subsolution. 
In order to tailor a similar proof for (\ref{BVP-A-1}) we will mollify our weak $A$-solution, add $K(r^2-|x-x_0|^2), K\gg 1$ and consider 
its restriction to $B_r(x_0)\subset \U$, a ball with sufficiently small radius. Such function turns out to be classical 
subsolution for  some large $K$ and small $r>0$. Consequently, from continuity method  one can obtain  existence of a smooth solution to Dirichlet's problem in $B_r(x_0)$. Finally employing the known $C^2$ a priori estimates and comparison principle, Proposition \ref{prop-compar},  the proof of Theorem D will follow, see Section \ref{sec-pr-Th2} for more details.
Our approach most closely follows that proposed by Xu-Jia Wang \cite{W-96}.

\smallskip

Let $\{b_i\}\subset \p D$ be a sequence of  points on the boundary of $D$ and
$\{a_i\}\subset D$. Let  $A_N=\{a_1, \dots, a_N\}$ and
$B_N=\{b_1, \dots, b_N\}\subset \partial D$, for each fixed $N\in \mathbb N$.
Furthermore, let  $\nu_k(x)$ be atomic measures supported at $a_k, 1\leq k\leq N$ and let 
\begin{equation}\label{eq-discret}
 \mathcal F [v](x)= \nu_k(x)\frac{f(x)}{g\circ\r_v(x)}.
\end{equation}

\begin{proposition}
 Let $\underline{u}\in \Wup^0(D, \Sigma)$ be a polyhedral subsolution of (\ref{eq-discret}),
 i.e. $\mathcal F[\underline u](x)\geq \nu_k(x)\frac{f(x)}{g\circ \r_{\underline u}(x)}$ at  $a_k\in A_N$.
  Then there is a unique $A$-type weak solution
to (\ref{eq-discret}) verifying the boundary condition $u=\underline{u}$ on $B_N$.
\end{proposition}

\pr We want to construct a sequence of subsolutions $\{u_m\}_{m=0}^\infty$ converging to the solution of discrete problem. Set $u_0=\underline{u}$ and define
$u_1$ such that $u_1\leq u_0$ in $A_N$, $u_1(b_i)=u_0(b_i), b_i\in B_N$ and
$\alpha_{u_1,g}(a_i)\leq \alpha_{u_0,g}(a_i)$ for $a_i\in A_N$.
It is convenient to introduce the  class of hyperboloids
\begin{eqnarray*}
  \Phi_{0, \delta}(a_1)=\left\{P\in \mathbb H^+_0(D, \Sigma): \begin{array}{lll}
    H(a_i)\geq u_0(a_i), i\not=1,\\
  H(a_1)\geq u_0(a_1)-\delta, \\  H(b_j)\geq u_0(b_j), 1\leq j\leq N
  \end{array}\right\}
\end{eqnarray*}
for $\delta>0$ and let
\begin{equation*}
  T^\delta_1u_0=\inf_{H\in \Phi_{0,\delta}(a_1)}H.
\end{equation*}
Let $\delta_1>0$ be the largest $\delta$ for which $T^{\delta_1}_1u_0$ is a subsolution to
(\ref{eq-discret}) on $A_N$. Consequently, by setting  $u_{0,1}=T^{\delta_1}_1$ we can proceed by induction and define the $k$-th class 
\begin{eqnarray*}
  \Phi_{0, \delta}(a_k)=\left\{H\in \mathbb H^+_0(D, \Sigma): \begin{array}{lll}
    H(a_i)\geq u_{0,k-1}(a_i), i\not=k,\\
  H(a_k)\geq u_{0,,k-1}(a_k)-\delta, \\  H(b_j)\geq u_{0,,k-1}(b_j), 1\leq j\leq N
  \end{array}\right\}
\end{eqnarray*}
and take $T^\delta_k u_0=\inf\limits_{H\in \Phi_{0,\delta}(a_k)}H$.
Therefore, one  can successively construct the functions $u_{0,k}=T^{\delta_k}_ku_{0,k-1}$ where
$\delta_k>0$ is the largest number for which $T^{\delta}_ku_{0,k-1}$ is a subsolution to
(\ref{eq-discret}) in $A_N$.
Taking the second subsolution in the approximating sequence to be 
$u_2(x)\eqdef T^{\delta_N}_Nu_{0,N-1}$ we get, by construction, that 
$\alpha_{u_0,g}(a_i)\leq \alpha_{u_1,g}(a_i)$, since we have the inclusions $\Phi_{l, \delta}(a_k)\subset \Phi_{l+1, \delta}(a_k)$ at $a_k$ as we proceed. Therefore we have a sequence of solutions $u_m$ to the
Dirichlet problem in $A_N$ such that
\begin{eqnarray*}
  \alpha_{u_{m},g}(a_i)\leq \alpha_{u_{m-1},g}(a_i),\\
  u_m(a_i)\leq u_{m-1}(a_i),\\
  u_m(b_i)= u_{m-1}(b_i).
\end{eqnarray*}
The first two inequalities are obvious.
As for the boundary condition we note that $u_0(b_i)\leq u_1(b_i)$ by construction. If
$u_0(b_i)< u_1(b_i)$ then by taking $\min[H_i(x),u_1(x)]$, where
$H_i(x)\in\mathbb H^+_0(D,\Sigma)$ is a supporting hyperboloid
of $u_0$ at $b_i$ we see that $\min[H_i(x),u_1(x)]$ belongs to the corresponding
$\Phi$ class.  Thus $u_0(b_i)= u_1(b_i)$.

From Lemma \ref{lem-support} we conclude that $u\in \Wup(D,\Sigma)$
and in view of Lemma \ref{lem-weakconv-A} $\alpha_{u_m,g}\rightharpoonup \alpha_{u,g}$
weakly. Thus $u=\lim\limits_{m\rightarrow \infty} u_m$ is a solution to the discrete
problem in $A_N$ with $u(b_i)=\underline{u}(b_i), b_i\in B_N$.
\qed

\medskip
\subsection{General case}
Perron's method, used in the proof of above proposition, can be strengthened in order to establish the solvability of  the general Dirichlet problem.
To do so we  take $\{a_i\}_{i=1}^\infty\subset D$ and $\{b_i\}_{i=1}^\infty\subset \p D$
to be dense subsets and $A_N=\{a_1, \dots, a_N\}\subset D, B_N=\{b_1, \dots, b_N\}\subset \partial D$.
\begin{proposition}\label{prop-gen_Dirichlet}
Let $\underline{u}\in\mathcal{AS}^+(D,\Sigma)$.
 Then there exists a unique weak  solution $u$ to the Dirichlet problem
\begin{eqnarray}
  \left\{
  \begin{array}{lll}
  \displaystyle \mathcal F [u] = \frac{f(x)}{g\circ\r_u(x)} & \mbox{in} \ D,\\
  \displaystyle u(x)=\underline{u}(x) & \mbox{on} \ \p D.
  \end{array}
  \right.
\end{eqnarray}
\end{proposition}

\pr For $\delta>0$ we denote $D_{\delta}=\{x\in D : \mbox{dist}(x,\p D)>\delta\}$ and
take $\eta(x)$ to be a smooth function such that $0\leq \eta(x)\leq 1$, $\eta\equiv 1$ in
$D_{2\delta}$ and $\eta\equiv 0$ in $D\setminus D_{\delta}$.
Consider the equation
\begin{equation}\label{xarn-hav}
 \mathcal F[v](x)=\nu_k(x)J(v(x))\eta_\delta(x)\frac{f(x)-\delta}{g\circ\r_v(x)}
\end{equation}
where $\nu_k(x)$ is a positive measure supported at $a_k\in A_N$ and
\begin{equation}
 J(t)=\left\{
\begin{array}{lll}
 1 & \quad \mbox{if} \quad \ 0\leq t \leq \sup\limits_D\underline{u},\\
\frac{2\sup\limits_D\underline{u}-t}{\sup\limits_D\underline{u}} & \quad \mbox{if} \quad \ \sup\limits_D\underline{u}
\leq t \leq 2\sup\limits_D\underline{u},\\
0 \ & \quad \mbox{if} \quad \ t> 2\sup\limits_D\underline{u}.
\end{array}
\right.
\end{equation}
Consider the class
\begin{equation}
 \W^+_{N,\underline{u}}=\left\{v\in \Wup^0(D, \Sigma) :
\mathcal F[v]\geq \nu_kJ(v)\eta_\delta\frac{f-\delta}{g\circ \r_u}\
\mbox{and} \ v\geq \underline{u}\ \mbox{on}\ B_N\right\}.
\end{equation}
Clearly $\W^+_{N,\underline{u}}$ is not empty
since $H(\cdot, a, Z)$ is in
this class if $a>0$ is sufficiently small.
We claim that if $v_{N, \delta}=\inf\limits_{\W^+_{N,\underline{u}}}v$ then 
$v_{N, \delta}$ solves (\ref{diff-op-A}) in the sense of Definition
\ref{def-AS+} and $v_{N, \delta}(b_i)=\underline{u}(b_i), b_i\in B_N$.

It is easy to see that $\alpha_{v_{N,\delta},g}(a_k)= v_k(a_k)J(v_{N,\delta})\eta_\delta(a_k)
\left(f(a_k)-\delta\right)$.
Indeed, if $v_{N, \delta}$ is a strict subsolution at $a_i$, i.e.
for some $a_i$ we have
$\alpha_{v_{N,\delta}, g}(a_i)>
v_k(a_i)J(v_{N,\delta})\eta_\delta(a_i)(f(a_i)-\delta)$, then
we can push $\Gamma_{v_{N,\delta}}$ downward by some $\delta>0$, decreasing the
$\alpha$ measure at $a_i$, which, however, will be in contradiction
with the definition of $v_{N,\delta}$.
Thus  $v_{N,\delta}$ is a solution of the equation (\ref{xarn-hav}).

Next, we check the boundary condition. Choose $H_i\in \mathbb H^+_0(\U, \Sigma)$
such that $H_i>v_{\delta}$ in
$\U_\delta$ and passes through $(b_i, \underline{u}(b_i))$. Such $H_i$
exists because by construction
$v_{N, \delta}(a_i)\leq \underline{u}(a_i)$ and $\delta>0$.

For $\widetilde H_i=\min[H_i, v_{N,\delta}]$, by construction,  we see that
$\mathcal F [\widetilde H_i]\geq \nu_kJ(\widetilde H_i)\eta_\delta\frac{f-\delta}{g\circ\r_{\widetilde H_i}}$
at $a_i$. Thus $\widetilde H_i\in \W^+_{N,\underline{u}}$. Hence

$$v_{N,\delta}(b_i)=\inf_{H\in \W^+_{N,\underline{u}}} H(b_i)\leq \widetilde H_i(b_i)
=\underline{u}(b_i).$$

Now the desired solution can be
obtained via a standard compactness argument that utilizes the
estimates of  Lemma \ref{lem-lips} and Lemma \ref{lem-weakconv-A}.
More precisely, for fixed $\delta>0$ we send
$N\rightarrow \infty$ and obtain a function $v_\delta$ that solves the
equation $\mathcal F[v_\delta]= J(v_\delta)\eta_\delta\frac{f-\delta}{g\circ\r_{v_\delta}}$.
To show that $v_\delta=\underline{u}$ on $\p D$ we take $x_0\in \p D$
and again use the comparison with
$\min[H_0, v_\delta]$ for a suitable $H_0\in \mathbb H^+_0(\U, \Sigma)$ such that
$H_0(x_0)=\underline{u}(x_0)$. Thus, from Proposition \ref{prop-compar}
we conclude that $v_\delta\leq \underline{u}$ in
$D$. Finally sending  $\delta\downarrow 0$ and employing the estimate
of Lemmas \ref{lem-lips} and  \ref{lem-weakconv-A}
we arrive at desired result.\qed

\medskip
%%%%%%%%%%%%%%%%%%%%%%%%%%%%%%%%%%%%%%%%%%%%%%%%%%%%%%%%%%%%%%%%%%%%
%                                                                  %
%                                                                  %
%                    SECTION                                       %
%                                                                  %
%                                                                  %
%%%%%%%%%%%%%%%%%%%%%%%%%%%%%%%%%%%%%%%%%%%%%%%%%%%%%%%%%%%%%%%%%%%%
\section{Proof of Theorem D}\label{sec-pr-Th2}

To fix the ideas we assume that $x_0=0 \in \U$ and $B_r=B_r(0)\subset \U$.
Let $u^\pm_{s,\delta}$ be the solutions to
\begin{equation}\label{verj-eqn}
 \left\{
 \begin{array}{ccc}
  \cal F[ u^\pm_{s, \delta}]=\frac{f\pm\delta}{h g\circ Z_{u^\pm_{s, \delta}}}& {\rm in} \ B_{r}\\
  u^\pm_{s,\delta}=\widetilde u_{s} & {\rm on} \ \partial B_{r}
 \end{array}
 \right.
\end{equation}
where $\widetilde u_s=u_s+K(r^2-|x|^2)$, $K>0$ and $u_s$ is a mollification of the weak solution $u$.
By Lemma \ref{lem-apprx} $\widetilde u_s$ is a subsolution (for appropriate choice of constants $K$ and $r$) and hence by
Proposition \ref{prop-gen_Dirichlet} the weak
solution to Dirichlet's problem exists. Note that for the Dirichlet
problem we have to consider the modified  receiver $\widetilde \Sigma$ to guarantee that 
$\widetilde u_s$ is admissible, see Lemma \ref{lem-apprx}.
In order to show the existence of smooth solutions we apply the continuity method:
%% continuity method
Let $\underline{w}\in \cal {AS}^+(B_r, \widetilde\Sigma)\cap C^{\infty}(B_r)$ and for $t\in[0,1]$ consider the solutions
to the Dirichlet problem
\begin{equation}\label{t-Dirichlet}
\left\{
\begin{array}{lll}
  \cal F[w^t]=t\frac{f}{hg\circ Z_w} +(1-t)\cal F[\underline w] \ \ & \mbox{in } \ B_r,\\
w^t=\underline w &\mbox{on} \ \p B_r,
\end{array}
\right.
\end{equation}
where $h$ is given by (\ref{def-rhs}).
%\td{ h is given by (\ref{def-rhs})}
Using the implicit function theorem, see \cite{Trud-Tokyo} Theorem 5.1,
we find that
the set of $t$'s for which (\ref{t-Dirichlet}) is solvable
is open.

Once $C^{1,1}$ global a priori estimates were established
in $\overline{B_r}$ then one can deduce that the set of such $t$'s is also closed.
Recall that if
$\p\Omega\in C^3, u\in C^4(\Omega)\cap C^3(\overline\Omega)$ and $\underline{u}\in C^4$
then from global $C^{1,1}$ estimates and  the elliptic regularity theory
we obtain that $w\in C^{2, \alpha}(\overline\Omega)$.
Therefore the existence of smooth solutions $u_{s, \delta}^\pm$ will follow once
we establish the global $C^{1,1}$ estimate for $w$.
The latter follows from 
\cite{Turd-chin-arxiv} and Theorem B.

Summarizing, we have that $u^\pm_{s, \delta}$ remain locally uniformly smooth in $B_r$. Letting $s\to 0$ and applying
the comparison principle (see Proposition \ref{prop-compar})
we have that $u_{0,\delta}^-\leq u \leq u^+_{0, \delta}$ and $u_{0,\delta}^\pm=u$ on $\p B_r$. It follows from
the a priori estimates in \cite{Turd-chin-arxiv} and  \cite{MTW}
that $u^\pm_{0, \delta}$ are locally uniformly $C^2$
in $B_r$ for any small $\delta>0$. After sending $\delta\to 0$ we will conclude the proof of Theorem D.

%\smallskip

%%%%%%%%%%%%%%%%%%%%%%%%%%%%%%%%%%%%%%%%%%%%%%%%%%%%%%%%%%%%
%                                                          %
%                SECTION                                   %
%                                                          %
%%%%%%%%%%%%%%%%%%%%%%%%%%%%%%%%%%%%%%%%%%%%%%%%%%%%%%%%%%%%
\section*{Acknowledgements}
I would like to thank the referees for their valuable comments and suggestions.

\end{document}